# The Balkans Continued Fraction

David Naccache[1] and Ofer Yifrach-Stav[1]

DIÉNS, ÉNS, CNRS, PSL University, Paris, France
45 rue d'Ulm, 75230, Paris CEDEX 05, France
ofer.friedman@ens.fr, david.naccache@ens.fr

**Abstract.** In a previous escapade we gave a collection of continued fractions involving Catalan's constant. This paper provides more general formulae governing those continued fractions. Having distinguished different cases associated to regions in the plan, we nickname those continued fractions "The Balkans" as they divide into areas which are related but still different in nature.
Because we do not provide formal proofs of those machine-constructed formulae we do not claim them to be theorems. Still, each and every proposed formula was extensively tested numerically.

## 1 Introduction

In a previous escapade [13] we gave a collection of continued fractions involving Catalan's constant. This paper provides more general formulae governing those continued fractions. Having distinguished different cases associated to regions in the plan, we nickname those continued fractions "The Balkans" as they divide into areas which are related but still different in nature.

Because we do not provide formal proofs of those machine-constructed formulae we do not claim them to be theorems. Still, each and every proposed formula was extensively tested numerically.

All the programs included in this article are *self-contained*, i.e. any code snippet can be run independently of the others to fully illustrate the encoded formula. This renders the code longer but has the great advantage of allowing the reader to run and modify each snippet directly by just cutting and pasting it into Mathematica without requiring any other module[1]. The code was compacted for the sake of concision but loading it into Mathematica's editor re-indents it automatically.

The code in this paper is necessary. Because we do not provide proofs explaining the discovered structures, any slight LaTeX misprint would be impossible to fix. Hence readers can consider the code as an unambiguously tested version of the proposed formulae.

---

[1] Each snippet ends by a `ClearAll["Global`*"];` command whose purpose is to make Mathematica "forget" all passed history.

## 2 Notations

We denote by $n!!$ the semifactorial of, i.e. the product of all the integers from 1 up to $n$ having the same parity as $n$:

$$n!! = \prod_{k=0}^{\lceil \frac{n}{2} \rceil - 1} (n-2k) = n(n-2)(n-4)\cdots$$

Because in all the following we will only apply semifactorials to odd numbers, this can be simplified as:

$$n!! = \prod_{k=1}^{\frac{n+1}{2}} (2k-1) = n(n-2)(n-4)\cdots 3\cdot 1$$

We denote by Catalan's constant by $G = 0.91596559\ldots$ and let $C_n$ be the $n$-th Catalan number:

$$C_n = \frac{1}{n+1}\binom{2n}{n} = \frac{(2n)!}{(n+1)!\,n!} = \prod_{k=2}^{n} \frac{n+k}{k} \qquad \text{for } n \geq 0$$

The first Catalan numbers are:

$$1, 1, 2, 5, 14, 42, 132, 429, 1430, 4862, 16796, 58786, \ldots$$

## 3 The target

We define for odd $j$ and $\kappa, c \in \mathbb{N}$ the following quantity nicknamed "*The Balkans continued fraction*":

$$Q_{j,\kappa,c} = j(2-j+2\kappa) + \mathop{\mathcal{K}}_{n=1}^{\infty}\left(\frac{-2n(c+n)(j+n-1)(1-j+2\kappa+n)}{j(2-j+2\kappa)+(3+4\kappa)n+3n^2}\right)$$

The question asked is that of finding a general process allowing to compute $Q_{j,\kappa,c}$ *without resorting* to numerical simulations or integer relation algorithms. The reason for this is that while integer relation algorithms allow us to "magically" discover relations, they do not provide *general information* about the underlying structure of the constants found. Why? Because we're in it for the thrill of discovery, not just the magic of shortcuts.

### 3.1 The contribution

The main contribution of this paper is a collection of formulae computing $Q_{j,\kappa,c}$ not involving any infinite sums and without requiring any numerical simulation for positive $j, \kappa, c$ and odd $j$.

### 3.2 Why this formula in particular?

The Ramanjuan Machine Project [15,5], lists a few continued fractions involving $G$ detected in 2020. We do not know why the project did not resort to (rather basic) integer relation algorithms to discover more relations. We hence decided to play detective and unleashed LLL that found a few hundreds of similar continued fractions involving $G$ in a few intensive calculation days. All continued fractions were of the form:

$$\epsilon + \mathop{\mathrm{K}}_{n=1}^{\infty}\left(\frac{-2n(n+\tau)(n+\eta)(n+\mu)}{\epsilon+\delta n+3n^2}\right) = \frac{a_0}{a_1+a_2G} \quad \text{where} \quad a_0, a_1, a_2 \in \mathbb{Z}$$

We performed two natural tests on the coefficient vectors $(\delta, \epsilon, \tau, \eta, \mu)$: a PCA to determine if the coefficients can be expressed as linear combinations of less than 5 variables and a Hough transform to detect affine relations in the dataset.

PCA revealed that, when projected on $(\delta, \tau, \eta, \mu)$, nearly all data was governed by three linear dimensions[2]. We hence understood that we were facing a linear behavior in a large region (Balkans) plus some sporadic cases (see Remark 12). This was also confirmed by the Hough transform that detected several parallel plans in the 3D-space. We thus decided to focus our efforts on the main plans in the 3D-space and understand them.

### 3.3 How formulae were reverse-engineered

Now comes the real adventure. The process that allowed us to reverse-engineer the formulae given in this paper is interesting by its own right. A quick look at many examples of the three quantities $a_0, a_1, a_2$ forming the fractions:

$$Q_{j,\kappa,c} = \frac{a_0}{a_1+a_2G} \quad \text{where} \quad a_0, a_1, a_2 \in \mathbb{Z}$$

showed that the $a_i$s are products of small prime factors and a few large prime factors. This suggested that $a_i$s were initially[3] of some form:

$$a_i = \text{expression}(j,\kappa,c) = \frac{\prod_{i=0}^{u-1} \phi_i'(j,\kappa,c)}{\prod_{i=0}^{v-1} \phi_i(j,\kappa,c)}$$

where the $\phi_i$ are "smooth" functions[4] such as $(an+b)!$, $2^{an+b}$, $(an+b)!!$, $C_{an+b}$, Pochhammer symbols of linear combinations of the parameters $j, \kappa, c$ etc. and a few unknown "mixing" functions causing the appearance of the large prime factors, e.g. polynomials or recurrence relations.

Fortunately, integer relation algorithms allow us to collect many instances of such forms for diverse $j, \kappa, c$ values. Hence the problem at hand consists in

[2] To which we gave the names $j, \kappa, c$.

[3] i.e., before simplification intervenes.

[4] i.e. functions generating smooth outputs.

identifying which $\phi_i$s are compatible with the cancellations due to the division. If a given $\phi$ is present in the expression then it is reasonably assumed that when tried for many $j, \kappa, c$ the new expression:

$$\text{expression}_1(j,\kappa,c) = \frac{\text{expression}(j,\kappa,c)}{\phi(j,\kappa,c)}$$

or

$$\text{expression}_1(j,\kappa,c) = \text{expression}(j,\kappa,c) \cdot \phi(j,\kappa,c)$$

will feature less small factors and hence stand-out as an outlier.

The process can hence be repeated with proper backtracking until all the combinatorial $\phi_i$s were peeled-off. Then it remains to detect what the remaining "mixing" functions are which is done by monitoring the average growth rate of those surviving constants to emit hypotheses on the type of recurrence relations (or polynomials) at hand or resorting to a variety of integer sequence recognition tools to identify the hidden culprits.

– We started our exploration with the simplest case of Bosnia & Herzegovina where $a_2 = 0$. We chose Bosnia & Herzegovina as a launchpad because it is the simplest *finite* Balkan continued fraction. We (rightfully) hoped that analyzing it will give us insight about the Balkans' structure before hell breaks loose with infinite continued fractions spitting $G$s into the convergence values.
– Having reverse-engineered Bosnia & Herzegovina we moved to Croatia for which $a_2 = 0$ as well. Because Croatia is a more complex version of Bosnia & Herzegovina, the insight gained in Bosnia & Herzegovina proved very helpful to derive further characteristics of Croatia . We thus stress that the distinction between Bosnia & Herzegovina and Croatia is paedagogic rather than technical.
– The insight gained in Croatia guided our software to the formula for Montenegro, which is simpler than Kosovo and Serbia given that Montenegro corresponds to $j = 1$.
– Having inferred Montenegro we moved on to Kosovo whose symmetry[5] with Serbia was quickly noted.

This work demonstrates the interest of statistical classifiers such as Maximum Likelihood Estimation (MLE) and Support Vector Machines (SVM) in mathematical exploration.

As will be shown, this "gradient descent" method proved itself very well, although it required a few thousands of computation hours on a very powerful cluster.

We estimate that 80% of the discovery effort was done by the machine. The remaining 20% being human "piloting" that, we are convinced, is already at the reach of today's most powerful LLMs, such as Gemini or AlphaEvolve.

[5] This symmetry simply comes from the fact that the equation $j(2-2j+2\kappa) = x(2-2x+2\kappa)$ has the two solutions $x = j$ and $x = \kappa - j + 1$.

It appears much better to read the coming sections first to understand what prey we are stalking and then refer to Section 12 describing the hunting process.

## 4 Kosovo, Serbia, Croatia, Montenegro, Bosnia & Herzegovina

.

As we will see, we distinguish five cases that we call Kosovo, Serbia, Croatia, Montenegro and Bosnia & Herzegovina after the regions of interest in the $(j,\kappa)$-space shown in Figure 1. In each region $c$ runs over the integers.

| Area | Domain | $Q_{j,\kappa,c}$ **is in** |
|---|---|---|
| Croatia | $j \geq 2\kappa+5$ | $\mathbb{Q}$ |
| Bosnia & Herzegovina | $j = 2\kappa+3$ | $\mathbb{Q}$ |
| Serbia | $2\kappa+1 \geq j \geq \kappa+3$ | $\mathbb{R}$ |
| Kosovo | $\kappa+2 \geq j \geq 3$ | $\mathbb{R}$ |
| Montenegro | $j = 1$ | $\mathbb{R}$ |

**Table 1.** Areas. We list fractions involving $G$ as being in $\mathbb{R}$ although it is currently unknown if $G \in \mathbb{Q}$ (this is a major open question).

A first restriction of our study will be to focus on odd $j$ that produces $Q_{j,\kappa,c}$ values involving $G$. Note that:

$$Q_{j,c,\kappa} = \begin{cases} \dfrac{a_0}{a_1 + a_2 G} & \text{if } j \text{ is odd} \\ \dfrac{a_0}{a_1 + a_2 \log 2} & \text{if } j \text{ is even} \end{cases}$$

Note as well that $j,\kappa,c$ can be negative. In the examples given in the appendix we adopt the notation:

$$\frac{a_0}{a_1 + a_2 \cdot \text{constant}} = T(0) + \mathop{\mathcal{K}}_{n=1}^{\infty} \left( \frac{P(n)}{T(n)} \right)$$

This paper does not treat the even $j$ case (a follow-up paper dealing with those cases is underway). Negative $j,\kappa,c$ are of little interest as the Balkan's numerator is defined as $-2n(c+n)(j+n-1)(1-j+2\kappa+n)$ hence, if $c<0$ the term $(c+n)$ will hit zero for $n=-c$ and subsequently render the summation finite. The same happens for negative $j$ with the term $(j+n-1)$ and for negative $\kappa$ with the term $(1-j+2\kappa+n)$.

## 5 The Montenegro Conjecture

We start with the first $j,\kappa$ space called Montenegro. Montenegro corresponds to the case $j=1$. In other words:

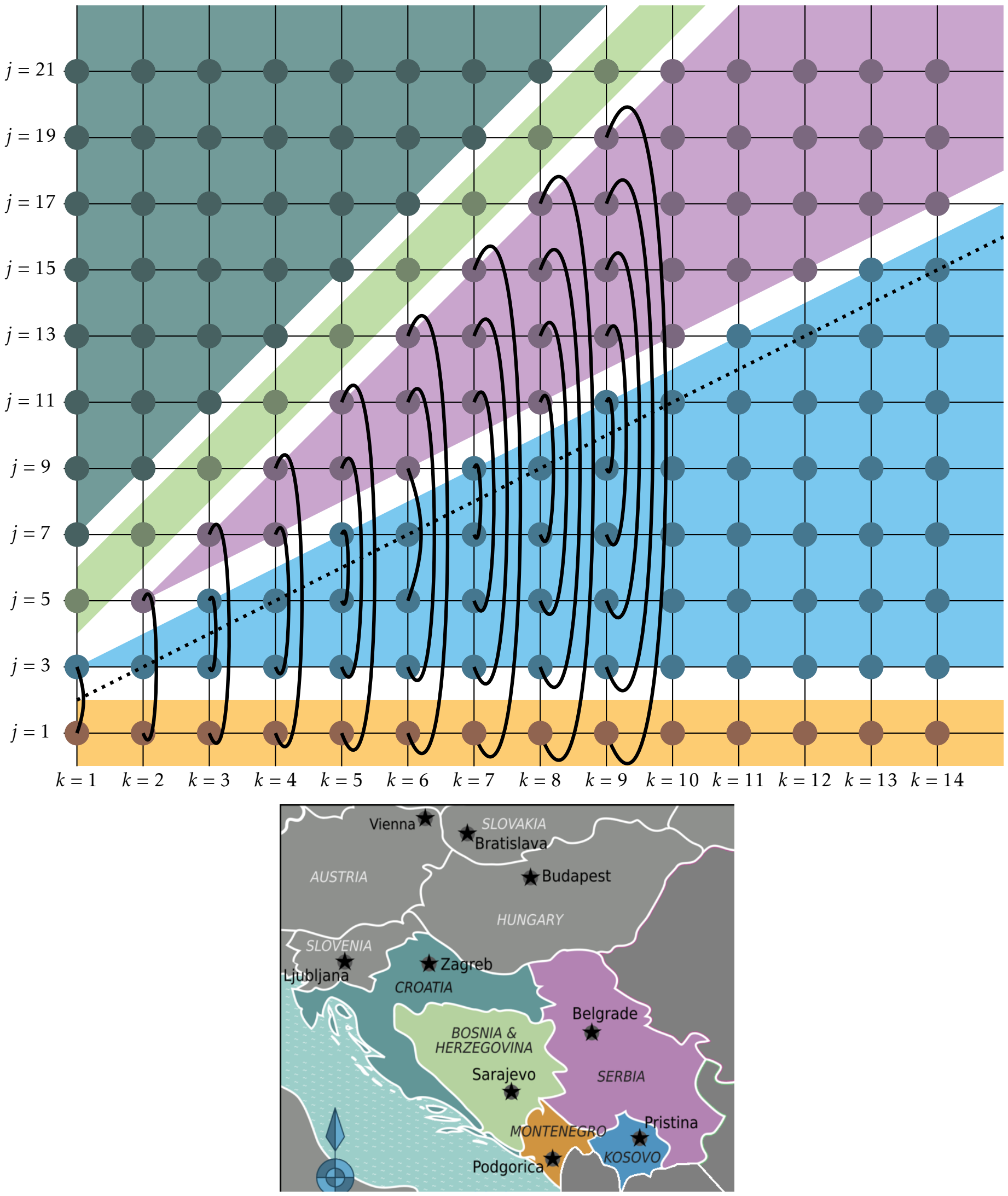

**Fig. 1.** The Five $j, k$ Areas. The meaning of the arches connecting Kosovo, Serbia, and Montenegro will be clarified later.

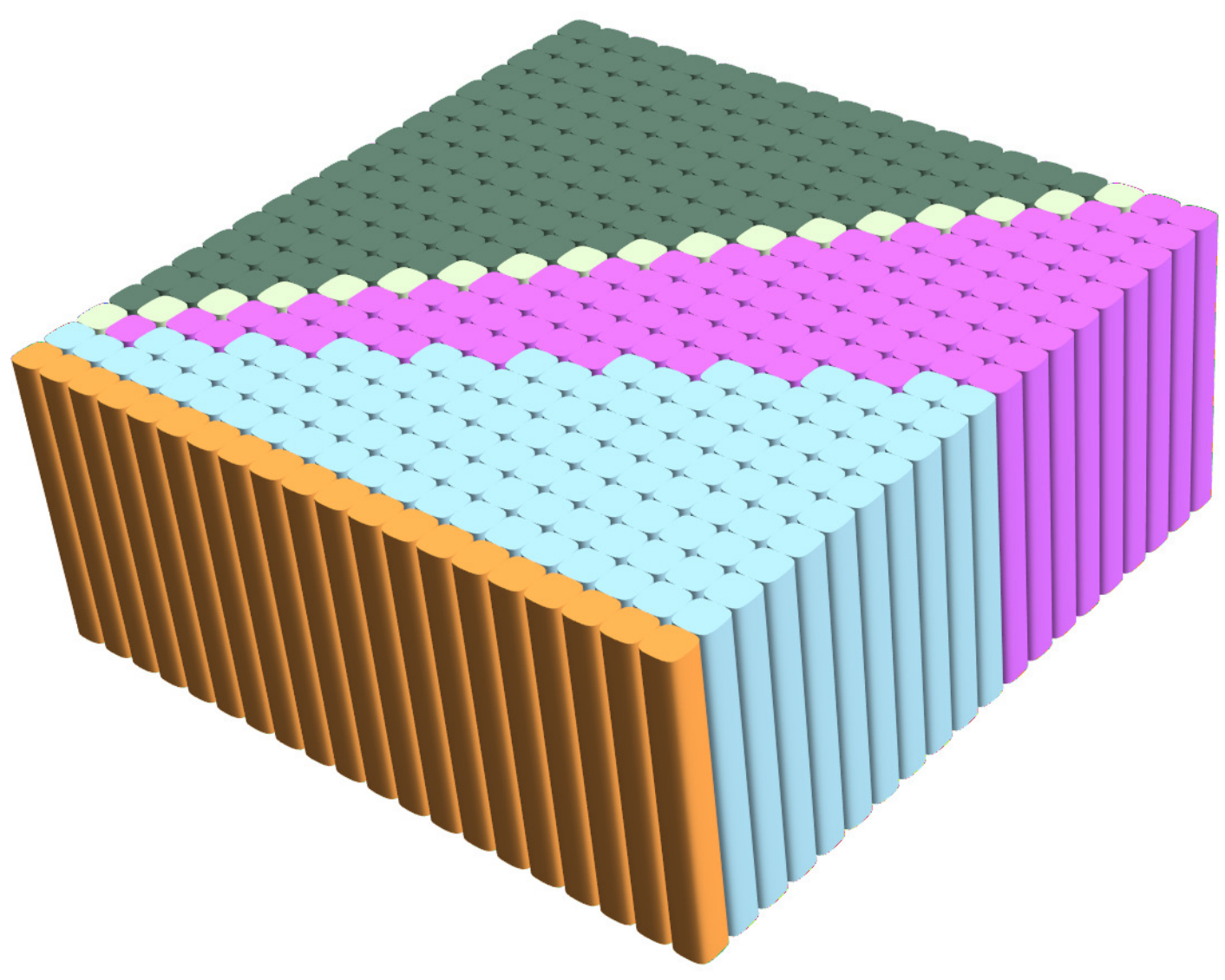

**Fig. 2.** The areas of Figure 1, with the dimension $c$ added. Each point in space corresponds to a $Q_{j,\kappa,c}$ value.

$$Q_{1,\kappa,c} = 2\kappa + 1 + \mathop{\mathcal{K}}_{n=1}^{\infty} \left( \frac{-2n^2(n+2\kappa)(n+c)}{3n^2 + (3+4\kappa)n + 2\kappa + 1} \right)$$

Define the functions:

$$\Delta_{\kappa,c}(\alpha,\beta) = \begin{cases} \alpha + \beta c & \text{if } c < 2 \\ -2c(2c-1)(2(c-\kappa)-1)^2 \Delta_{\kappa,c-2}(\alpha,\beta) & \text{if } c \geq 2 \\ \quad + (8c^2 + (2-8\kappa)c - 2\kappa + 1)\Delta_{\kappa,c-1}(\alpha,\beta) & \end{cases}$$

$$\Gamma_{\kappa,c}(\alpha,\beta) = (2c-1)!!^2 G + \Delta_{\kappa,c-1}(\alpha,\beta) \cdot \prod_{i=0}^{\kappa-1} (2(c-i)-1)$$

$$\delta_\kappa = \frac{4^{\kappa-1}}{(2\kappa-1)C_{\kappa-1}} \quad \text{and} \quad \rho_\kappa = \frac{\delta_\kappa (-1)^\kappa (1-2\kappa)}{(2\kappa)!(2\kappa-3)!!}$$

$$\alpha_\kappa = \rho_\kappa \Delta_{1,\kappa-1}(1,-2) \;\; \text{and} \;\; \beta_\kappa = -\rho_\kappa (2\kappa-3)^2 \Delta_{2,\kappa-1}(1,12) - \alpha_\kappa$$

$$\text{Then } \; \forall \kappa, c \in \mathbb{N}^2, \;\; Q_{1,\kappa,c} = \frac{\delta_\kappa (2c)!}{\Gamma_{\kappa,c}(\alpha_\kappa, \beta_\kappa)}$$

$Q_{1,\kappa,c}$ is hence an explicitly computable fraction of the form:

$$Q_{1,\kappa,c} = \frac{a_0}{a_1 + a_2 G} \quad \text{where} \quad a_0, a_1, a_2 \in \mathbb{Z}$$

The code snippet testing this formula over the square $1 \leq \kappa, c \leq 14$ is entitled `"1. Montenegro"`.

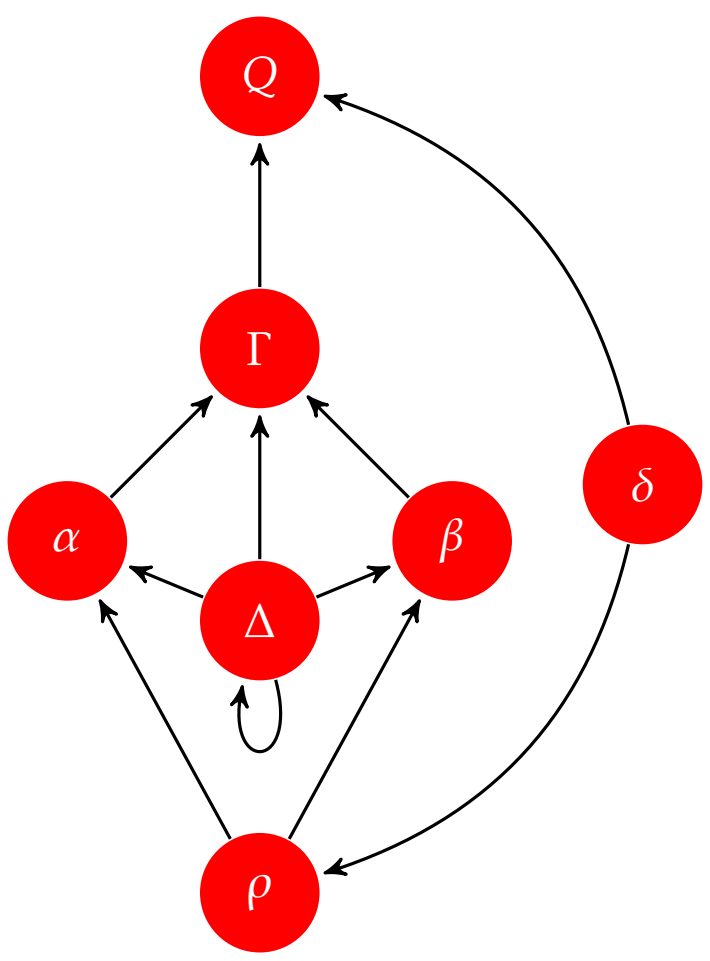

**Fig. 3.** Functional dependency between the functions computing $Q_{1,\kappa,c}$

In summery in Montenegro we do not need to resort to any integer relation algorithms nor to infinite summations to compute $Q_{1,\kappa,c}$ for all $\kappa, c$ values.

*Remark 1.* Note that, as described here, the complexity of $\Delta_{\kappa,c}$ is exponential in $c$, however, using classical Fibonacci memoization, this complexity can be reduced to $O(c \log c)$ thereby resulting in a very efficient algorithm for computing $Q_{1,\kappa,c}$.

## 6 The Bosnia & Herzegovina Conjecture

Bosnia & Herzegovina corresponds to the line $2\kappa = j - 3$, that is:

$$Q_{j,\frac{j-3}{2},c} = -j + \underset{n=1}{\overset{\infty}{\mathcal{K}}}\left(\frac{-2n(n-2)(n+c)(n+j-1)}{3n^2+(2j-3)n-j}\right) = 2+2c-j$$

We start by defining:

$$\alpha_j = 1 \quad \text{and} \quad \beta_j = 15 - 4j$$

And let:

$$\Delta_{j,c}(\alpha_j, \beta_j) = \begin{cases} \alpha_j + \beta_j c & \text{if } c < 2 \\ -2c(2c-j)(2c+1)(2c-j+2)\Delta_{j,c-2}(\alpha_j,\beta_j) & \text{if } c \geq 2 \\ \quad + (8c^2 + (14-4j)c - 3(j-2))\Delta_{j,c-1}(\alpha_j,\beta_j) & \end{cases}$$

$$g(j,c)=\frac{(2c)!}{2}\prod_{i=1}^{\frac{j-1}{2}}(2+2i-j) \text{ and } h(j,c)=\prod_{i=1}^{\frac{j-5}{2}}(2c-2i-1)$$

Then:

$$Q_{j,\frac{j-3}{2},c}=\frac{g(j,c)}{\Delta_{j,c-1}(\alpha_j,\beta_j)\cdot h(j,c)}\in\mathbb{Q}$$

Hence, in Bosnia & Herzegovina as well we do not need to resort to any integer relation algorithms to compute $Q_{j,\frac{j-3}{2},c}$ which is an explicitly computable fraction in $\mathbb{Q}$.

The corresponding code snippet is "`2. Bosnia`".

Note that $Q_{j,\frac{j-3}{2},c}=2+2c-j$ gives a non-recursive formula for $\Delta_{j,c}(1,15-4j)$.

This "long detour" for computing $Q_{j,(j-3)/2,c}$ which is a finite continued fraction (exactly computable by summation) is extremely useful as it unveiled the $\Delta$ structure and the basic $\phi_i$s that repeatedly intervene in other Balkan areas.

# 7 Roadmap

We are now ready to describe the roadmap that will govern the rest of this paper.

## 7.1 Areas governed by one running variable

As we have just seen, Montenegro and Bosnia & Herzegovina follow similar behaviors. This pattern revolves around the "magic" values $\alpha,\beta$ which are formally known for Montenegro and Bosnia & Herzegovina. Both regions are lines parameterized by a single running variable ($\kappa$ for Montenegro and $j$ for Bosnia & Herzegovina).

## 7.2 The $c$-level master formula for all Balkans except Montenegro

Except Montenegro (whose case was settled) all other areas obey a common $c$-level master formula that we will provide below. Computing $Q_{j,\kappa,c}$ for any $c$ using this $c$-level master formula requires knowledge of two rational parameters: $\alpha_{j,\kappa}$ and $\beta_{j,\kappa}$.

The cornerstone of the rest of this paper is thus the quest for $\alpha_{j,\kappa}$ and $\beta_{j,\kappa}$ for all Balkan areas except Montenegro.

$\alpha_{j,\kappa},\beta_{j,\kappa}$ can always be inferred by computing numerically[6] $Q_{j,\kappa,c_1}$ and $Q_{j,\kappa,c_2}$ for two values $c_1\neq c_2$ and solving a system of two equations in the unknowns $\alpha_{j,\kappa},\beta_{j,\kappa}$.

[6] i.e. using an integer relation algorithm.

When $\alpha_{j,\kappa}, \beta_{j,\kappa}$ are found $Q_{j,\kappa,c}$ can be computed for any other $c \notin \{c_1, c_2\}$.

Using this process requires resorting twice to integer relation algorithms such as LLL [12], HJLS [11], PSOS [1,7] or PSLQ [8] for each $(j,\kappa)$-pair. This works perfectly in practice but is not entirely satisfactory because the process has a "blind spot" which is the integer relation oracle. We would like to avoid such blind spots as much as possible and dispose of a fully algebraic process for computing each $Q_{j,\kappa,c}$. The ideal situation being, of course, a direct algebraic computation of $Q_{j,\kappa,c}$ from the data $j, \kappa, c$ alone.

We completely achieve the goal of computing all $Q_{j,\kappa,c}$ values in an exact algebraic way all over the Balkans. Note that we do not need to bother with Croatia where $Q_{j,\kappa,c}$ is a finite continued fraction computable exactly by summation. This simply stems from the fact that as we increase $n$ the $(1-j+2\kappa+n)$ component of Balkan's continued fraction will necessarily hit 0 in Bosnia & Herzegovina and in Croatia.

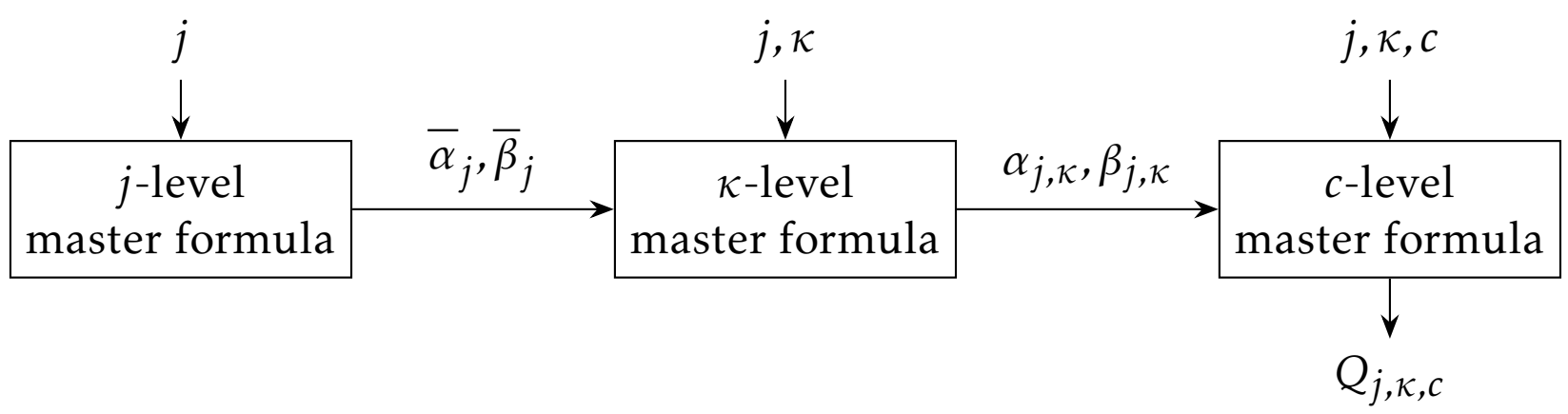

**Fig. 4.** $j$, $\kappa$ and $c$-level master formulae. This paper provides exact algebraic processes for inferring $Q_{j,\kappa,c}$ for all areas.

The snippet "`3. Northern Balkans`" illustrates the process with all areas except Montenegro. The snippet validates the $c$-level master formula on the 3D volume $3 \le j \le 13$, $1 \le \kappa \le 6$, $1 \le c \le 7$ (Figure 17) stretching over parts of Croatia, Bosnia & Herzegovina, Serbia and Kosovo. In each case, the program derives the corresponding $\alpha_{j,\kappa}, \beta_{j,\kappa}$ and allows the computation of $Q_{j,\kappa,c}$ for any $c$. The formally computed results are then successfully compared to numerical ones.

As we will later see, for Serbia and Kosovo we have more powerful master formulae operating at the $\kappa$-level and the $j$-level.

The $c$-level master formula[7] is defined as follows, assuming that we are somehow given the magic constants $\alpha_{j,\kappa}, \beta_{j,\kappa}$.

Let:

$$\Delta_{j,\kappa,c} = \begin{cases} \alpha_{j,\kappa} + \beta_{j,\kappa} c & \text{if } c < 2 \\ -2c(2c-j)(2c-2\kappa+j-2)(2c-2\kappa-1)\Delta_{j,\kappa,c-2} & \text{if } c \ge 2 \\ \quad + (8c^2 + (2-8\kappa)c + (j-2)(2\kappa-j))\Delta_{j,\kappa,c-1} & \end{cases}$$

[7] valid for all areas except Montenegro whose case was anyhow previously settled.

$$f_{j,\kappa,c} = C_{\frac{j-3}{2}} C_{\kappa-1} (j-2)(2\kappa-1)(2c-1)!!^2 \prod_{i=1}^{\frac{j-1}{2}} (2c-2\kappa+2i-1)(\kappa-i+1)$$

$$g_{j,\kappa,c} = (2c)! 2^{\frac{j+4\kappa-7}{2}} \prod_{i=1}^{\frac{j-1}{2}} (2c-2i+1)(2\kappa-2i+1)$$

$$h_{j,\kappa,c} = \prod_{i=0}^{\frac{j-3}{2}} (2c-2i-1) \prod_{i=0}^{\kappa-1} (2c-2i-1)$$

Then:

$$Q_{j,\kappa,c} = \frac{g(j,\kappa,c)}{\Delta_{j,\kappa,c-1} \cdot h(j,\kappa,c) + f(j,\kappa,c) \cdot G}$$

*Remark 2.* For Bosnia & Herzegovina and Croatia $f_{j,\kappa,c} = 0$. This happens automatically given the definition of $f_{j,\kappa,c}$.

## 8 The Serbia conjecture

Montenegro, Kosovo and Serbia feature internal symmetries illustrated by the arches in Figures 1 where connected $Q$ values are identical. We will also provide a formula connecting their $\alpha, \beta$ values over the plan. This means that any formula applicable to Kosovo will also settle the case of Serbia (given that Montenegro is settled and given that Montenegro provides the values for the upper border line of Serbia).

Given that:

$$Q_{1+2i,\kappa,c} = Q_{2\kappa+1-2i,\kappa,c} \quad \text{valid for} \quad \kappa = 2,3,\ldots \quad \text{and} \quad 0 \le i \le \left\lfloor \frac{k}{2} \right\rfloor - 1$$

to compute any value in Serbia, we can just compute its symmetric correspondance in Kosovo.

## 9 The Kosovo Conjecture

Kosovo is very specific in that it has a $\kappa$-level master formula. This means that for each $j$, knowledge of four upper-level rational constants[8]:

[8] This somewhat unusual notation is used to note the "$\alpha$ of $\alpha$", the "$\beta$ of $\alpha$" etc., as the same type of formula is applied at both $c$ and $\kappa$ levels.

$$\{\overline{\alpha}_j, \overline{\beta}_j\} = \left\{ \{\overset{\alpha}{\alpha}_j, \overset{\beta}{\alpha}_j\}, \{\overset{\alpha}{\beta}_j, \overset{\beta}{\beta}_j\} \right\}$$

allows to generate $\alpha_{j,\kappa}, \beta_{j,\kappa}$ for all $\kappa$ values along a $j$-line within Kosovo.

Finding $\overline{\alpha}_j, \overline{\beta}_j$ for a given $j$ requires solving a system of equations using four known $\alpha_{j,\kappa_1}, \beta_{j,\kappa_1}, \alpha_{j,\kappa_2}, \beta_{j,\kappa_2}$ for some $\kappa_1 \neq \kappa_2$. To determine the two $\alpha_{j,\kappa_i}, \beta_{j,\kappa_i}$ pairs[9] we can solve two systems of equations (each in two unknowns) using any $Q_{j,\kappa_1,c_1}, Q_{j,\kappa_1,c_2}, Q_{j,\kappa_2,c_3}$ and $Q_{j,\kappa_2,c_4}$. Note that there is no opposition to take $c_1 = c_3$ and $c_2 = c_4$. This process is illustrated in Figure 5.

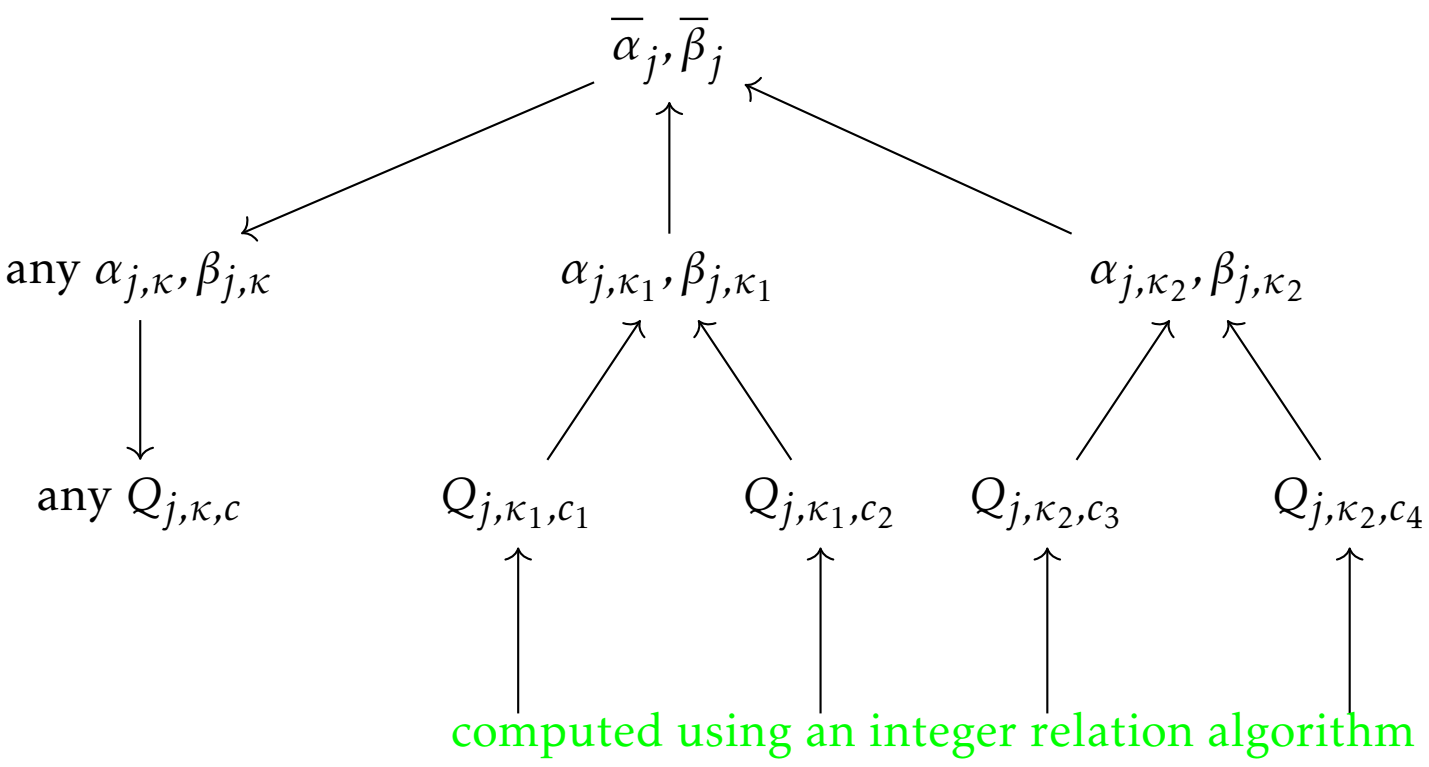

**Fig. 5.** $\kappa$-level resolution process for Kosovo.

### 9.1 The Kosovo $\kappa$-level master formula

The $\kappa$-level master formula for Kosovo allows to infer $\alpha_{j,\kappa}, \beta_{j,\kappa}$ from $\overline{\alpha}_j, \overline{\beta}_j$ for a fixed $j$ and a variable $\kappa$.

In other words, the Kosovo $\kappa$-level master formula generates for a fixed $j$ and for a variable $\kappa$ the data $\alpha_{j,\kappa}, \beta_{j,\kappa}$ necessary to operate the general $c$-level master formula given in subsection 7.2.

Recall that the master formula of the $c$ level of Subsection 7.2 generates for a fixed $j, \kappa$ and a variable $c$ a formal expression of the continued fraction $Q_{j,\kappa,c}$, given the auxiliary input $\alpha_{j,\kappa}, \beta_{j,\kappa}$.

Define:

$$s(j,\kappa) = \prod_{i=0}^{\frac{j-3}{2}} (\kappa - i)(2\kappa - 2i - 1)^2$$

[9] for $i \in \{1,2\}$.

$$\ell(n,j,\kappa) = \frac{(-1)^{\kappa+1}(2\kappa)!^2}{\kappa! 2^{3\kappa-2}(2\kappa-j)(2\kappa-1)(n((2\kappa-j-2)(3-2\kappa)-1)+1)\cdot s(j,\kappa)}$$

$$\eta(n,j,\kappa) = (2\kappa+2j-9-2n)(2\kappa+j-8-2n)(-2\kappa+5-j)(2\kappa+j-6)$$

$$\phi(n,j,\kappa) = 8\kappa^2+\kappa(10j-48-8n)+3j^2-(28+4n)j+68+18n$$

$$\bar{\Delta}_{n,j,\kappa}(\alpha,\beta) = \begin{cases} \alpha+\beta\kappa & \text{if } \kappa < 2 \\ \eta(n,j,\kappa)\cdot\bar{\Delta}_{n,j,\kappa-2}(\alpha,\beta)+\phi(n,j,\kappa)\cdot\bar{\Delta}_{n,j,\kappa-1}(\alpha,\beta) & \text{if } \kappa \geq 2 \end{cases}$$

We assume that we are given the four constants:

$$\{\overline{\alpha}_j,\overline{\beta}_j\} = \left\{\{\overset{\alpha}{\alpha}_j,\overset{\beta}{\alpha}_j\},\{\overset{\alpha}{\beta}_j,\overset{\beta}{\beta}_j\}\right\} \in \mathbb{Q}^4$$

Then:

$$\alpha_{j,\kappa} = \frac{\bar{\Delta}_{0,j,\kappa-j+2}(\overset{\alpha}{\alpha}_j,\overset{\beta}{\alpha}_j)}{\ell(0,j,\kappa)} \text{ and } \beta_{j,\kappa} = \frac{\bar{\Delta}_{1,j,\kappa-j+2}(\overset{\alpha}{\beta}_j,\overset{\beta}{\beta}_j)}{\ell(1,j,\kappa)} - \alpha_{j,\kappa}$$

The process is illustrated by the code snippet `"4. Kosovo"`.

The computation of the "magic" lists of constants $\{\overline{\alpha}_j,\overline{\beta}_j\}$ hard-coded in the snippet `"4. Kosovo"` is done by resorting to integer relation resolution in snippet `"5. resolution"`. A formal computation of $\overline{\alpha}_{j+2},\overline{\beta}_{j+2}$ from $\overline{\alpha}_j,\overline{\beta}_j$ is given in the next subsection.

*Remark 3.* It is interesting to note that thanks to the inner symmetry within Kosovo, it is possible to determine $\overline{\alpha}_{j+2},\overline{\beta}_{j+2}$ if $\overline{\alpha}_j,\overline{\beta}_j,\overline{\alpha}_{j+4},\overline{\beta}_{j+4}$ are known. Taking as an example $j=5$, Figure 1 shows that $Q_{5,3,c} = Q_{3,3,c}$. Hence knowledge of $\overline{\alpha}_3,\overline{\beta}_3$ will be used to compute $Q_{3,3,c}$ which is identical to $Q_{5,3,c}$.

Note that $Q_{5,5,c} = Q_{7,5,c}$. Hence knowledge of $\overline{\alpha}_7,\overline{\beta}_7$ will be used to compute $Q_{7,5,c}$ which is identical to $Q_{5,5,c}$. Finally, having in hand $Q_{5,3,c},Q_{5,5,c}$ we can solve a system in two unknowns and determine $\overline{\alpha}_5,\overline{\beta}_5$. Unfortunately this process has an information flow that only operates in "sandwich mode" allowing to determine $\overline{\alpha}_{j+2},\overline{\beta}_{j+2}$ from $\overline{\alpha}_j,\overline{\beta}_j,\overline{\alpha}_{j+4},\overline{\beta}_{j+4}$. This information flow cannot be reversed into an "escalator" allowing to ascend to level $j+4$ from levels $j$ and $j+2$. The situation is illustrated in Figure 6 where $x \longrightarrow y$ denotes the relation "$y$ is computable from $x$". This limitation is solved at the next subsection where we infer $\overline{\alpha}_{j+2},\overline{\beta}_{j+2}$ from $\overline{\alpha}_j,\overline{\beta}_j$.

### 9.2 The Kosovo $j$-level master formula

If $j = 3$ or $j = 5$ let:

$$\{\overset{\alpha}{\alpha}_3, \overset{\beta}{\alpha}_3, \overset{\alpha}{\beta}_3, \overset{\beta}{\beta}_3\} = \{-1, 4, -\frac{1}{3}, -\frac{14}{3}\} \text{ and } \{\overset{\alpha}{\alpha}_5, \overset{\beta}{\alpha}_5, \overset{\alpha}{\beta}_5, \overset{\beta}{\beta}_5\} = \{19, 234, -17, -8\}$$

If $j \geq 7$ define:

$$\varrho_j = \frac{2^{j+1}(j-1)!}{(j-2)(j-4)} \text{ and } b_j = (j-6)(j-2)(j-1)j(2j-7)(2j-5)$$

$$a_j = 4(j-1)j(2j-5) \text{ and } p_j = \frac{(j-6)(j-4)(j-1)(j+1)}{4}$$

$$d_j = 6j^6 - 15j^5 - 68j^4 + 74j^3 + 89j^2 - 44j - 18 \text{ and } e_j = (3j+1)(j^2-7) + 3$$

and iterate using the following formulae to compute $\{\overset{\alpha}{\alpha}_j, \overset{\beta}{\alpha}_j, \overset{\alpha}{\beta}_j, \overset{\beta}{\beta}_j\}$:

$$\overset{\alpha}{\alpha}_{j+2} = \frac{\overset{\alpha}{\alpha}_j a_j(2j-3) - (3j-2)\varrho_j}{(j-1)(j+1)} \text{ and } \overset{\beta}{\alpha}_{j+2} = \frac{\overset{\beta}{\alpha}_j a_j(2j+1) - e_{j-2}\varrho_j}{(j-3)(j+1)}$$

$$\overset{\alpha}{\beta}_{j+2} = \frac{\overset{\alpha}{\beta}_j b_j - 3\left((j-3)(j-1)-1\right)\varrho_j}{p_j} \text{ and } \overset{\beta}{\beta}_{j+2} = \frac{\overset{\beta}{\beta}_j b_j(j(2j-3)-1) - d_{j-2}\varrho_j}{p_j((j-2)(2j-7)-1)}$$

The test of those formulae proceeds in two steps. A first snippet ("`16. files`") computes the $\overline{\alpha}_j$ and $\overline{\beta}_j$ for $j = 3, 5, \ldots, 501$ using numerical simulation and records them in two files called `file1.txt` and `file2.txt`. Those files are read by snippet "`17. j-level`" that compares them to the values derived using the formal $j$-level formulae.

*Remark 4.* It is also possible to directly infer the $\alpha_{j,\kappa}, \beta_{j,\kappa}$ through symmetry. Note that the relation below also works for negative $j, \kappa$ values.

Define:

$$\zeta(j,u) = \frac{1}{2^u}\prod_{i=0}^{u-1}(2+2i-j)$$

$$\tau_{j,u} = \begin{cases} \dfrac{\zeta(j,u)}{|\zeta(j,u)|} \cdot (2j-2u-3)!!(2u-j)!!(-2)^u & \text{if } 2u > j-1 \\ \zeta(j,u) \cdot (-4)^u & \text{otherwise} \end{cases}$$

Then:

$$\frac{\alpha_{j,j-u-1}}{\alpha_{j-2u,j-u-1}} = \frac{\beta_{j,j-u-1}}{\beta_{j-2u,j-u-1}} = \tau_{j,u}$$

In particular $\tau_{j,j-1} = \frac{1}{2j-4}$. The code is snippet "`6. symmetry`".

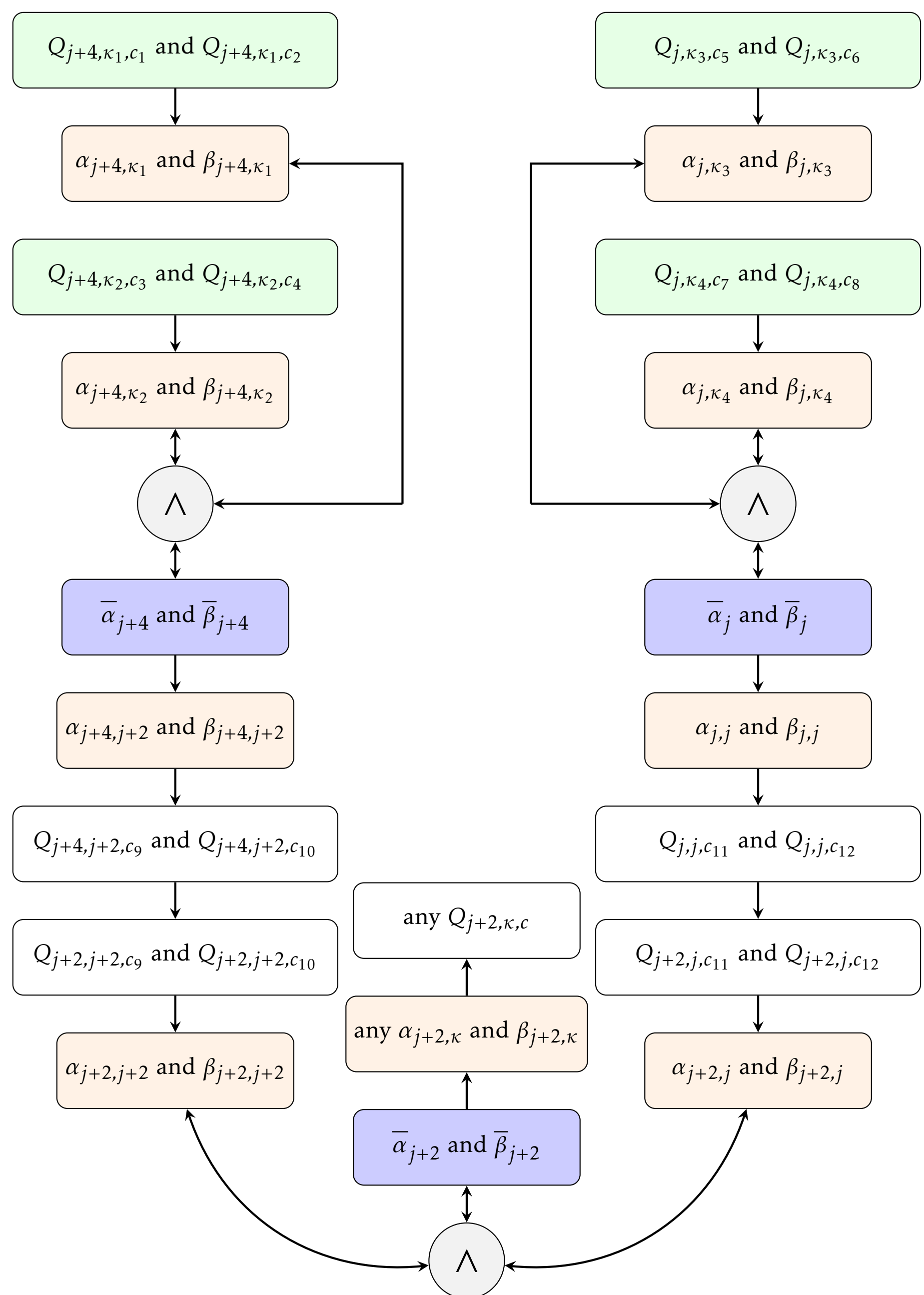

**Fig. 6.** Inferring $\overline{\alpha}_{j+2}, \overline{\beta}_{j+2}$ from $\overline{\alpha}_j, \overline{\beta}_j$ and $\overline{\alpha}_{j+4}, \overline{\beta}_{j+4}$. The $Q_{j,\kappa,c}$ in green boxes are determined using integer relation algorithms.

## 10 The Croatia Conjecture

As Bosnia & Herzegovina is a particular border case of Croatia the following is valid for both Bosnia & Herzegovina and Croatia.

Our automated software detected the following stunning behavior providing a $\kappa$-level formula for Croatia.

Let:

$$\mu_{i,j} = -(-2)^{\frac{3j-11-4i}{2}} \prod_{q=1}^{i} (j-2q-2)$$

For every $i = 0, 1, 2, \ldots$ there exist two polynomials in $j$, denoted $\psi_1(i,j)$ and $\psi_2(i,j)$ such that for $j \geq 2i+5$ we have:

$$\alpha_{j,\frac{j-2i-3}{2}} = \frac{\psi_1(i,j)}{\mu(i,j)} \quad \text{and} \quad \beta_{j,\frac{j-2i-3}{2}} = \frac{\psi_2(i,j)}{\mu(i,j)}$$

Tables 2 and 3 provide the first values of the polynomials $\psi_1, \psi_2$.

$\psi_1, \psi_2$ can be computed algebraically because when $(j, \kappa, c) \in$ Croatia, $Q_{j,\kappa,c}$ can be obtained by finite summation. We can hence derive $\alpha_{j,\frac{j-2i-3}{2}}$ and $\beta_{j,\frac{j-2i-3}{2}}$.

Then, by deriving enough $(\alpha_{j,\bullet}, \beta_{j,\bullet})$ pairs and knowing that $\deg_j \psi_1(i,j) = i$ and $\deg_j \psi_2(i,j) = i+1$, we can compute $\psi_1(i,j)$ and $\psi_2(i,j)$ by interpolation. We did not code this tedious yet straightforward process.

See code snippet `"7. Croatia"`.

| $i$ | $\psi_1(i,j)$ |
|---|---|
| 0 | $-1$ |
| 1 | $14-j$ |
| 2 | $-464+58j-3j^2$ |
| 3 | $27936-4692j+432j^2-15j^3$ |
| 4 | $-2659968+542256j-67836j^2+4260j^3-105j^4$ |
| 5 | $367568640-86278560j+13203480j^2-1139700j^3+51450j^4-945j^5$ |

**Table 2.** $\psi_1(i,j)$ for $0 \leq i \leq 5$.

*Remark 5.*

The leading coefficients of $\psi_1(i,j)$ (that is, $1, 1, 3, 15, 105, 945, \ldots$) are $(2i-1)!!$ whereas the leading coefficients of $\psi_2(i,j)$ (i.e. $4, 4, 12, 60, 420, 3780, \ldots$) are $4(2i-1)!!$.

*Remark 6.* The $\psi$ polynomials can always be written under a nested form, e.g.:

$$\begin{aligned}\psi_1(j,6) = &-14487726825-(104826150+(452605725\\ &+(121200300+(13697775+(640710\\ &+10395\cdot(j-27))(j-25))(j-23))(j-21))(j-19))(j-17)\end{aligned}$$

| $i$ | $\psi_2(i,j)$ |
|---|---|
| 0 | $-15+4j$ |
| 1 | $306-95j+4j^2$ |
| 2 | $-13360+4646j-357j^2+12j^3$ |
| 3 | $999648-379692j+40368j^2-2457j^3+60j^4$ |
| 4 | $-113885568+46449360j-6124164j^2+513228j^3-22935j^4+420j^5$ |
| 5 | $18333538560-7933530720j+1224286440j^2-126833100j^3+7864950j^4-266175j^5+3780j^6$ |

**Table 3.** $\psi_2(i,j)$ for $0 \le i \le 5$.

$$\begin{aligned}\psi_2(j,6) =& 3198013886925+(145296572850+(5207427225\\ &+(4353102000+(877052475+(78210090+(3023055\\ &+41580\cdot(j-29))(j-27))(j-25))(j-23))(j-21))(j-19))(j-17)\end{aligned}$$

*Remark 7.* The GCD between the $u$-th coefficient of $\psi_1(i,j)$ and the $u$-th coefficient of $\psi_2(i,j)$ is always smooth as illustrated in the code snippet "`8. coefficients`".

## 11 Balkans Knowledge Summary

In summary:

- Any $Q_{j,\kappa,c}$ value in Bosnia & Herzegovina is algebraically computable either by summation or by the Bosnian $\kappa$-level formula.
- In Croatia $Q_{j,\kappa,c}$ is directly computable by summation (hence making a $j$-level formula superfluous) and, in addition, has $\kappa$-level formulae.
- Any $Q_{j,\kappa,c}$ value in Montenegro is algebraically computable.
- For Kosovo we have $j$-level formulae.
- Serbia is fully determined by our knowledge of Montenegro and Kosovo.
- Negative $j,\kappa,c$ values are not of interest as they provide finite summations.

We therefore have algebraic formulae for computing all $Q_{j,\kappa,c}$ values over all the Balkans.

| **Area** | **$j$-level formula** | **$\kappa$-level formula** | **$c$-level formula** |
|---|---|---|---|
| Croatia | not required | ✓ | ✓ |
| Bosnia & Herzegovina | not required | ✓ | ✓ |
| Kosovo+Serbia | ✓ | ✓ | ✓ |
| Montenegro | not required | ✓ | ✓ |

**Table 4.** Knowledge Summary

A challenge on which the authors are currently working is characterizing the case of even values $j$ (which produce formulae involving $\log 2$).

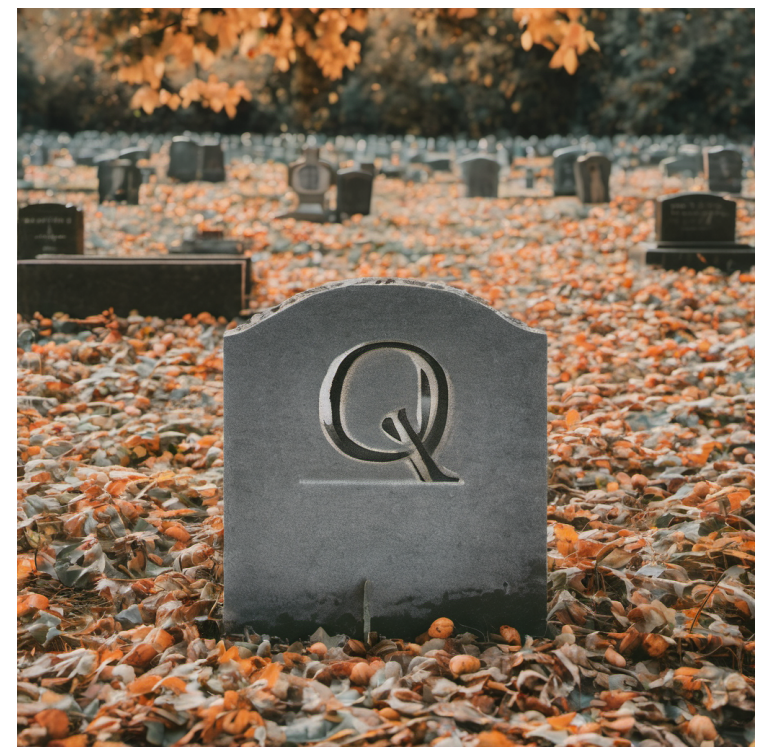

**Fig. 7.** $Q_{j,\kappa,c}$ *Requiescat in pace.*

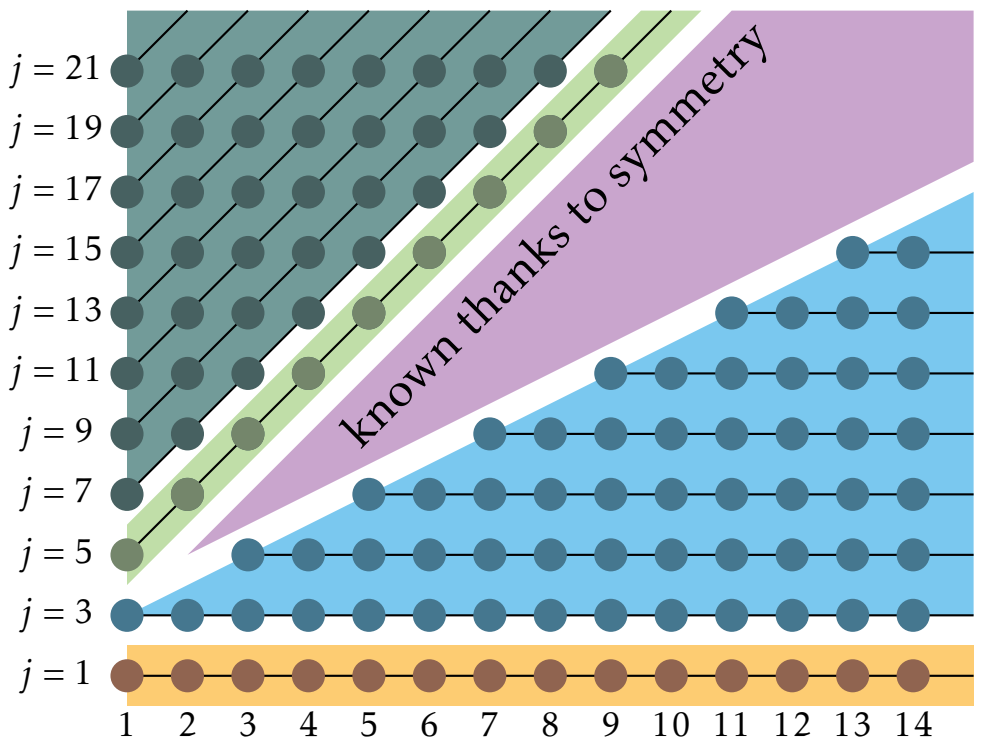

**Fig. 8.** The axes along which the $\kappa$-level master formulae operate in each area. We hence see that a finite amount of "magic" information in one dimension (that we fully provide in this paper) allows to algebraically compute $Q_{j,\kappa,c}$ over the two remaining dimensions.

As a motivational example let us unveil the simple $\log 2$ example $Q_{2,\kappa,0}$. In this case:

$$Q_{2,\kappa,0} = \frac{(-1)^{\kappa+1} a_\kappa}{-b_\kappa + a_\kappa \log 2}$$

Where $a_\kappa$ is the LCM of the list of $\kappa$ integers starting with $\kappa$ (in Mathematica: `Table[Apply[LCM,Table[i,{i,k,2k-1}]],{k,1,100]}`) while the $b_\kappa$ are the numerators of the coefficients in the power series for $-\log(1+x)\log(1-x)$.

Recall that the coefficients in the expansion of $\log(1+x)\log(1-x)$ are given by

$$\frac{1}{2}\binom{2n}{n}\int_{x=0}^{1} (x(1-x))^{n-1} \log(x)\,\mathrm{d}x$$

and

$$\log(1+x)\log(1-x) = \frac{1}{2}\int_{z=0}^{1} \frac{\log(z)}{z(1-z)}\left(\frac{1}{\sqrt{1-4x^2 z(1-z)}} - 1\right)\mathrm{d}z$$

Evidently, providing formal proofs of the formulae provided in this paper is a challenge by its own right. The Conservative Matrix Field approach of [6] is one promising direction to investigate. To date, attempts to code automated substitution-simplification-induction proofs using variants of Gosper's algorithm [10] were unsuccessful. Yet another interesting question is that of *reversal*: Given a $Q_{j,\kappa,c}$ (i.e. $a_0, a_1, a_2$) find $j, \kappa, c$[10].

The main open questions remaining are very simple to formulate:

**Open Question 1**

Prove the formulae given in this paper.

**Open Question 2**

What happens in the Balkans for even $j$ values?

**Open Question 3**

What happens in Inostranstvo (cf. Remark 12)?

The following observations are hints that may serve in future quests:

*Remark 8.*

$$Q_{1,0,c} = \frac{(2c)!}{2(2c-1)!!^2 G - \Delta_{c-1,0}}$$

[10] Reversal is not always possible over Montenegro because $\forall \kappa,\ Q_{1,\kappa,1} = Q_{1,1,\kappa}$.

where:

$$\Delta_{1,0,c} = \begin{cases} 1+10c & \text{if } c<2 \\ 2c(1-2c)^3\Delta_{1,0,c-2} + (8c^2+2c+1)\Delta_{1,0,c-1} & \text{if } c \geq 2 \end{cases}$$

or under an equivalent more compact form:

$$\Delta_{1,0,c} = \begin{cases} 1 & \text{if } c=0 \\ (2c)! + (2c+1)^2\Delta_{1,0,c-1} & \text{if } c>0 \end{cases}$$

and even [4]:

$$\Delta_{1,0,c-1} = (2c)!\left(\frac{2G\binom{2c}{c}}{4^c} - \int_0^\infty \frac{t}{\cosh^{2c+1}(t)}\,\mathrm{d}t\right)$$

*Remark 9.*

$$\lim_{c\to\infty} Q_{j,\kappa,c+1} - Q_{j,\kappa,c} = 2$$

$$\lim_{\kappa\to\infty} Q_{j,\kappa+1,c} - Q_{j,\kappa,c} = 2j$$

$$\lim_{\kappa\to\infty} Q_{j+1,j+2r+1,c} - Q_{j,j+2r+1,c} = 4r+1$$

*Remark 10.* Denoting:

$$Q_{j,\kappa,c} = \frac{a_0}{a_1 + a_2 G} \quad \text{where} \quad a_0, a_1, a_2 \in \mathbb{Z}$$

The following formula (code snippet "`9. ratio`") is valid all over areas:

$$\rho_{j,\kappa,c} = \prod_{i=1}^{\frac{j-1}{2}} \frac{(2c-2\kappa+2i-1)(\kappa-i+1)}{(2c-2i+1)(2\kappa-2i+1)} \quad \text{and} \quad \varepsilon_{j,\kappa} = 2\kappa + \frac{j-7}{2} + \left\lfloor \frac{1}{j} \right\rfloor$$

$$\frac{a_0}{a_2} = \frac{(2c)! \cdot 2^{\varepsilon_{j,\kappa}}}{(2c-1)!!^2 \cdot C_{\kappa-1} \cdot C_{\frac{j-3}{2}} \cdot (2\kappa-1)\cdot(j-2)\cdot\rho_{j,\kappa,c}}$$

Where $\rho_{1,\kappa,c} = 1$ by definition.

*Remark 11.* Although possibly unrelated, we note that low-degree continued fractions involving $\log 2$ can be also obtained with lower degree polynomials, e.g. (See code snippet "`10. log2-a`".):

$$\frac{2}{L(\frac{1}{2},1,c-1)} = \frac{2}{\sum_{n=0}^{\infty} \frac{e^{\pi i n}}{(n+c-1)}} = \frac{1}{2^{c-2}\log(2) - \sum_{j=1}^{c-2} \frac{2^{c-j-2}}{j}} = c + \mathop{\mathcal{K}}_{n=1}^{\infty}\left(\frac{-2n^2}{3n+c}\right)$$

We get a similar behavior for:

$$R_c = c + \underset{n=1}{\overset{\infty}{\mathrm{K}}} \left( \frac{-2n^2 - 2n}{3n + c} \right)$$

where $2^{c-4} \cdot (c-3) \cdot a_0 = a_2$ and for which we provide numerical examples in Table 12.

See code snippet "`11. log2-b`".

*Remark 12.* We also discovered other continued fractions involving $G$ outside the Balkans. This suggests the existence of a more general formula encompassing both the Balkans and those other territories (called "Inostranstvo").

We denote those relations:

$$Q'_{\delta,\epsilon,\tau,\eta,\mu} = \epsilon + \underset{n=1}{\overset{\infty}{\mathrm{K}}} \left( \frac{-2n(n+\tau)(n+\eta)(n+\mu)}{\epsilon + \delta n + 3n^2} \right) = \frac{a_0}{a_1 + a_2 G} \quad \text{where} \quad a_0, a_1, a_2 \in \mathbb{Z}$$

Luckily, given that we have two coefficients[11] in the denominator of the continued fractions, we can compare in one plot the coefficients $\epsilon$ and $\delta$ of the Balkans and of Inostranstvo in one figure (Figure 9). Similarly we can visualize in 3D the $\tau, \eta, \mu$ of both regions (Figure 10). The alignments of red points show that there is clearly another structured family hiding out beyond the Balkans.

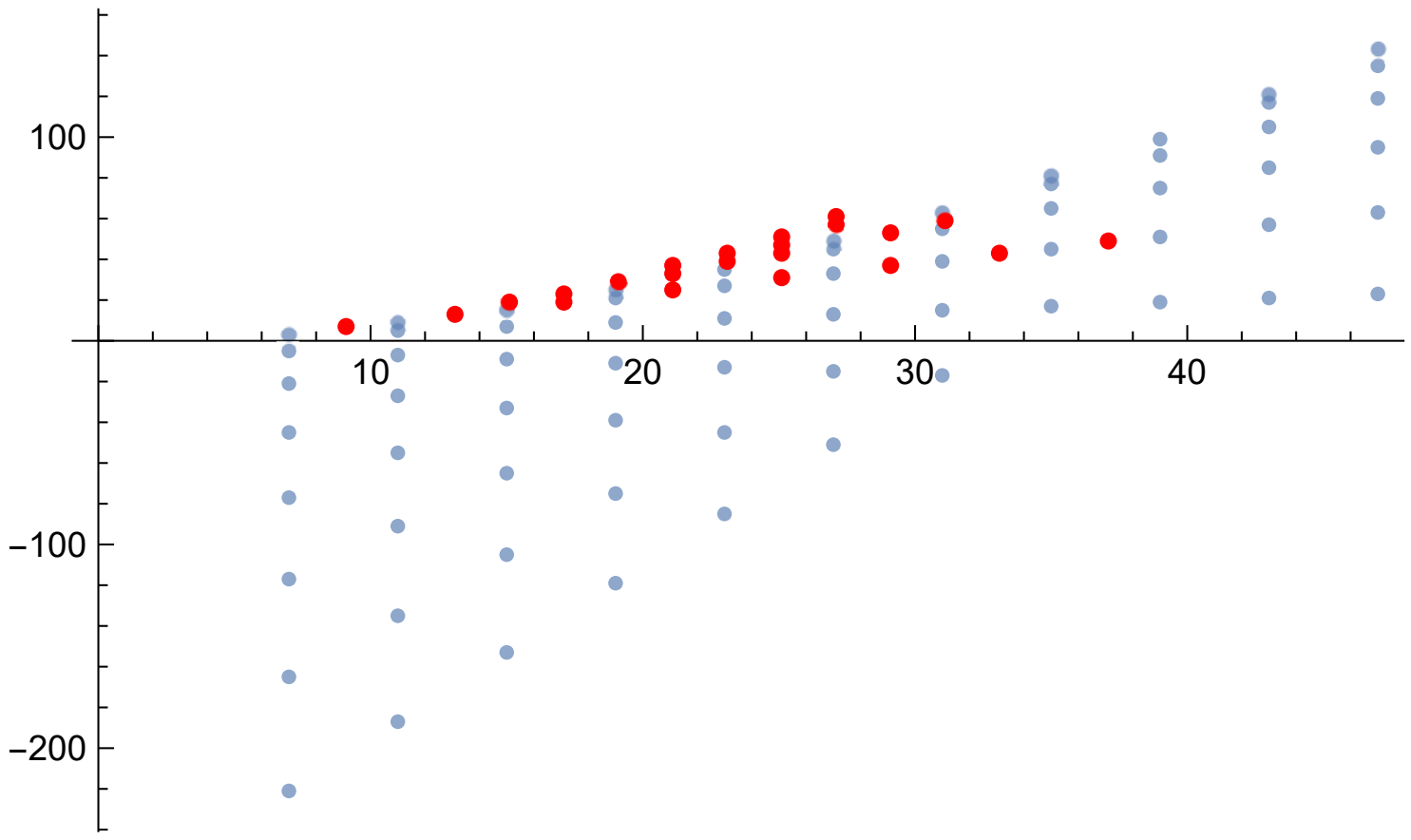

**Fig. 9.** $\epsilon, \delta$ for the Balkans (in blue) and Inostranstvo (in red).

PCA and automated matching revealed that $\epsilon, \delta$ are dependent on $\tau, \eta, \mu$ and:

[11] We exclude the 3 of $3n^2$ in the denominator, which is common to both the Balkans and to Inostranstvo.

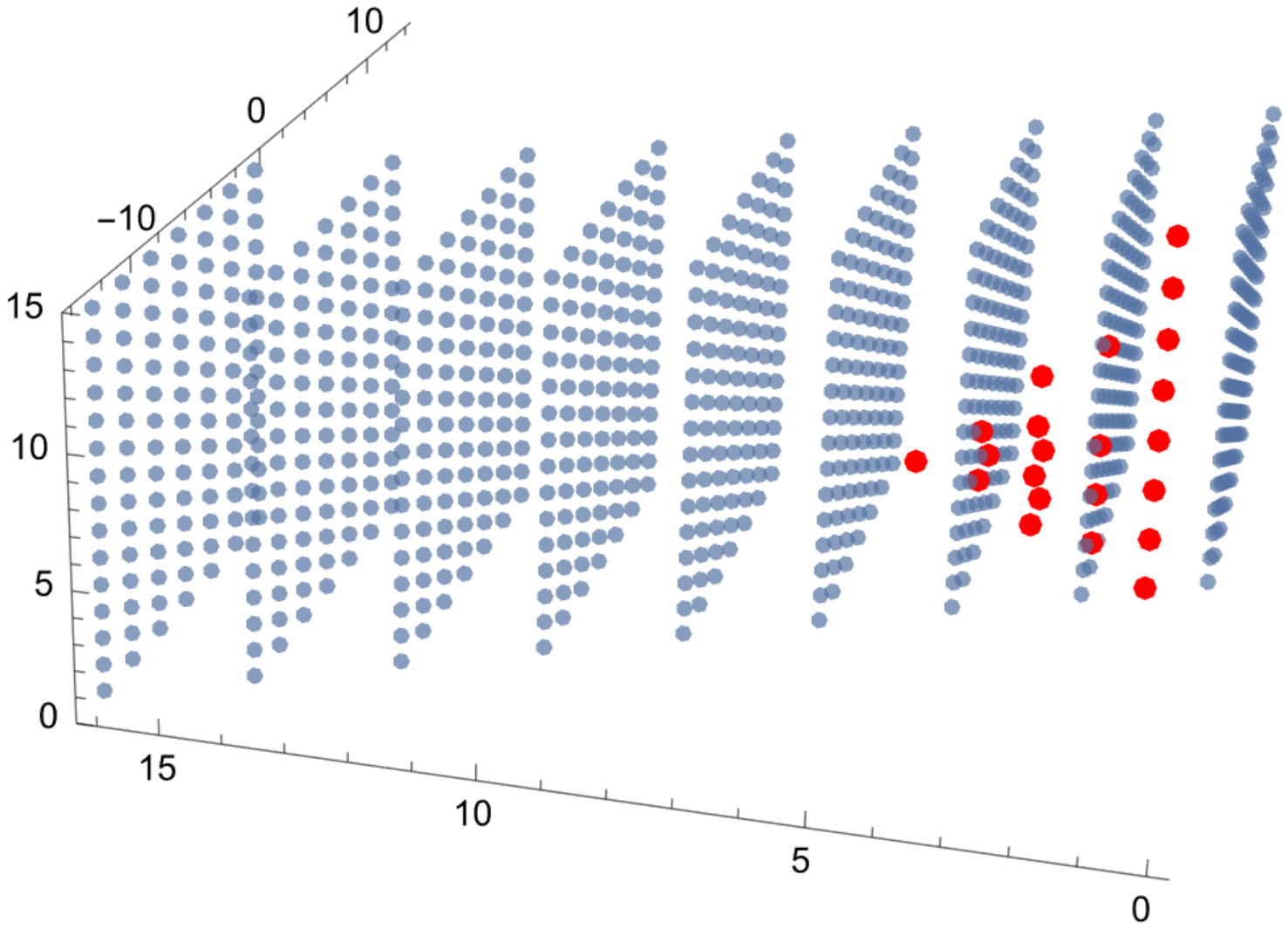

**Fig. 10.** $\tau, \eta, \mu$ for the Balkans (in blue) and Inostranstvo (in red).

$$Q'_{\tau,\eta,\mu} = x + 2(\tau+\eta+1)i + \underset{n=1}{\overset{\infty}{\mathcal{K}}}\left(\frac{-2n(n+\tau)(n+\eta)(n+2i+\mu)}{x+2(\tau+\eta+1)i+(2(\tau+\eta+\mu+2i)+3)n+3n^2}\right)$$

Where $x = (1+\eta)(1+\mu) + \tau(1+\eta+\mu)$.

This formula works when all variables[12] have identical parity, i.e.:

$$\tau \equiv \eta \equiv \mu \bmod 2$$

For instance:

$$Q'_i = 7 + 6i + \underset{n=1}{\overset{\infty}{\mathcal{K}}}\left(\frac{-2n(n+1)^2(n+2i+1)}{7+6i+(9+4i)n+3n^2}\right) \quad \text{for} \quad i = 0, 1 \ldots$$

Whose formal expression (code snippet "`12. Inostranstvo1`") turns out to be:

$$Q'_i = \frac{(2i+1)!}{\Delta'_i - 2G(2i+1)!!^2}$$

$$\Delta'_i = \begin{cases} 2+15i & \text{if } i < 2 \\ 2(2i-1)^3(1-i)\cdot\Delta'_{i-2} + (8i^2-2i+3)\cdot\Delta'_{i-1} & \text{if } i \geq 2 \end{cases}$$

[12] (code snippet "`14. Inostranstvo`")

or[13]:

$$Q_i'' = 23 + 10i + \mathop{\mathcal{K}}_{n=1}^{\infty}\left(\frac{-2n(n+1)(n+3)(n+2i+3)}{23+10i+(17+4i)n+3n^2}\right) \text{ for } i = 0, 1\dots$$

We did not investigate further the various Inostranstvo families but conjecture that they share the same behaviors as the Balkans.

*Remark 13.* The $\Delta$ functions appearing in this paper are particular cases of "Generalized Fibonacci Polynomials" (GFBs) studied by various authors, e.g. [9]. GFBs have numerous properties that might shed light on the open questions listed *supra*. We did not investigate this further.

Yet another route consists in considering the Inostranstvo and/or even-$j$ targets as algebraic equations and attempting on them an algebraic sieving approach such as [3]. This was not investigated this process given its theoretical and logistical complexity.

Finally, the similarity between the infinite sums given in [14] and the continued fractions investigated in this paper may reveal connections allowing to prove our conjectures. While analyzing [14] we noted a probable misprint in the formulae given for $y_3$ and $y_5$ (bottom of page 9 of [14]). We hence conducted our own experiments and discovered the relations described in Table 5 where:

$$w_2c + w_3 + w_1\sum_{n=1}^{\infty}(-1)^{n+1}\prod_{i=1}^{7}(2n+2i-3+\epsilon)^{-e_i} = 0$$

## 12 The Gradient Descent Process

The Gradient Descent process that generated the formulae will be presented in steps. We sample execution at critical points using the arbitrary example $j = 11, \kappa = 6$ and $40 \leq c \leq 47$ to explain the automated exploration process.

### 12.1 The starting point

We first recall our notations:

$$Q_{j,\kappa,c} = \frac{a_0}{a_1 + a_2 G} \text{ where } a_0, a_1, a_2 \in \mathbb{Z}$$

We have by now "seen the end of the movie", and we know that for all areas:

[13] for which

$$\frac{a_0}{a_2} = \frac{(2j+5)!}{(2j+4)(2j+5)!!^2}$$

(code snippet "`13. Inostranstvo2`")

| $\epsilon$ | $c$ | $e_1, e_2, \ldots, e_7$ | $w_1, w_2, w_3$ |
|---|---|---|---|
| 0 | $\pi$ | 0, 0, 2, 2, 2, 2, 0 | 57153600, −33075, 103904 |
| 0 | $\pi$ | 0, 1, 1, 2, 2, 1, 1 | 69854400, −24255, 76192 |
| 0 | $\pi$ | 0, 1, 2, 1, 1, 2, 1 | 62868960, −14553, 45712 |
| 0 | $\pi$ | 0, 2, 2, 1, 2, 2, 0 | −65318400, −4725, 14848 |
| 0 | $\pi$ | 0, 2, 2, 2, 2, 0, 0 | −6350400, −3675, 11552 |
| 0 | $\pi$ | 1, 1, 2, 2, 1, 1, 0 | −907200, −315, 992 |
| 0 | $\pi$ | 1, 2, 1, 1, 2, 1, 0 | −635040, −147, 464 |
| 0 | $\pi$ | 1, 2, 2, 2, 2, 1, 0 | 16934400, −245, 768 |
| 0 | $\pi$ | 2, 1, 1, 1, 1, 2, 0 | −59535000, −6615, 21292 |
| 0 | $\pi$ | 2, 2, 0, 0, 2, 1, 1 | −93139200, −2695, 9344 |
| 0 | $\pi$ | 2, 2, 0, 0, 2, 2, 0 | −52920000, −2695, 9056 |
| 0 | $\pi$ | 2, 2, 1, 2, 2, 0, 0 | −50803200, 3675, −11264 |
| 0 | $\pi$ | 2, 2, 2, 2, 0, 0, 0 | 129600, −75, 224 |
| 0 | $G$ | 0, 1, 2, 2, 2, 1, 0 | −50803200, 66150, −60577 |
| 0 | $G$ | 0, 3, 2, 3, 0, 0, 0 | −3456000, −6750, 6197 |
| 0 | $G$ | 1, 2, 2, 2, 1, 0, 0 | −2419200, −3150, 2909 |
| 0 | $G$ | 1, 3, 0, 3, 1, 0, 0 | 8064000, 8750, −8109 |
| 0 | $G$ | 2, 2, 0, 2, 2, 0, 0 | −33868800, −22050, 21131 |
| 0 | $G$ | 3, 2, 3, 0, 0, 0, 0 | −27648, 54, −25 |
| 1 | $\log 2$ | 0, 0, 0, 2, 2, 2, 2 | −33177600, −38400, 26617 |
| 1 | $\log 2$ | 0, 0, 2, 2, 2, 2, 0 | 1382400, −1600, 1109 |
| 1 | $\log 2$ | 0, 1, 1, 2, 2, 1, 1 | 11059200, −7680, 5323 |
| 1 | $\log 2$ | 0, 1, 2, 1, 1, 2, 1 | −22118400, 10240, −7097 |
| 1 | $\log 2$ | 0, 2, 2, 0, 0, 2, 2 | 17280000, −1760, 1219 |
| 1 | $\log 2$ | 0, 2, 2, 2, 2, 0, 0 | 442368, 512, −355 |
| 1 | $\log 2$ | 1, 1, 2, 0, 0, 2, 2 | −88473600, 5120, −3539 |
| 1 | $\log 2$ | 1, 1, 2, 2, 1, 1, 0 | −230400, −160, 111 |
| 1 | $\log 2$ | 1, 2, 1, 1, 2, 1, 0 | −22118400, −10240, 7109 |
| 1 | $\log 2$ | 2, 2, 0, 0, 2, 1, 1 | −88473600, −5120, 3627 |
| 1 | $\log 2$ | 2, 2, 0, 0, 2, 2, 0 | 69120000, 7040, −4951 |
| 1 | $\log 2$ | 2, 2, 1, 2, 2, 0, 0 | 7077888, −1024, 707 |
| 1 | $\log 2$ | 2, 2, 2, 2, 0, 0, 0 | −13824, 16, −11 |

**Table 5.** Relations for $\pi$, $G$ and $\log 2$. See code snippet "`15. Series`".

$$n_0(j,\kappa,c)=\frac{a_0}{a_2}=\frac{(2c)!\cdot 2^{\varepsilon_{j,\kappa}}}{(2c-1)!!^2\cdot C_{\kappa-1}\cdot C_{\frac{j-3}{2}}\cdot(2\kappa-1)\cdot(j-2)\cdot\rho_{j,\kappa,c}}$$

Where:

$$\rho_{j,\kappa,c}=\prod_{i=1}^{\frac{j-1}{2}}\frac{(2c-2\kappa+2i-1)(\kappa-i+1)}{(2c-2i+1)(2\kappa-2i+1)} \quad\text{and}\quad \varepsilon_{j,\kappa}=2\kappa+\frac{j-7}{2}+\left\lfloor\frac{1}{j}\right\rfloor$$

We started our journey by manually inspecting $n_0(j,\kappa,c)$ for several $j,\kappa,c$ values, noting that $n_0(j,\kappa,c)$ is always very smooth.

This suggests that $n_0(j,\kappa,c)$ is the product of basic "smooth" combinatorial functions such as factorials, binomials, semifactorials, Catalan numbers, Pochhammer symbols etc.

But the question is – of course – **which functions?**

We now know that $n_0(j,\kappa,c)$ is an exotic zoo containing the following animals:

$$\phi_0=(2c)!,\quad \phi_1=(2c-1)!!,\quad \phi_2=(2c-1)!!,\quad \phi_3=C_{\kappa-1},\quad \phi_4=C_{\frac{j-3}{2}}$$

$$\phi_5=(2\kappa-1),\quad \phi_6=(j-2),\quad \phi_7=2^{2\kappa},\quad \phi_8=2^{\frac{j-7}{2}},\quad \phi_9=2^{\lfloor\frac{1}{j}\rfloor}$$

$$\phi_{10}=\prod_{i=1}^{\frac{j-1}{2}}(2c-2\kappa+2i-1),\quad \phi_{11}=\prod_{i=1}^{\frac{j-1}{2}}(\kappa-i+1),\quad \phi_{12}=\prod_{i=1}^{\frac{j-1}{2}}(2c-2i+1)$$

$$\phi_{13}=\prod_{i=1}^{\frac{j-1}{2}}(2\kappa-2i+1)$$

$$n_0(j,\kappa,c)=\frac{\phi_0\cdot\phi_7\cdot\phi_8\cdot\phi_9\cdot\phi_{12}\cdot\phi_{13}}{\phi_1\cdot\phi_2\cdot\phi_3\cdot\phi_4\cdot\phi_5\cdot\phi_6\cdot\phi_{10}\cdot\phi_{11}}$$

However, at start, we had no idea what the $\phi_i$s were nor did we know how many $\phi_i$s are there.

*Remark 14.* The basic functions in our catalog are not independent as some multiplicatively generate others. e.g., Catalan numbers, binomials, multinomials and Pochhammer symbols are all products of factorials. Adding to the catalog $2^x$ we reach semifactorials etc. The code can hence *successfully* follow different paths for a given target $n_0(j,\kappa,c)$. While those functional dependencies do not impact the final result, they do impact complexity: e.g., using Catalan numbers *reduces the depth of search*[14] but *increases its width*.

Conversely, the code can also remove successfully components of $\phi_i$s and get stuck later if the remaining (un-removed part) of the $\phi_i$ is absent from the catalog.

[14] When a Catalan number is identified three factorial identifications are avoided.

## 12.2 Mutating functions

The algorithm performs a gradient descent on $n_i(j,\kappa,c)$, using an LLM to guide the descent. Catalog functions are not used in "bare metal" mode. They appear with specific linear combinations of $j,\kappa,c,1$. To capture those combinations let:

$$\sigma(\bar{u}) = u_0 j + u_1 \kappa + u_2 c + u_3 \quad \text{where} \quad \bar{u} = \{u_0, u_1, u_2, u_3\}$$

In other words[15]:

$$\phi_0 = \sigma(0,0,2,0)!, \quad \phi_1 = \phi_2 = \sigma(0,0,2,-1)!!, \quad \phi_3 = C_{\sigma(0,1,0,-1)}$$

$$\phi_4 = C_{\sigma(\frac{1}{2},0,0,-\frac{3}{2})}, \quad \phi_5 = \sigma(0,2,0,-1), \quad \phi_6 = \sigma(1,0,0,-2), \quad \phi_7 = 2^{\sigma(0,2,0,0)}$$

$$\phi_8 = 2^{\sigma(\frac{1}{2},0,0,-\frac{7}{2})}, \quad \phi_{10} = \prod_{i=1}^{\sigma(\frac{1}{2},0,0,-\frac{1}{2})} (\sigma(0,-2,2,-1) + 2i)$$

$$\phi_{11} = \prod_{i=1}^{\sigma(\frac{1}{2},0,0,-\frac{1}{2})} (\sigma(0,1,0,1) - i), \quad \phi_{12} = \prod_{i=1}^{\sigma(\frac{1}{2},0,0,-\frac{1}{2})} (\sigma(2,0,0,1) - 2i)$$

$$\phi_{13} = \prod_{i=1}^{\sigma(\frac{1}{2},0,0,-\frac{1}{2})} (\sigma(0,2,0,1) - 2i)$$

## 12.3 What information do we have?

The integer relations oracle[16] provides tens of thousands of $n_0(j,\kappa,c)$ instances. We record these in a database $B$ where each entry has the form:

$$B_i = \{j_i, \kappa_i, c_i, t_i\} = \{j_i, \kappa_i, c_i, n_0(j_i, \kappa_i, c_i)\}$$

.

The code needs to infer the $\bar{u}$s that intervene in $n_0(j,\kappa,c)$ from those numerous numerical examples.

## 12.4 A first step

Assume that the code is discovering the $\phi_i$s and their $\bar{u}$s one by one. The code is at some step $\omega$ at which it has already discovered the $\phi_i$s shown in red in the formulae:

$$n_0(j,\kappa,c) = \frac{a_0}{a_2} = \frac{(2c)! \cdot 2^{\varepsilon_{j,\kappa}}}{(2c-1)!!^2 \cdot C_{\kappa-1} \cdot C_{\frac{j-3}{2}} \cdot (2\kappa - 1) \cdot (j-2) \cdot \rho_{j,\kappa,c}}$$

---

[15] $\phi_9 = 2^{\lfloor \frac{1}{j} \rfloor}$ was not in the catalog and was added manually to unify two families of general formulae discovered by two runs.

[16] In our case LLL.

And:

$$\rho_{j,\kappa,c} = \frac{\prod_{i=1}^{\frac{j-1}{2}}(2c-2\kappa+2i-1)(\kappa-i+1)}{\prod_{i=1}^{\frac{j-1}{2}}(2c-2i+1)(2\kappa-2i+1)}$$

Divide the $t_i$ of each record $B_i$ by the red components evaluated at $j_i, \kappa_i, c_i$. Update the target to:

$$n_\omega(j,\kappa,c) = (2c-1)!!^2 \cdot C_{\kappa-1} \cdot C_{\frac{j-3}{2}} \cdot (2\kappa-1) \cdot (j-2) \cdot \prod_{i=1}^{\frac{j-1}{2}}(2c-2\kappa+2i-1)(\kappa-i+1)$$

A successfully mutated function, e.g. $\phi_1 = (2c-1)!! = \sigma(0,0,2,-1)!!$ stands-out because:

$$n_\omega(j,\kappa,c) \bmod \sigma(0,0,2,-1)!! = 0 \quad \text{for all} \quad j,\kappa,c \quad \text{values}$$

Evidently, this criterion generates false positives. For instance:

$$n_\omega(j,\kappa,c) \bmod \sigma(0,0,2,-1)!! = 0 \Rightarrow n_\omega(j,\kappa,c) \bmod \sigma(0,0,2,-3)!! = 0$$

*Remark 15.* False positives come in two flavors: “False false positives” and “True false positives”. $\sigma(0,0,2,-3)!!$ is a “false false positive”. Detecting $\sigma(0,0,2,-3)!!$ is useful as it decreases the target but allows progress. In the example above, instead of peeling-off $(2c-1)!!$ in one round, a first round will peel-off $(2c-3)!!$ and leave an extra $(2c-1)$ to some subsequent round. By opposition a “True false positive” is a candidate appearing as a factor of the target by the sole effect of chance over the available dataset.

To visualize efficiently execution we introduce *Backgammon diagrams*.

### 12.5 Backgammon diagrams

Our gradient descent monitoring tool is called “Backgammon diagrams” because of its visual similarity to a backgammon board (Figure 11).

In a Backgammon diagram the $x$-axis shows a search performed over a coordinate $u_s$ within an interval $u_s \in \{u_{\text{start}}, \ldots, u_{\text{end}}\}$.

In all the following, consider that all $\bar{u}$ coordinates other than $s$ were fixed to correct values denoted by ✓.

The draughts represent experiments with fixed $j,\kappa$ values and different $c$ values. In the example $j = 11, \kappa = 6$ and $40 \leq c \leq 47$. Each draught color corresponds to a different $c$ value (●,●,●,●,●,●,●,●). The legend is not repeated to save space. Draughts were slightly lifted to avoid covering each other.

If for a given $c$ value:

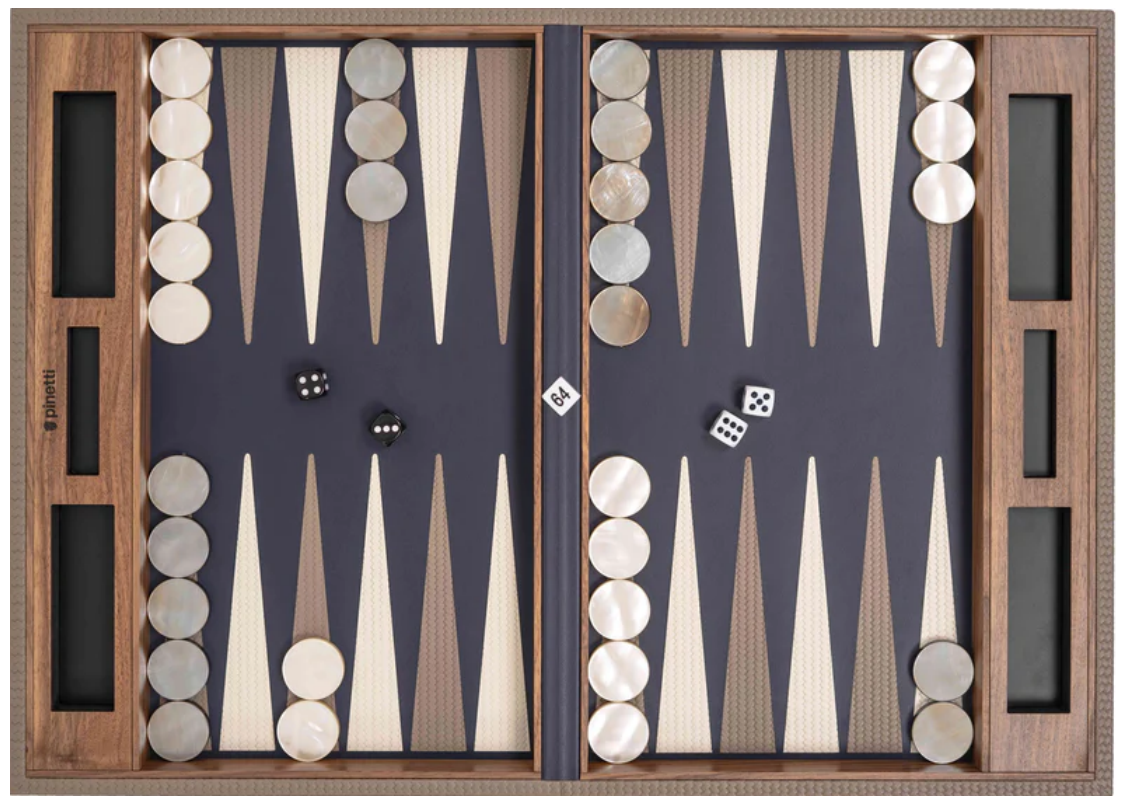

**Fig. 11.** A typical Backgammon board.

$$n_\omega(j,\kappa,c) \bmod \phi_i(\sigma(\checkmark,\ldots,\checkmark,u_s,\checkmark,\ldots,\checkmark)) = 0$$

then the draught of $c$'s color is lowered to the bottom of diagram, else it it raised to the top. Hence, a glance at the diagram shows which $u_s$ values are compatible with the target $n_\omega(j,\kappa,c)$ for each $c$.

Two additional features enhance reading: a red line showing the correct answer[17] and little red triangles (▲) denoting $u_s$ values for which tests succeeded <u>for all $c$ values</u>.

This view is grandly simplified with respect to reality. In the code $j,\kappa$ vary as well (resulting in multi-dimensional and hence unvisualizable backgammon diagrams) and several $u_s$s are simultaneously tried at each round.

### 12.6 A worked-out example

We want to disassemble:

$$n_\omega(j,\kappa,c) = (2c-1)!!^2 \cdot C_{\kappa-1} \cdot C_{\frac{j-3}{2}} \cdot (2\kappa-1) \cdot (j-2) \cdot \prod_{i=1}^{\frac{j-1}{2}} (2c-2\kappa+2i-1)(\kappa-i+1)$$

[17] This is, of course, the sought result – not provided by the software.

| **step** | $\omega+1$ |
|---|---|
| **target** | $n_{\omega}(j,\kappa,c)$ |
| **candidate** | $(uc-1)!!$ |
| | |
| **new target** | $n_{\omega+1}(j,\kappa,c) = n_{\omega}(j,\kappa,c)/(2c-1)!!$ |
| **database update** | $\forall i$ do $t_i = t_i/(2c_i-1)!!$ |
| **remarks** | We have one solution which is $u=2$. |

| **step** | illustrating option 1: $\omega+2$ |
|---|---|
| **target** | $n_{\omega+1}(j,\kappa,c)$ |
| **candidate** | illustrating option 1: $(uc-1)!!$ |
| | |
| **new target** | $n_{\omega+2}(j,\kappa,c) = n_{\omega+1}(j,\kappa,c)/(2c-1)!!$ |
| **database update** | $\forall i$ do $t_i = t_i/(2c_i-1)!!$ |
| **remarks** | We have one solution which is $u=2$. |

| **step** | illustrating option 2: $\omega+2$ |
|---|---|
| **target** | $n_{\omega+1}(j,\kappa,c)$ |
| **candidate** | illustrating option 2: $(2c+u)!!$ |
| | |
| **new target** | $n_{\omega+2}(j,\kappa,c) = n_{\omega+1}(j,\kappa,c)/(2c-1)!!$ |
| **database update** | $\forall i$ do $t_i = t_i/(2c_i-1)!!$ |
| **remarks** | 4 solutions added to backtracking list. |

| **step** | illustrating option 1: $\omega + 3$ |
| --- | --- |
| **target** | $n_{\omega+2}(j, \kappa, c)$ |
| **candidate** | illustrating option 1: $(uc - 1)!!$ |
| | −5 0 5 |
| **new target** | None, wrong guess. |
| **database update** | None, wrong guess. |
| **remarks** | No solutions: repeat step $\omega + 3$ with another candidate. |

| **step** | illustrating option 2: $\omega + 3$ |
| --- | --- |
| **target** | $n_{\omega+2}(j, \kappa, c)$ |
| **candidate** | illustrating option 2: $(u\kappa - 1)$ |
| | −5 0 5 |
| **new target** | $n_{\omega+3}(j, \kappa, c) = n_{\omega+2}(j, \kappa, c)/(2\kappa - 1)$ |
| **database update** | $\forall i$ do $t_i = t_i/(2\kappa_i - 1)$ |
| **remarks** | 7 solutions added to backtracking list. |

| **step** | illustrating option 3: $\omega + 3$ |
| --- | --- |
| **target** | $n_{\omega+2}(j, \kappa, c)$ |
| **candidate** | illustrating option 3: $(2\kappa + u)$ |
| | −5 0 5 |
| **new target** | $n_{\omega+3}(j, \kappa, c) = n_{\omega+2}(j, \kappa, c)/(2\kappa - 1)$ |
| **database update** | $\forall i$ do $t_i = t_i/(2\kappa_i - 1)$ |
| **remarks** | 14 solutions added to backtracking list. |

| **step** | illustrating option 1: $\omega+4$ |
|---|---|
| **target** | $n_{\omega+3}(j,\kappa,c)$ |
| **candidate** | illustrating option 1: $(j+u)$ |
| | |
| **new target** | $n_{\omega+4}(j,\kappa,c) = n_{\omega+3}(j,\kappa,c)/(j-2)$ |
| **database update** | $\forall i$ do $t_i = t_i/(j_i - 2)$ |
| **remarks** | 13 solutions added to backtracking list. |

| **step** | illustrating option 2: $\omega+4$ |
|---|---|
| **target** | $n_{\omega+3}(j,\kappa,c)$ |
| **candidate** | illustrating option 2: $(uj-2)$ |
| | |
| **new target** | $n_{\omega+4}(j,\kappa,c) = n_{\omega+3}(j,\kappa,c)/(j-2)$ |
| **database update** | $\forall i$ do $t_i = t_i/(j_i - 2)$ |
| **remarks** | 9 solutions added to backtracking list. |

| **step** | illustrating option 1: $\omega+5$ |
|---|---|
| **target** | $n_{\omega+4}(j,\kappa,c)$ |
| **candidate** | illustrating option 1: $C_{\kappa+u}$ |
| | |
| **new target** | $n_{\omega+5}(j,\kappa,c) = n_{\omega+4}(j,\kappa,c)/C_{\kappa-1}$ |
| **database update** | $\forall i$ do $t_i = t_i/C_{\kappa_i-1}$ |
| **remarks** | 7 solutions added to backtracking list. |

| **step** | illustrating option 2: $\omega+5$ |
|---|---|
| **target** | $n_{\omega+4}(j,\kappa,c)$ |
| **candidate** | illustrating option 2: $C_{u\kappa-1}$ |
| | −5 0 5 |
| **new target** | $n_{\omega+5}(j,\kappa,c) = n_{\omega+4}(j,\kappa,c)/C_{\kappa-1}$ |
| **database update** | $\forall i$ do $t_i = t_i/C_{\kappa_i-1}$ |
| **remarks** | 2 solutions added to backtracking list. |

| **step** | illustrating option 1: $\omega+6$ |
|---|---|
| **target** | $n_{\omega+5}(j,\kappa,c)$ |
| **candidate** | illustrating option 1: $\prod_{x=1}^{\frac{j-1}{2}}(2c+uk+2x-1)$ |
| | −5 0 5 |
| **new target** | $n_{\omega+6}(j,\kappa,c) = n_{\omega+5}(j,\kappa,c)/\prod_{x=1}^{\frac{j-1}{2}}(2c-2k+2x-1)$ |
| **database update** | $\forall i$ do $t_i = t_i/\prod_{x=1}^{\frac{j_i-1}{2}}(2c_i-2k_i+2x-1)$ |
| **remarks** | We have one solution which is $u=-2$. |

| step | illustrating option 2: $\omega+6$ |
|---|---|
| target | $n_{\omega+5}(j,\kappa,c)$ |
| candidate | illustrating option 2: $\prod_{x=1}^{\frac{j-1}{2}}(2c-2k+2x+u)$ |
| | −5 0 5 |
| new target | $n_{\omega+6}(j,\kappa,c)=n_{\omega+5}(j,\kappa,c)/\prod_{x=1}^{\frac{j-1}{2}}(2c-2k+2x-1)$ |
| database update | $\forall i$ do $t_i=t_i/\prod_{x=1}^{\frac{j_i-1}{2}}(2c_i-2k_i+2x-1)$ |
| remarks | We have one solution which is $u=-1$. |

| step | illustrating option 3: $\omega+6$ |
|---|---|
| target | $n_{\omega+5}(j,\kappa,c)$ |
| candidate | illustrating option 3: $\prod_{x=1}^{\frac{j-1}{2}}(uc-2k+2x-1)$ |
| | −5 0 5 |
| new target | $n_{\omega+6}(j,\kappa,c)=n_{\omega+5}(j,\kappa,c)/\prod_{x=1}^{\frac{j-1}{2}}(uc-2k+2x-1)$ |
| database update | $\forall i$ do $t_i=t_i/\prod_{x=1}^{\frac{j_i-1}{2}}(uc_i-2k_i+2x-1)$ |
| remarks | We have one solution which is $u=2$. |

| **step** | illustrating option 4: $\omega + 6$ |
|---|---|
| **target** | $n_{\omega+5}(j, \kappa, c)$ |
| **candidate** | illustrating option 4: $\prod_{x=1}^{\frac{2uj+j-1}{2}} (2c - 2\kappa + 2x - 1)$ |
| | |
| **new target** | $n_{\omega+6}(j, \kappa, c) = n_{\omega+5}(j, \kappa, c) / \prod_{x=1}^{\frac{2uj+j-1}{2}} (2c - 2\kappa + 2x - 1)$ |
| **database update** | $\forall i$ do $t_i = t_i / \prod_{x=1}^{\frac{2uj_i+j_i-1}{2}} (2c_i - 2\kappa_i + 2x - 1)$ |
| **remarks** | We have one solution which is $u = 0$. |

| **step** | $\omega + 7$ |
|---|---|
| **target** | $n_{\omega+6}(j, \kappa, c)$ |
| **candidate** | $C_{\frac{j+2u+1}{2}}$ |
| | |
| **new target** | $n_{\omega+7}(j, \kappa, c) = n_{\omega+6}(j, \kappa, c) / C_{\frac{j-3}{2}}$ |
| **database update** | $\forall i$ do $t_i = t_i / C_{\frac{j_i-3}{2}}$ |
| **remarks** | 3 solutions added to backtracking list. |

| **step** | $\omega + 8$ |
|---|---|
| **target** | $n_{\omega+7}(j, \kappa, c)$ |
| **candidate** | $\prod_{x=1}^{\frac{j-1}{2}} (\kappa - x + u)$ |
| | |
| **new target** | $n_{\omega+8}(j, \kappa, c) = n_{\omega+7}(j, \kappa, c) / \prod_{x=1}^{\frac{j-1}{2}} (\kappa - x + 1)$ |
| **database update** | $\forall i$ do $t_i = t_i / \prod_{x=1}^{\frac{j_i-1}{2}} (\kappa_i - x + 1)$ |
| **remarks** | 4 solutions added to backtracking list. |

| **step** | $\omega + 8$ |
|---|---|
| **target** | $n_{\omega+7}(j, \kappa, c)$ |
| **candidate** | $\prod_{x=1}^{\frac{j-1}{2}} (\kappa + ux + 1)$ |
| | |
| **new target** | $n_{\omega+8}(j, \kappa, c) = n_{\omega+7}(j, \kappa, c) / \prod_{x=1}^{\frac{j-1}{2}} (\kappa - x + 1)$ |
| **database update** | $\forall i$ do $t_i = t_i / \prod_{x=1}^{\frac{j_i-1}{2}} (\kappa_i - x + 1)$ |
| **remarks** | 2 solutions added to backtracking list. |

At this point one of the backtracking branches hits the constant function 1. $n_\omega(j, \kappa, c)$ was hence disassembled.

### 12.7 The Pathfinder

As we have just seen, at any step the algorithm can take several paths each of which offering a different backtracking fan-out. At first we considered resorting to a 2D-backtracking where one parameter is the offered fan-out and the second is the $\phi$-backtracking *per se*. The bookkeeping associated to this procedure seemed prohibitive, therefore we just opted to take the candidate with the least fan-out at each step. This appeared sufficient for our purpose and compatible with the computational means at hand.

### 12.8 The Decimator

A very significant speed-up is achieved by a software module called "the Decimator". The Decimator restricts the $\bar{u}$ space by removing $u_i$ affine combinations incompatible with the target.

Consider the target:

$$n_{\omega+2}(j,\kappa,c) = C_{\kappa-1} \cdot C_{\frac{j-3}{2}} \cdot (2\kappa-1) \cdot (j-2) \cdot \prod_{i=1}^{\frac{j-1}{2}} (2c-2\kappa+2i-1)(\kappa-i+1)$$

Create the database $B$ given in Table 6.

| $i$ | $j_i$ | $\kappa_i$ | $c_i$ | $t_i = n_{\omega+2}(j_i,\kappa_i,c_i)$ |
|---|---|---|---|---|
| 1 | 11 | 6 | 40 | 86562004597992000 |
| 2 | 11 | 6 | 41 | 99107222655672000 |
| 3 | 11 | 6 | 42 | 113065986409992000 |
| 4 | 11 | 6 | 43 | 128554477699032000 |
| 5 | 11 | 6 | 44 | 145695074725569600 |
| 6 | 11 | 6 | 45 | 164616513001617600 |
| 7 | 11 | 6 | 46 | 185454046292961600 |
| 8 | 11 | 6 | 47 | 208349607563697600 |

**Table 6.** The database $B$

Assume that we want to test a candidate $y(j,\kappa,c) = u_0 j + u_1\kappa + u_2 c + u_3$. We are typically interested in exploring each $u_i$ over a small interval, e.g. $[-8,8]$. Start by creating a list $\mathcal{L}$ of all $(2\times 8+1)^4$ $\vec{u}$-values.

For each $\{u_0,u_1,u_2,u_3\}$ value, if a $t_i$ is not divisible by at least one $y(j_i,\kappa_i,c_i)$ then remove $\{u_0,u_1,u_2,u_3\}$ from $\mathcal{L}$.

Table 7 gives the number of $\vec{u}$s eliminated from $[-8,8]^4$ by each $c$ value.

When we merge all the forbidden $\vec{u}$s, removing duplicates, we get a collection of 78567 combinations to skip. We are hence left with $(2\times 8+1)^4 - 78567 = 4954$ survivors to test, i.e. 5.9%.

| $c$ | 40 | 41 | 42 | 43 | 44 | 45 | 46 | 47 |
|---|---|---|---|---|---|---|---|---|
| **eliminated $\vec{u}$s** | 51721 | 54371 | 55635 | 56161 | 52771 | 51203 | 51609 | 46383 |

**Table 7.** Number of eliminated $\vec{u}$ combinations per $c$ tried.

The above example is restricted to eight $c$ toy-values. In reality we perform calculations on $\simeq 10^5$ combinations of $(j,\kappa,c)$. This reduces drastically the search space[18]. Note however that as $c$ increases the percentage of newly removed $\bar{u}$ values per round decreases.

Note that, in practice even an exploration in $[-3,3]^4$ would have sufficed[19]. For $\vec{u} \in [-3,3]^4$ we get 332 survivors out of 2401 (14%).

For $[-8,8]^4$ and $40 \leq c \leq 47$, decimating for $y(j,\kappa,c) = C_{u_0 j+u_1\kappa+u_2 c+u_3}$ removes 83233 candidates, leaving only 288 possibilities (0.35%).

Decimating over $[-5,5]^5$ and $20 \leq c \leq 60$ for the candidate:

$$y(j,\kappa,c) = \prod_{i=1}^{\frac{j-1}{2}} (u_0 j + u_1\kappa + u_2 c + u_3 + u_4 i)$$

results in 1496 survivors out of 161051 (0.93%).

Counter-intuitively, the more "complex"[20] the candidate is the more efficiently it is decimated. We hence select the candidates in decreasing complexity order. Measuring the complexity of mathematical formulas can be a subjective task, as it depends on various factors such as the number of terms, the presence of functions, exponents, and variables, as well as the overall structure. There's no definitive metric to quantify formula complexity. The following criteria can nonetheless be used to "measure" complexity and hence get rid of $\phi_i$s as fast as possible:

**Counting Elements**: We can count the number of distinct elements in each candidate, such as variables, constants, operators, and functions. For example, the candidate $\kappa - 1$ has two elements while in others we may have multiple variables, exponents, a product symbol, and a summation, which increases its complexity.

**Nesting and Hierarchy**: Analyze the nesting of operations and functions within the candidates. A candidate with multiple levels of nesting or hierarchy can be regarded as more complex.

**Mathematical Operations**: Consider the types of mathematical operations present in the candidates. More complex operations, such as exponentiation and summation, contribute to higher complexity compared to simpler operations like addition or multiplication.

[18] e.g., exploring for $20 \leq c \leq 60$ reduces the survivors' pool to 3702 (4.4%).

[19] The 7 in $\phi_8$ would have been decimated in two rounds.

[20] i.e., closer to a random oracle.

**Function Complexity**: If the candidate includes functions, their complexity should be taken into account. For instance, a Catalan number can add complexity compared to linear or constant functions.

**Symbolic Representation**: Represent each candidate in a symbolic format, such as a parse tree or abstract syntax tree. Compare the depth and branching of the trees as a rough measure of complexity.

**Information Theory**: Explore concepts from information theory, such as Kolmogorov complexity or algorithmic information theory, to quantify the amount of information needed to describe each candidate. This approach can be quite theoretical and may not provide a practical measure for all cases.

We resorted to a much more brutal approach [2]. Using the API[21], we asked GPT-4 to compare candidates pairwise, obtained a subjective complexity comparison ($\prec$) of each pair and translated the ternary results[22] to a directed graph. Because the LLM does not provide consistent answers (i.e. it might say that $A \prec B \prec C$ and... $C \prec A$) we performed a random walk of $10^6$ steps on the graph and counted the number of times each candidate was visited. The most visited candidate was considered as the "most complex". We then removed this candidate from the graph and started over again.

As a final note, we remark that some $\phi_i$s are "process killers". This is the case of candidates such as $(2c)!$. We called $\phi_i$s silencers "opioids" are they efficiently remove the symptoms but do not remove any lurgy. Because, in essence, $(2c)!$ contains just about any number of interest, it silences the modular tests. We hence start by launching the process with $\bar{u}!$ and start the backtracking afresh for each $n_0(j,\kappa,c)/\bar{u}!$. Luckily, the early investigation of Bosnia & Herzegovina by which we started allowed to uncover the opioids before formulae get too complex.

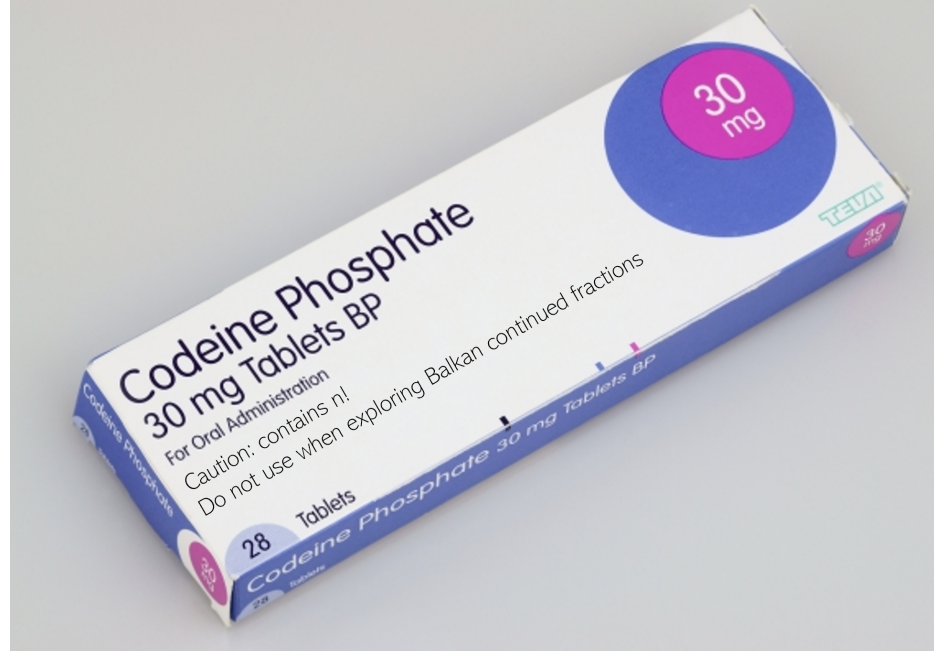

**Fig. 12.** *$Q_{j,\kappa,c}$ Process killers.*

[21] endpoint `https://api.openai.com/v1/engines/davinci-codex/completions`

[22] "$A$ is more complex that $B$": $A \to B$, "$B$ is more complex than $A$": $B \to A$ or "unsure": no edge.

### 12.9 Ascending to descend

As underlined in [16], some descents require intermediate ascents. The process described so far advances only in the case of a monotonous descent. In other words, if the initial $\phi_i$s were wrongly chosen the process will not converge.

We thus need a process allowing temporary ascents to get out of bowls. To that end we use a modified version of Broyden-Fletcher-Goldfarb-Shanno's (BFGS) algorithm. This requires a more refined measure of the penalty/profit of each move which cannot just be the Boolean "$x$ divides $y$".

Denote by $p_i$ the $i$-th prime.

We introduce a measure called "brittleness"[23] denoted $\Xi(n)$.

$$\text{Let: } n = \prod_{i=0}^{a-1} p_i^{m_i} \in \mathbb{Q} \text{ be a simplified fraction. Then } \Xi(n) = \sum_{i=0}^{a-1} |m_i|$$

In other words, $\Xi(n)$ counts, with repetition, the number of distinct factors appearing in either the numerator or the denominator of $n$.

The following example illustrates the evolution of brittleness ($y$-axis) during the peeling-off process.

Dotted lines $\Rightarrow$ Target's brittleness at the concerned $j, \kappa, c$ points.
Draughts $\Rightarrow$ Raised to $\Xi$(target/candidate) for different $u$ values.

As before, ▲ denotes $u$s at which all draughts are lower than their same-color dotted lines. The thin red vertical line is the correct answer.

As an example start with the target:

$$n_1(j,\kappa,c) = C_{\kappa-1} \cdot C_{\frac{j-3}{2}} \cdot (2\kappa-1) \cdot (j-2) \cdot \prod_{i=1}^{\frac{j-1}{2}} (2c - 2\kappa + 2i - 1)(\kappa - i + 1)$$

Figure 13 represents the following steps:

| **step** | 1 of Figure 13 |
|---|---|
| **target** | $n_0(j,\kappa,c)$ |
| **candidate** | illustrating first option: $u\kappa - 1$ |
| **new target** | $n_1(j,\kappa,c) = n_0(j,\kappa,c)/(2\kappa-1)$ |

| **step** | 1′ of Figure 13 |
|---|---|
| **target** | $n_0(j,\kappa,c)$ |
| **candidate** | illustrating second option: $2\kappa + u$ |
| **new target** | $n_1(j,\kappa,c) = n_0(j,\kappa,c)/(2\kappa-1)$ |

[23] Brittleness is a generalization the prime $\Omega$ function to $\mathbb{Q}$: $\Xi(n) = \Omega(\text{numerator}(n)) + \Omega(\text{denominator}(n))$.

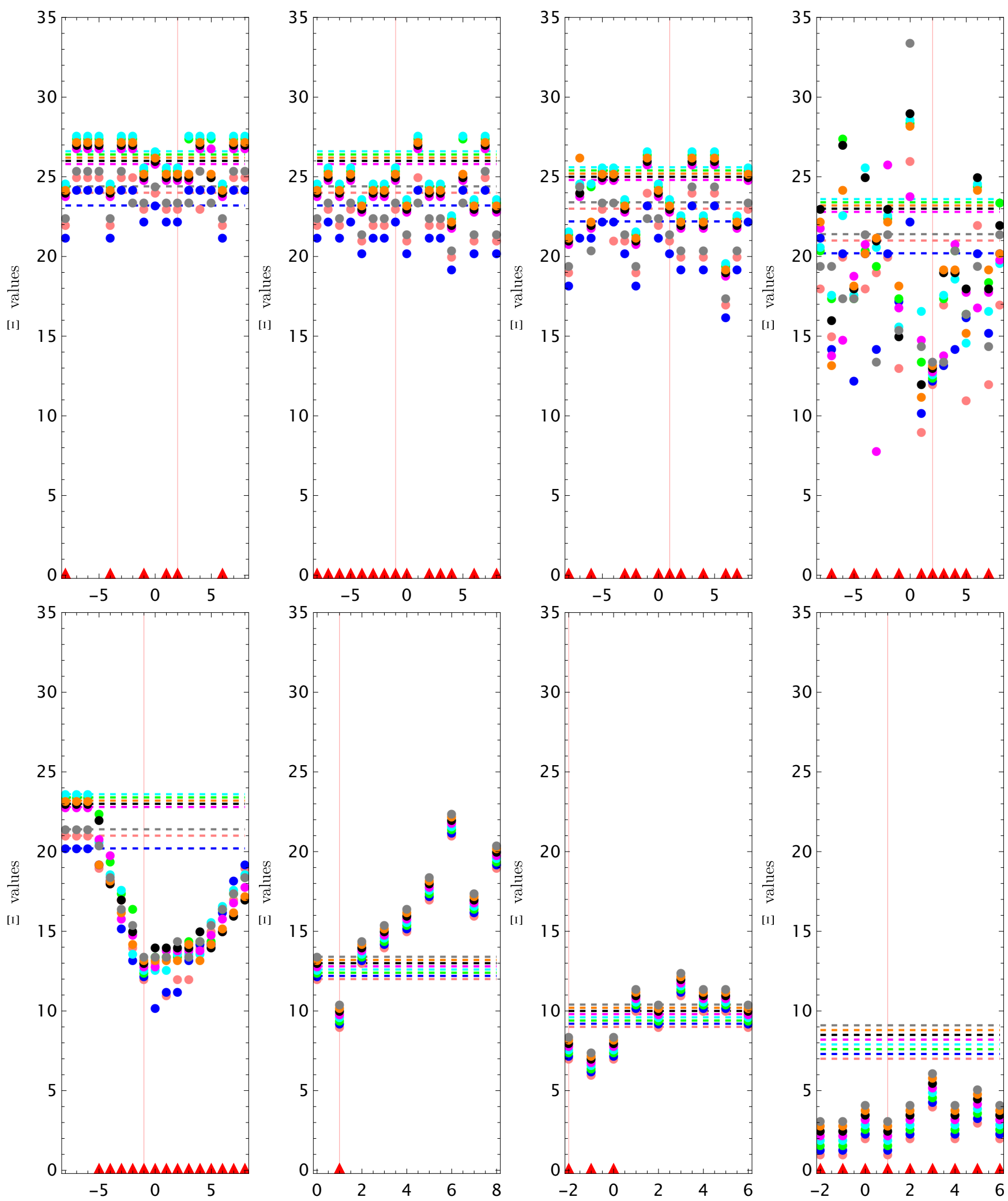

**Fig. 13.** Steps are ordered from left→right and up→down, i.e. $\begin{smallmatrix}1 & 1' & 2 & 3\\ 3' & 4 & 5 & 6\end{smallmatrix}$. Note the steady decrease of $\Xi$ levels. Steps 1 and $1'$ are two division options occurring at the same step and so are 3 and $3'$.

| step | 2 of Figure 13 |
|---|---|
| target | $n_1(j,\kappa,c)$ |
| candidate | $uj-2$ |
| new target | $n_2(j,\kappa,c) = n_1(j,\kappa,c)/(j-2)$ |

| step | 3 of Figure 13 |
|---|---|
| target | $n_2(j,\kappa,c)$ |
| candidate | illustrating first option: $\prod_{x=1}^{\frac{j-1}{2}}(2c-2k+ux-1)$ |
| new target | $n_3(j,\kappa,c) = n_2(j,\kappa,c)/\prod_{x=1}^{\frac{j-1}{2}}(2c-2k+2x-1)$ |

| step | 3′ of Figure 13 |
|---|---|
| target | $n_2(j,\kappa,c)$ |
| candidate | illustrating second option: $\prod_{x=1}^{\frac{j+2u+1}{2}}(2c-2k+2x-1)$ |
| new target | $n_3(j,\kappa,c) = n_2(j,\kappa,c)/\prod_{x=1}^{\frac{j-1}{2}}(2c-2k+2x-1)$ |

| step | 4 of Figure 13 |
|---|---|
| target | $n_3(j,\kappa,c)$ |
| candidate | $C_{u\kappa-1}$ |
| new target | $n_4(j,\kappa,c) = n_3(j,\kappa,c)/C_{\kappa-1}$ |

| step | 5 of Figure 13 |
|---|---|
| target | $n_4(j,\kappa,c)$ |
| candidate | $C_{\frac{j+2u+1}{2}}$ |
| new target | $n_5(j,\kappa,c) = n_4(j,\kappa,c)/C_{\frac{j-3}{2}}$ |

| step | 6 of Figure 13 |
|---|---|
| target | $n_5(j,\kappa,c)$ |
| candidate | $\prod_{x=1}^{\frac{j-1}{2}}(k-x+u)$ |
| new target | None: Processed finished. |

## A Tables of Constants

| $a_0$ | $a_1$ | $a_2$ | $P(n)/(-2n)$ | $T(n)-3n^2$ |
|---|---|---|---|---|
| −7351344 | −32375839 | 46558512 | $n(3+n)(17+n)$ | $72+43n$ |
| −2450448 | −1768477 | 2450448 | $n(1+n)(17+n)$ | $36+39n$ |
| −1081080 | −16147379 | 23279256 | $n(5+n)(15+n)$ | $96+43n$ |
| −793800 | −232217 | 322560 | $(4+n)(9+n)^2$ | $100+39n$ |
| −504504 | −16140515 | 23279256 | $n(7+n)(13+n)$ | $112+43n$ |
| −436590 | 7989199 | −11531520 | $(3+n)(9+n)(11+n)$ | $120+43n$ |
| −180180 | −1001393 | 1441440 | $(2+n)(5+n)(13+n)$ | $84+39n$ |
| −174636 | −8069449 | 11639628 | $n(9+n)(11+n)$ | $120+43n$ |
| −72072 | −850133 | 1225224 | $n(5+n)(13+n)$ | $84+39n$ |
| −72072 | −251099 | 360360 | $n(3+n)(13+n)$ | $56+35n$ |
| −72072 | −52279 | 72072 | $n(1+n)(15+n)$ | $32+35n$ |
| −45045 | 124048 | −180180 | $(1+n)(3+n)(13+n)$ | $56+35n$ |
| −28028 | −999391 | 1441440 | $(2+n)(7+n)(13+n)$ | $112+43n$ |
| −27720 | −20417 | 27720 | $n(1+n)(11+n)$ | $24+27n$ |
| −22050 | 27649 | −40320 | $(3+n)(7+n)(9+n)$ | $80+35n$ |
| −19305 | 424423 | −612612 | $(1+n)(5+n)(15+n)$ | $96+43n$ |
| −17640 | −250007 | 360360 | $n(7+n)(9+n)$ | $80+35n$ |
| −14700 | −153907 | 221760 | $(2+n)(7+n)(9+n)$ | $80+35n$ |
| −10395 | 124741 | −180180 | $(1+n)(5+n)(11+n)$ | $72+35n$ |
| −9450 | 76691 | −110880 | $(1+n)(5+n)(9+n)$ | $60+31n$ |
| −7350 | −62563 | 90090 | $n(7+n)^2$ | $64+31n$ |
| −7007 | 424566 | −612612 | $(1+n)(7+n)(13+n)$ | $112+43n$ |
| −6300 | −14087 | 20160 | $(2+n)(5+n)(9+n)$ | $60+31n$ |
| −6006 | 99839 | −144144 | $(1+n)(5+n)(13+n)$ | $84+39n$ |
| −4900 | −21043 | 30240 | $(2+n)(7+n)^2$ | $64+31n$ |
| −2520 | −1879 | 2520 | $n(1+n)(9+n)$ | $20+23n$ |
| −2450 | 28781 | −41580 | $(1+n)(7+n)^2$ | $64+31n$ |
| −1225 | 367 | −560 | $(3+n)(7+n)^2$ | $64+31n$ |
| −525 | −2413 | 3465 | $n(5+n)(7+n)$ | $48+27n$ |
| −450 | 1151 | −1680 | $(1+n)(5+n)^2$ | $36+23n$ |
| −350 | 1739 | −2520 | $(1+n)(5+n)(7+n)$ | $48+27n$ |
| −180 | −299 | 420 | $n(3+n)(5+n)$ | $24+19n$ |
| −105 | 142 | −210 | $(1+n)(3+n)(7+n)$ | $32+23n$ |
| −18 | 7 | −12 | $(1+n)(3+n)^2$ | $16+15n$ |
| −6 | −5 | 6 | $n(1+n)(3+n)$ | $8+11n$ |

**Table 8.** Examples of convergence to $\frac{a_0}{a_1+a_2\log 2}$

| $a_0$ | $a_1$ | $a_2$ | $P(n)/(-2n)$ | $T(n)-3n^2$ |
|---|---|---|---|---|
| 1 | 1 | −1 | $n(1+n)^2$ | $4+7n$ |
| 9 | 11 | −15 | $n(3+n)^2$ | $16+15n$ |
| 50 | 147 | −210 | $n(5+n)^2$ | $36+23n$ |
| 60 | 47 | −60 | $n(1+n)(5+n)$ | $12+15n$ |
| 90 | −79 | 120 | $(1+n)(3+n)(5+n)$ | $24+19n$ |
| 420 | 319 | −420 | $n(1+n)(7+n)$ | $16+19n$ |
| 420 | 887 | −1260 | $n(3+n)(7+n)$ | $32+23n$ |
| 900 | 361 | −480 | $(2+n)(5+n)^2$ | $36+23n$ |
| 1890 | −3443 | 5040 | $(1+n)(3+n)(9+n)$ | $40+27n$ |
| 2100 | 2377 | −3360 | $(2+n)(5+n)(7+n)$ | $48+27n$ |
| 5544 | 50035 | −72072 | $n(5+n)(11+n)$ | $72+35n$ |
| 6468 | −249713 | 360360 | $(1+n)(7+n)(11+n)$ | $96+39n$ |
| 7560 | 19409 | −27720 | $n(3+n)(9+n)$ | $40+27n$ |
| 8316 | −19031 | 27720 | $(1+n)(3+n)(11+n)$ | $48+31n$ |
| 15444 | −49705 | 72072 | $(1+n)(3+n)(15+n)$ | $64+39n$ |
| 22050 | −499279 | 720720 | $(1+n)(7+n)(9+n)$ | $80+35n$ |
| 24255 | −76586 | 110880 | $(3+n)(7+n)(11+n)$ | $96+39n$ |
| 25740 | 200107 | −288288 | $(2+n)(5+n)(15+n)$ | $96+43n$ |
| 37800 | 250427 | −360360 | $n(5+n)(9+n)$ | $60+31n$ |
| 38808 | 849671 | −1225224 | $n(7+n)(11+n)$ | $96+39n$ |
| 39690 | −1997851 | 2882880 | $(1+n)(9+n)^2$ | $100+39n$ |
| 41580 | 154327 | −221760 | $(2+n)(5+n)(11+n)$ | $72+35n$ |
| 58212 | 3997025 | −5765760 | $(2+n)(9+n)(11+n)$ | $120+43n$ |
| 79380 | 2123957 | −3063060 | $n(9+n)^2$ | $100+39n$ |
| 83160 | 251561 | −360360 | $n(3+n)(11+n)$ | $48+31n$ |
| 87318 | −8491859 | 12252240 | $(1+n)(9+n)(11+n)$ | $120+43n$ |
| 97020 | 1999321 | −2882880 | $(2+n)(7+n)(11+n)$ | $96+39n$ |
| 132300 | 3997907 | −5765760 | $(2+n)(9+n)^2$ | $100+39n$ |
| 198450 | −1227581 | 1774080 | $(3+n)(9+n)^2$ | $100+39n$ |
| 216216 | 852707 | −1225224 | $n(3+n)(15+n)$ | $64+39n$ |
| 360360 | 263111 | −360360 | $n(1+n)(13+n)$ | $28+31n$ |
| 630630 | −3990557 | 5765760 | $(3+n)(7+n)(13+n)$ | $112+43n$ |
| 918918 | −3383801 | 4900896 | $(1+n)(3+n)(17+n)$ | $72+43n$ |
| 1746360 | 2475007 | −3548160 | $(4+n)(9+n)(11+n)$ | $120+43n$ |
| 46558512 | 33464927 | −46558512 | $n(1+n)(19+n)$ | $40+43n$ |

**Table 9.** Examples of convergence to $\frac{a_0}{a_1+a_2\log 2}$. The entry in blue is the one reported by the Ramanujan Project.

| $a_0$ | $a_1$ | $a_2$ | $P(n)/(-2n)$ | $T(n)-3n^2$ |
|---|---|---|---|---|
| −80281600 | −10675439 | 9459450 | $(2+n)(3+n)(14+n)$ | $45+35n$ |
| −15728640 | 392683 | −727650 | $(4+n)^2(12+n)$ | $65+35n$ |
| −13107200 | −263867 | −28350 | $(4+n)(5+n)(10+n)$ | $55+31n$ |
| −10485760 | −93699 | −103950 | $(4+n)(5+n)(12+n)$ | $65+35n$ |
| −7372800 | −884203 | 727650 | $(2+n)(3+n)(12+n)$ | $39+31n$ |
| −6291456 | 149419 | −257250 | $(4+n)(8+n)^2$ | $81+35n$ |
| −5242880 | −86807 | 9450 | $(5+n)(6+n)(10+n)$ | $77+35n$ |
| −3932160 | −116317 | 3150 | $(4+n)(5+n)(8+n)$ | $45+27n$ |
| −3276800 | −158859 | 22050 | $(2+n)(5+n)(10+n)$ | $33+27n$ |
| −2621440 | −48609 | −1050 | $(5+n)(6+n)(8+n)$ | $63+31n$ |
| −2359296 | −168445 | 103950 | $(2+n)(4+n)(12+n)$ | $39+31n$ |
| −491520 | 50593 | −66150 | $(3+n)(4+n)(10+n)$ | $55+31n$ |
| −327680 | −21271 | 9450 | $(2+n)(4+n)(10+n)$ | $33+27n$ |
| −230400 | −1909 | −22050 | $(2+n)^2(5+n)$ | $9+11n$ |
| −196608 | −184547 | 198450 | $(3+n)(6+n)(10+n)$ | $77+35n$ |
| −163840 | −4981 | −22050 | $n(5+n)(8+n)$ | $9+19n$ |
| −131072 | −2951 | −630 | $(4+n)^2(8+n)$ | $45+27n$ |
| −122880 | −13079 | 9450 | $(2+n)(3+n)(10+n)$ | $33+27n$ |
| −61440 | −2467 | −3150 | $(2+n)(4+n)(5+n)$ | $15+15n$ |
| −61440 | 791 | −9450 | $n(5+n)(6+n)$ | $7+15n$ |
| −51200 | 2839 | −9450 | $n(4+n)(5+n)$ | $5+11n$ |
| −49152 | −1919 | 90 | $(4+n)^2(6+n)$ | $35+23n$ |
| −36864 | −2693 | −450 | $(2+n)(4+n)(6+n)$ | $21+19n$ |
| −18432 | −419 | −3150 | $n(4+n)(6+n)$ | $7+15n$ |
| −18432 | −419 | −3150 | $n(3+n)(8+n)$ | $9+19n$ |
| −8192 | −487 | −54 | $(4+n)^3$ | $25+19n$ |
| −3072 | 121 | −630 | $n(4+n)^2$ | $5+11n$ |
| −2048 | −129 | −90 | $(2+n)(4+n)^2$ | $15+15n$ |
| −2048 | −43 | −30 | $(3+n)(4+n)(6+n)$ | $35+23n$ |
| −768 | −77 | −18 | $(2+n)(3+n)(4+n)$ | $15+15n$ |
| −288 | 31 | −90 | $n(2+n)(3+n)$ | $3+7n$ |

**Table 10.** Examples of convergence to $\frac{a_0}{a_1+a_2G}$

| $a_0$ | $a_1$ | $a_2$ | $P(n)/(-2n)$ | $T(n)-3n^2$ |
|---|---|---|---|---|
| 192 | 13 | 18 | $(2+n)^2(3+n)$ | $9+11n$ |
| 384 | 1 | 90 | $n(3+n)(4+n)$ | $5+11n$ |
| 2304 | 389 | 450 | $n(3+n)(6+n)$ | $7+15n$ |
| 3072 | 179 | −18 | $(3+n)(4+n)^2$ | $25+19n$ |
| 4608 | 133 | 450 | $(2+n)^2(4+n)$ | $9+11n$ |
| 4608 | 383 | −90 | $(2+n)(3+n)(6+n)$ | $21+19n$ |
| 11520 | −1373 | 3150 | $n(2+n)(4+n)$ | $3+7n$ |
| 12288 | 973 | −750 | $(3+n)(6+n)^2$ | $49+27n$ |
| 12288 | 1145 | −630 | $(2+n)(3+n)(8+n)$ | $27+23n$ |
| 16384 | −543 | 1050 | $(3+n)(4+n)(8+n)$ | $45+27n$ |
| 81920 | 3983 | 1350 | $(4+n)^2(5+n)$ | $25+19n$ |
| 89600 | −10891 | 22050 | $n(2+n)(5+n)$ | $3+7n$ |
| 98304 | 2263 | 150 | $(4+n)(6+n)^2$ | $49+27n$ |
| 98304 | 35389 | −36750 | $(3+n)(6+n)(8+n)$ | $63+31n$ |
| 122880 | 6563 | 3150 | $(2+n)(5+n)(6+n)$ | $21+19n$ |
| 147456 | 21365 | 22050 | $n(4+n)(8+n)$ | $9+19n$ |
| 163840 | 6789 | 450 | $(4+n)(5+n)(6+n)$ | $35+23n$ |
| 262144 | 710401 | −771750 | $(3+n)(8+n)^2$ | $81+35n$ |
| 294912 | 18013 | −3150 | $(2+n)(4+n)(8+n)$ | $27+23n$ |
| 524288 | 27787 | −22050 | $(4+n)(6+n)(10+n)$ | $77+35n$ |
| 786432 | 19099 | −5250 | $(4+n)(6+n)(8+n)$ | $63+31n$ |
| 983040 | 25979 | −450 | $(5+n)(6+n)^2$ | $49+27n$ |
| 1310720 | 1723 | 28350 | $(4+n)^2(10+n)$ | $55+31n$ |
| 2949120 | 168821 | 22050 | $(2+n)(5+n)(8+n)$ | $27+23n$ |
| 9830400 | −1833409 | 2182950 | $(3+n)(4+n)(12+n)$ | $65+35n$ |
| 12582912 | 184025 | −7350 | $(5+n)(8+n)^2$ | $81+35n$ |
| 39321600 | 1965547 | −727650 | $(2+n)(5+n)(12+n)$ | $39+31n$ |
| 165150720 | 12969199 | −9459450 | $(2+n)(4+n)(14+n)$ | $45+35n$ |
| 330301440 | 17687791 | −9459450 | $(2+n)(5+n)(14+n)$ | $45+35n$ |

**Table 11.** Examples of convergence to $\frac{a_0}{a_1+a_2G}$

| $c$ | $a_0$ | $a_1$ | $a_2$ |
|---|---|---|---|
| 3 | 2 | 1 | 0 |
| 4 | 1 | 1 | −1 |
| 5 | −1 | −3 | 4 |
| 6 | −2 | −17 | 24 |
| 7 | −3 | −67 | 96 |
| 8 | 12 | 667 | −960 |
| 9 | 5 | 666 | −960 |
| 10 | 30 | 9319 | −13440 |
| 11 | −105 | −74537 | 107520 |
| 12 | −280 | −447187 | 645120 |
| 13 | 126 | 447173 | −645120 |
| 14 | 63 | 491884 | −709632 |
| 15 | 231 | 3935051 | −5677056 |
| 16 | 2772 | 102311095 | −147603456 |
| 17 | −1287 | −102310996 | 147603456 |
| 18 | −6006 | −1023109531 | 1476034560 |
| 19 | 45045 | 16369749493 | −23616552960 |
| 20 | 720720 | 556571437717 | −802962800640 |

**Table 12.** Examples of convergence to $\frac{a_0}{a_1+a_2\log 2}$ for $R_c$

| $j$ | $\overset{\alpha}{\alpha}_j$ | $\overset{\beta}{\alpha}_j$ |
| --- | --- | --- |
| 3 | −1 | 4 |
| 5 | 19 | 234 |
| 7 | 5818/3 | 254456/3 |
| 9 | 667115 | 60003486 |
| 11 | 467946090 | 71121907440 |
| 13 | 554143204110 | 127451285438100 |
| 15 | 994115449382940 | 322092692148962160 |
| 17 | 2516347061651130075 | 1092094185270706446150 |
| 19 | 8546069024090201027250 | 4785798287838257081935200 |
| 21 | 37508692924557081882027450 | 26331102038134635548392485900 |
| 23 | 206659254109760483703789089700 | 177726957997323983116663150902000 |
| 25 | 1396637676485497608584841260027550 | 1444123356588023432320434243315206700 |
| 27 | 11361110319787394788017568214856502500 | 13905999029609441333101619589964946580000 |
| 29 | 109509742351999832489255793094925601037500 | 156598931559029451368898717824937174831465000 |
| 31 | 1234320809247763942235545044494798498436195000 | 2039097976865181167119056627863149102390546140000 |
| 33 | 16085205915675471439195309128783843538512283666875 | 30401039180587356456007967587920548312623820610393750 |
| 35 | 239989379884263177615577263747245812249369757283461250 | 514537230471714428505965482811829861362838523445500920000 |

| $j$ | $\overset{\alpha}{\beta}_j$ | $\overset{\beta}{\beta}_j$ |
| --- | --- | --- |
| 3 | −1/3 | −14/3 |
| 5 | −17 | −8 |
| 7 | −758 | −27820 |
| 9 | −302117 | −23010044 |
| 11 | −1091480994/5 | −146282046156/5 |
| 13 | −262476468810 | −54596049230880 |
| 15 | −475443072646380 | −141682352738003640 |
| 17 | −1211573031414907725 | −489475664504671450500 |
| 19 | −4135193781750207709650 | −2175112041708995560914300 |
| 21 | −18218507239728799899288030 | −12097088912487772715204794320 |
| , 23 | −100680148628028059378172563700 | −82356361704096372069838207986600 |
| 25 | −682078864161239229949889893754850 | −673893917980353010760236819146271800 |
| 27 | −5559692282317104149119150499246482500 | −6527078739785105011529098668023829975000 |
| 29 | −53681246247288656939970174534335708392500 | −73865394837022289182570623863734010339760000 |
| 31 | −605939306349175124039948450466713432304279000 | −965867254322525126192328035746702493817188406000 |
| 33 | −7906287653442943862409767973858652552097126548125 | −14452693830499521315903006473321900713406050369972500 |
| 35 | −118090012323922699712409299094969252935070156336941250 | −245391045609131483699190960852072336578359955374403917500 |

**Table 13.** The first $\overline{\alpha}_j = \left\{ \{\overset{\alpha}{\alpha}_j, \overset{\beta}{\alpha}_j\} \right\}$ values and the first $\overline{\beta}_j = \left\{ \{\overset{\alpha}{\beta}_j, \overset{\beta}{\beta}_j\} \right\}$ values.

| $j$ | $\kappa$ | $c$ | $j'$ | $\kappa'$ | $c'$ | $Q_{j,\kappa,c}/Q_{j',\kappa',c'}$ |
|---|---|---|---|---|---|---|
| 3 | 3 | 4 | 5 | 5 | 3 | 1/3 |
| 3 | 3 | 4 | 7 | 5 | 3 | 1/3 |
| 3 | 5 | 6 | 5 | 5 | 7 | 5/7 |
| 3 | 5 | 6 | 7 | 5 | 7 | 5/7 |
| 3 | 6 | 6 | 5 | 7 | 5 | 3/7 |
| 3 | 6 | 8 | 5 | 9 | 6 | 1/3 |
| 3 | 9 | 4 | 7 | 10 | 5 | 7/15 |
| 3 | 9 | 9 | 5 | 9 | 10 | 3/5 |
| 3 | 9 | 10 | 5 | 11 | 8 | 5/13 |
| 3 | 9 | 12 | 5 | 13 | 9 | 1/3 |
| 3 | 15 | 8 | 7 | 15 | 8 | 7/15 |
| 5 | 3 | 4 | 5 | 5 | 3 | 1/3 |
| 5 | 3 | 4 | 7 | 5 | 3 | 1/3 |
| 5 | 5 | 5 | 7 | 6 | 4 | 5/9 |
| 5 | 5 | 6 | 5 | 6 | 4 | 3/5 |
| 5 | 6 | 4 | 7 | 5 | 6 | 5/3 |
| 5 | 7 | 10 | 5 | 10 | 7 | 25/49 |
| 5 | 8 | 10 | 7 | 8 | 11 | 9/11 |
| 5 | 10 | 6 | 9 | 11 | 6 | 3/5 |
| 5 | 10 | 11 | 7 | 10 | 12 | 7/9 |
| 5 | 10 | 12 | 5 | 13 | 10 | 3/5 |
| 7 | 5 | 5 | 7 | 6 | 4 | 5/9 |
| 7 | 14 | 9 | 11 | 14 | 8 | 99/133 |
| 9 | 10 | 9 | 11 | 11 | 8 | 27/35 |
| 9 | 10 | 9 | 13 | 11 | 8 | 27/35 |

**Table 14.** Examples of nontrivial ratios within Kosovo. The values linked to each other are shown in Figure 14. There seems to be no straightforward relation between those (except that ratios never feature even numbers).

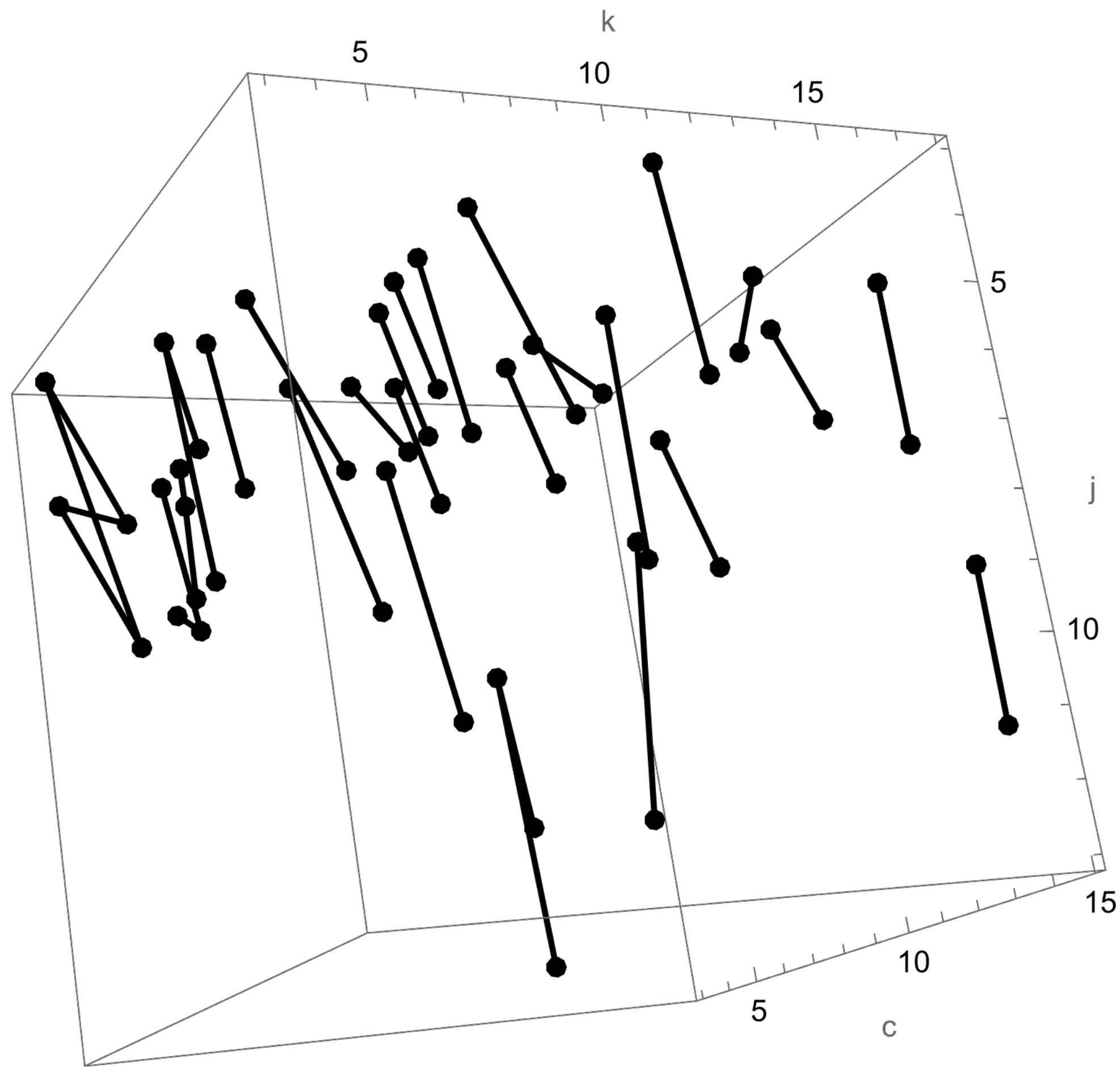

**Fig. 14.** Any connected points are $Q_{j,\kappa,c}$ values whose ratio is a rational.

## B Mathematica Code

```
(* code snippet "1. Montenegro" *)
NumQ[k_,c_]:=2k+1+ContinuedFractionK[-2n^2(n+2k)(n+c),3n^2+(3+4k)n+2k+1,{n,1,5000}];
vD[a_,b_,k_,c_]:=If[c<2,a+b*c,-2c(2c-1)(2(c-k)-1)^2vD[a,b,k,c-2]+(8c^2+(2-8k)c-2k+1)vD[a,b,k,c-1]];
vG[a_,b_,k_,c_]:=((2c-1)!!)^2Catalan+Product[(2(c-j)-1),{j,0,k-1}]*vD[a,b,k,c-1];
vd[k_]:=4^(k-1)/(2k-1)/CatalanNumber[k-1];vr[k_]:=vd[k](-1)^(k)(1-2k)/((2k)!(2k-3)!!);
va[k_]:=vr[k]vD[1,-2,1,k-1];vb[k_]:=-vr[k](2k-3)^2vD[1,12,2,k-1]-va[k];
QFor[k_,c_]:=vd[k](2c)!/vG[va[k],vb[k],k,c];
Print[Union[Flatten[Table[N[QFor[k,c],200]==N[NumQ[k,c],200],{c,1,14},{k,1,14}]]]];ClearAll["Global‘*"];
```

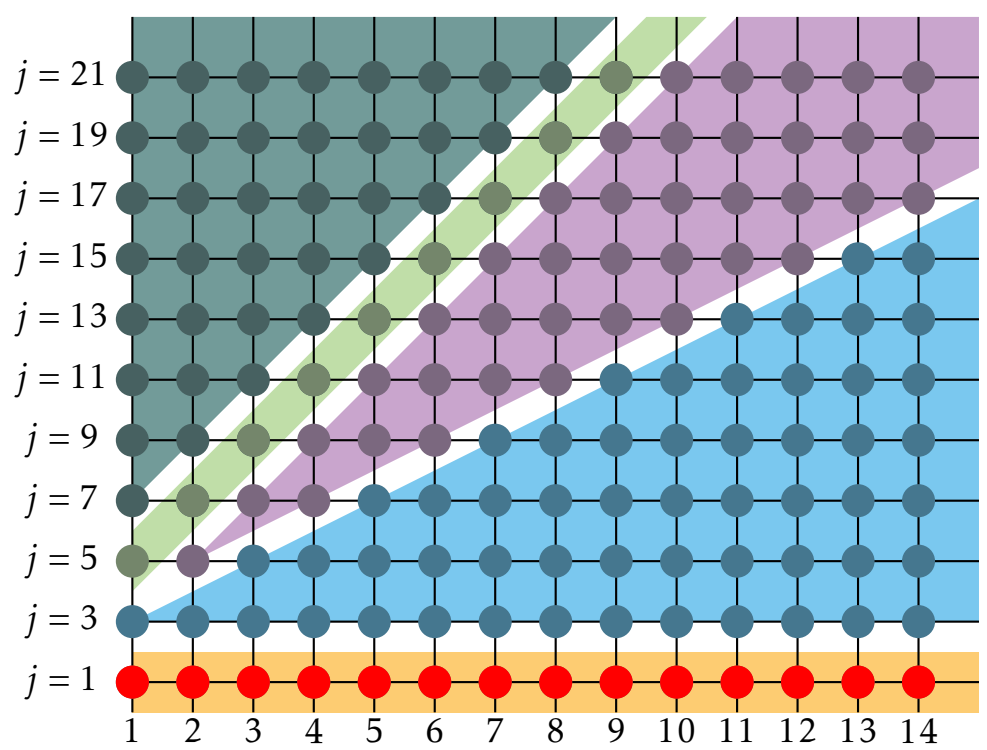

**Fig. 15.** The points on which Montenegro is tested for $1 \leq c \leq 14$ by the snippet "1. Montenegro" are shown in red.

```
(* code snippet "2. Bosnia" *)
NumQ[j_,c_]:=-j+ContinuedFractionK[-2n(n-2)(n+c)(n+j-1),3n^2+(2j-3)n-j,{n,1,5000}];
vD[j_,c_]:=If[c<2,1+(15-4j)c,-2c(2c-j)(2c+1)(2c-j+2)vD[j,c-2]+(8c^2+(14-4j)c-3(j-2))vD[j,c-1]];
g[j_,c_]:=Product[(2+2i-j),{i,1,(j-1)/2}](2c)!/2;
h[j_,c_]:=Product[(2c-2i-1),{i,0,(j-5)/2}];
Q[j_,c_]:=Simplify[g[j,c]/(vD[j,c-1]*h[j,c])];
Print[Union[Flatten[Table[N[NumQ[j,c]==Q[j,c],200],{j,5,13,2},{c,1,14}]]]];ClearAll["Global‘*"];
```

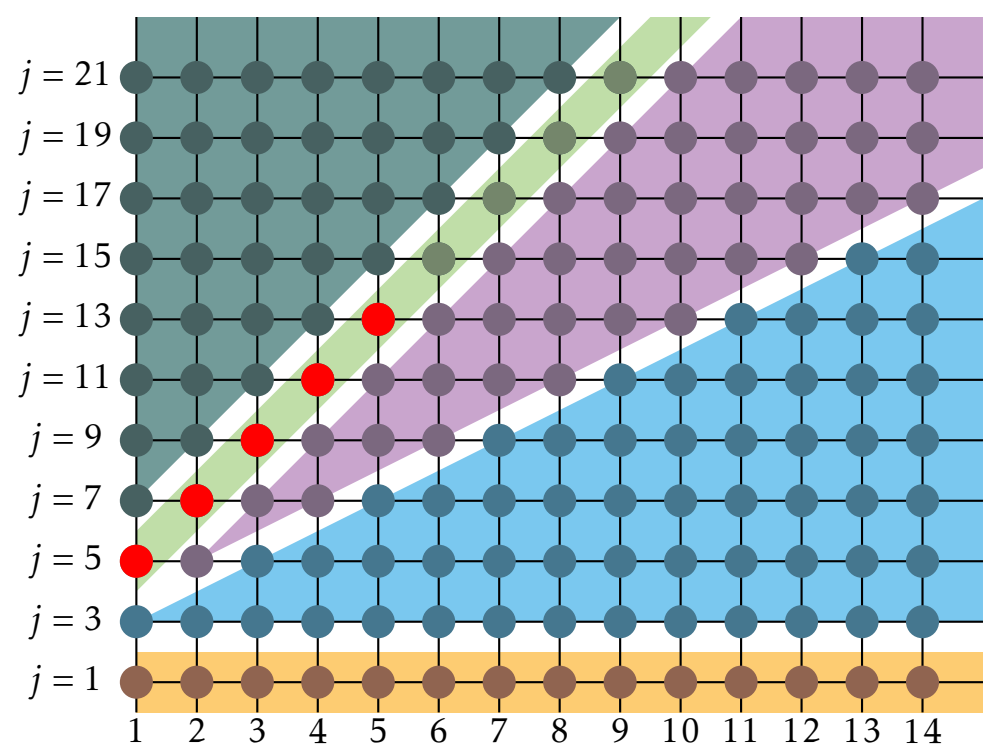

**Fig. 16.** The points on which Bosnia & Herzegovina is tested for $1 \leq c \leq 14$ by the snippet "2. Bosnia" are shown in red.

```
(* code snippet "3. Northern Balkans"*)
NumQ[j_,k_,c_]:=j(2-j+2k)+ContinuedFractionK[-2n(c+n)(j+n-1)(1-j+2k+n),3n^2+(3+4k)n+j(2-j+2k),{n,1,15000}];
vD[ab_,j_,k_,c_]:=If[c<2,ab[[1]]+ab[[2]]c,-2c(2c-j)(2c-2k+j-2)(2c-2k-1)vD[ab,j,k,c-2]+(8c^2+(2-8k)c+(j-2)(2k-j))vD[ab,j
    ,k,c-1]];
```

```
4  f[j_,k_,c_]:=Product[(2c-2k+2i-1)(k-i+1),{i,1,(j-1)/2}]CatalanNumber[(j-3)/2](j-2)(2k-1)*(2c-1)!!^2*CatalanNumber[k-1];
5  g[j_,k_,c_]:=Product[(2c-2i+1)(2k-2i+1),{i,1,(j-1)/2}](2c)!2^(2k+(j-7)/2);
6  h[j_,k_,c_]:=Product[2c-2i-1,{i,0,(j-3)/2}]Product[2c-2i-1,{i,0,k-1}];
7  Q[j_,k_,c_,ab_]:=Simplify[g[j,k,c]/(vD[ab,j,k,c-1]h[j,k,c]+f[j,k,c]Catalan)];
8  Ratio:=Function[{j,k,c},r=N[NumQ[j,k,c],2000];v=FindIntegerNullVector[{1,r,N[Catalan*r,2000]}];-v[[1]]/v[[2]]];
9  GenAB[j_]:=Table[({a,b}/.Solve[Table[{a,b}.{1,c-1}==g[j,k,c]/Ratio[j,k,c]/h[j,k,c],{c,1,2}],{a,b}])[[1]],{k,1,6}];
10 result={};For[j=-7,j<=13,If[j==1,j=3];AB=GenAB[j];AppendTo[result,Union[Flatten[Table[N[Q[j,k,c,AB[[k]]]==NumQ[j,k,c
       ],20],{k,1,6},{c,1,7}]]]];j+=2];Print[Union[Flatten[result]]];ClearAll["Global‘*"];
```

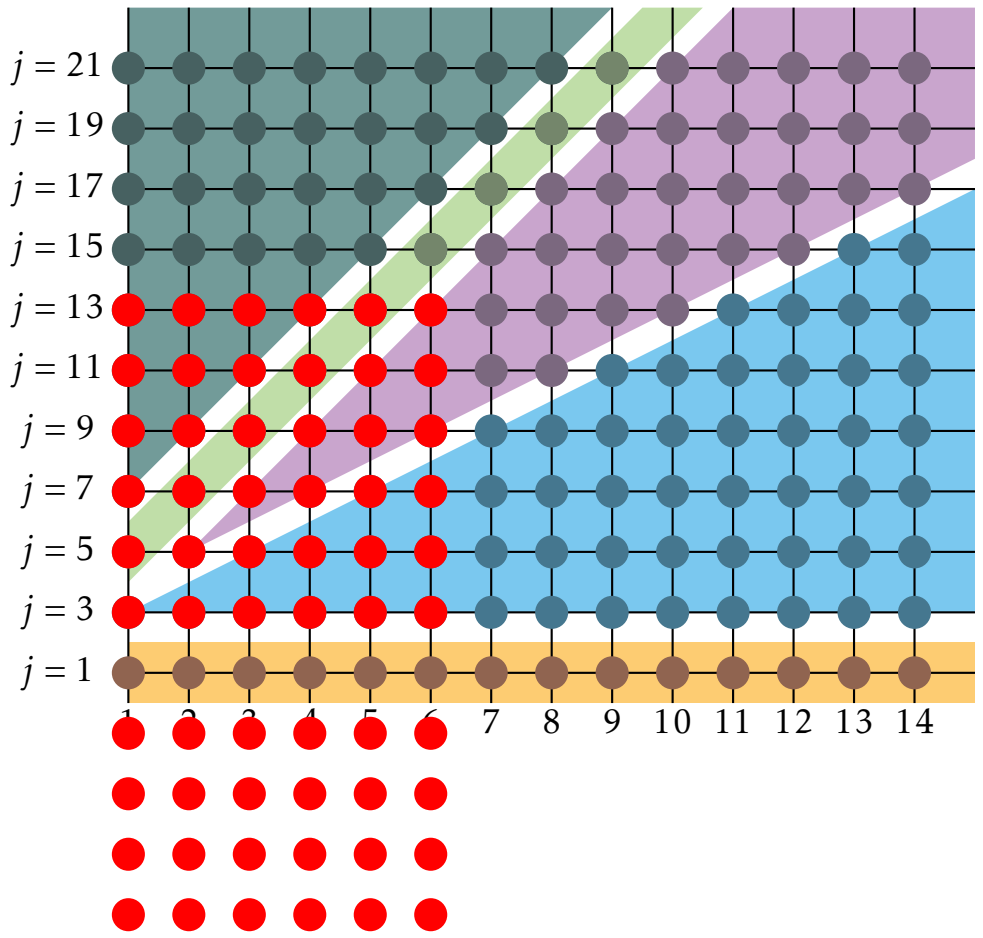

**Fig. 17.** The points on which the $c$-level master formula is tested for $1 \leq c \leq 7$ by the snippet "`3. Northern Balkans`" are shown in red. Note that the Montenegro line does not obey that same $c$-level formula but was validated by snippet "`1. Montenegro`". Note that the formula is also tested for negative $j$ values (not shown as red dots in the diagram).

```mathematica
(* code snippet "4. Kosovo" *)
NumQ[j_,k_,c_]:=j(2-j+2k)+ContinuedFractionK[-2n(c+n)(j+n-1)(1-j+2k+n),3n^2+(3+4k)n+j(2-j+2k),{n,1,15000}];
vD:=Function[{ab,j,k,c},If[c<2,ab[[1]]+ab[[2]]c,-2c(2c-j)(2c-2k+j-2)(2c-2k-1)vD[ab,j,k,c-2]+(8c^2+(2-8k)c+(j-2)(2k-j))
    vD[ab,j,k,c-1]]];
f[j_,k_,c_]:=Product[(2c-2k+2i-1)(k-i+1),{i,1,(j-1)/2}]CatalanNumber[(j-3)/2](j-2)(2k-1)*(2c-1)!!^2*CatalanNumber[k-1];
g[j_,k_,c_]:=Product[(2c-2i+1)(2k-2i+1),{i,1,(j-1)/2}](2c)!2^(2k+(j-7)/2);
h[j_,k_,c_]:=Product[2c-2i-1,{i,0,(j-3)/2}]Product[2c-2i-1,{i,0,k-1}];
Q[j_,k_,c_,ab_]:=Simplify[g[j,k,c]/(vD[ab,j,k,c-1]*h[j,k,c]+f[j,k,c]Catalan)];
l[n_,j_,k_]:=(-1)^(k+1)(2k)!^2/k!/2^(3*k-2)/Product[(k-i)(2k-2i-1)^2,{i,0,(j-3)/2}]/(2k-j)/(2k-1)/(n((2k-j-2)(3-2k)-1)
    +1);
kD:=Function[{n,ab,j,k},If[k<2,ab[[1]]+ab[[2]]k,(2k+2j-9-2n)(2k+j-8-2n)(-2k+5-j)(2k+j-6)*kD[n,ab,j,k-2]+(8k^2+k(10j
    -48-8n)+(3j^2-(28+4n)j+68+18n))kD[n,ab,j,k-1]]];
Descend:=Function[{j,k},abh[0]={{-1,4},{19,234},{5818/3,254456/3},{667115,60003486},{467946090,71121907440},
{554143204110,127451285438100},{994115449382940,322092692148962160},{2516347061651130075,1092094185270706446150},
{8546069024090201027250,4785798287838257081935200},{37508692924557081882027450,26331102038134635548392485900},
{206659254109760483703789089700,177726957997323983116663150902000},
{1396637676485497608584841260027550,1444123356588023432320434243315206700},
{11361110319787394788017568214856502500,13905999029609441333101619589964946580000},
{109509742351999832489255793094925601037500,156598931559029451368898717824937174831465000},
{1234320809247763942235545044494798498436195000,2039097976865181167119056627863149102390546140000},
{16085205915675471439195309128783843538512283666875,30401039180587356456007967587920548312623820610393750},
{239989379884263177615577263747245812249369757283461250,514537230471714428505965482811829861362838523445500920000}};
abh[1]={{-1/3,-14/3},{-17,-8},{-758,-27820},{-302117,-23010044},{-1091480994/5,-146282046156/5},
{-262476468810,-54596049230880},{-475443072646380,-141682352738003640},{-1211573031414907725,-489475664504671450500},
{-4135193781750207709650,-2175112041708995560914300},{-18218507239728799899288030,-12097088912487772715204794320},
{-100680148628028059378172563700,-82356361704096372069838207986600},
{-682078864161239229949889893754850,-673893917980353010760236819146271800},
{-5559692282317104149119150499246482500,-6527078739785105011529098668023829975000},
{-53681246247288656939970174534335708392500,-73865394837022289182570623863734010339760000},
{-605939306349175124039948450466713432304279000,-965867254322525126192328035746702493817188406000},
```

```
28 {-790628765344294386240976797385865255209712654812 5,-1445269383049952131590300647332190071340605036997 2500},
```

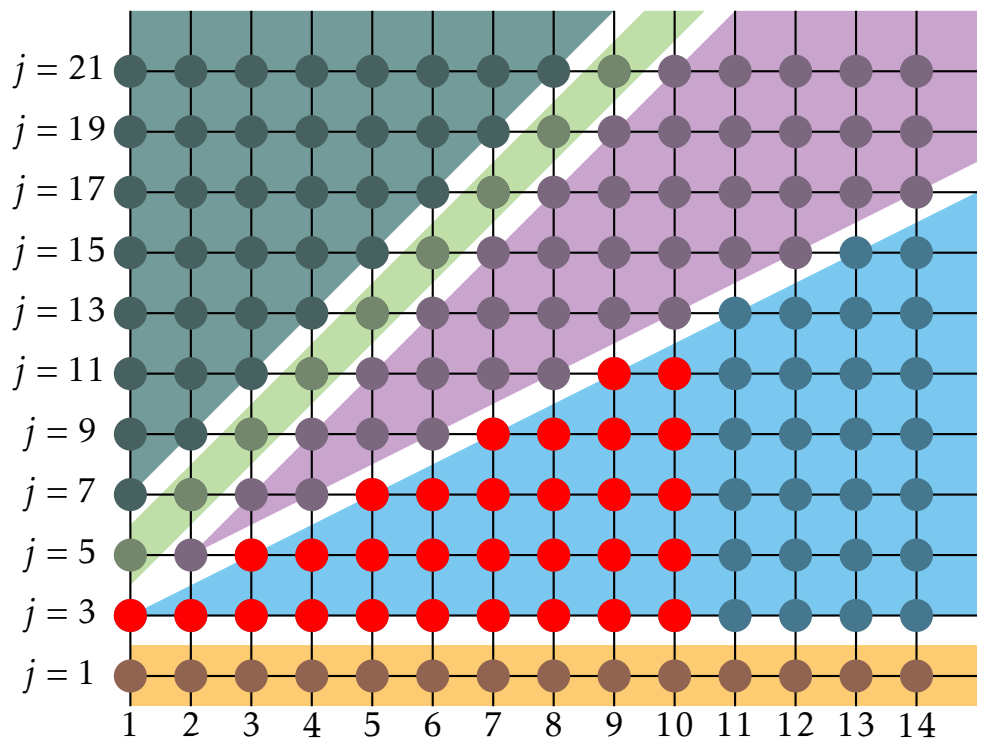

**Fig. 18.** The points on which Kosovo is tested for $1 \leq c \leq 7$ by the snippets `"4. Kosovo"` and `"18. Altogether"` are shown in red.

```
1 (* code snippet "5. resolution" *)
2 NumQ[j_,k_,c_]:=j(2-j+2k)+ContinuedFractionK[-2n(c+n)(j+n-1)(1-j+2k+n),3n^2+(3+4k)n+j(2-j+2k),{n,1,10000}];
3 f[j_,k_,c_]:=Product[(2c-2k+2i-1)(k-i+1),{i,1,(j-1)/2}]CatalanNumber[(j-3)/2](j-2)(2k-1)*(2c-1)!!^2*CatalanNumber[k-1];
4 g[j_,k_,c_]:=Product[(2c-2i+1)(2k-2i+1),{i,1,(j-1)/2}](2c)!2^(2k+(j-7)/2);
5 h[j_,k_,c_]:=Product[2c-2i-1,{i,0,(j-3)/2}]Product[2c-2i-1,{i,0,k-1}];
```

```
6  l:=Function[{n,j,k},(-1)^(k+1)(2k)!^2/k!/2^(3k-2)/Product[(k-i)(2k-2i-1)^2,{i,0,(j-3)/2}]/(2k-j)/(2k-1)/(n((2k-j-2)(3-2
       k)-1)+1)];
7  AB:=Function[{j,k},For[c=1,c<=2,r=N[NumQ[j,k,c],2000];v=FindIntegerNullVector[{1,r,N[Catalan*r,2000]}];d[c]=(-g[j,k,c](
       v[[3]]Catalan+v[[2]])/v[[1]]-f[j,k,c]Catalan)/h[j,k,c];c++];{d[1],d[2]-d[1]}];ABlists:=Function[lim,z[1]=z[0]={};
8  For[j=3,j<=lim,For[w=0,w<=1,e[w]=Table[AB[j,k].{1,w}l[w,j,k],{k,j-2,j-1}];
9  ab[w]={a,b}/.Solve[Table[a+b(u-1)==e[w][[u]],{u,1,2}],{a,b}];AppendTo[z[w],ab[w][[1]]];w++];j+=2];{z[1],z[0]}];
10 Print[ABlists[15]];ClearAll["Global‘*"];
```

```
1  (* code snippet "6. symmetry" *)
2  NumQ[j_,k_,c_]:=j(2-j+2k)+ContinuedFractionK[-2n(c+n)(j+n-1)(1-j+2k+n),3n^2+(3+4k)n+j(2-j+2k),{n,1,15000}];
3  g[j_,k_,c_]:=Product[(2c-2i+1)(2k-2i+1),{i,1,(j-1)/2}](2c)!2^(2k+(j-7)/2);
4  h[j_,k_,c_]:=Product[2c-2i-1,{i,0,(j-3)/2}]Product[2c-2i-1,{i,0,k-1}];
5  zeta[j_,u_]:=Product[2+2i-j,{i,0,u-1}]/2^u;Ratio:=Function[{j,k,c},r=N[NumQ[j,k,c],2000];
6  v=FindIntegerNullVector[{1,r,N[Catalan*r,2000]}];-v[[1]]/v[[2]]];
7  GenAB[j_,k_]:=({a,b}/.Solve[Table[{a,b}.{1,c-1}==g[j,k,c]/Ratio[j,k,c]/h[j,k,c],{c,1,2}],{a,b}])[[1]];
8  tau[j_,u_]:=If[u>(j-1)/2,Sign[zeta[j,u]](2j-2u-3)!!(2u-j)!!(-2)^u,zeta[j,u](-4)^u];
9  For[j=3,j<=13,Print[MatrixForm[Table[{{j,j-u-1},{j-2u,j-u-1},(GenAB[j,j-u-1]/GenAB[j-2u,j-u-1])/tau[j,u]=={1,1}},{u,1,j
       +3}]]];j+=2];ClearAll["Global‘*"];
```

```
1  (* code snippet "7. Croatia" *)
2  NumQ[j_,k_,c_]:=j(2-j+2k)+ContinuedFractionK[-2n(c+n)(j+n-1)(1-j+2k+n),3n^2+(3+4k)n+j(2-j+2k),{n,1,2000}];
3  g[j_,k_,c_]:=Product[(2c-2i+1)(2k-2i+1),{i,1,(j-1)/2}](2c)!2^(2k+(j-7)/2);
4  h[j_,k_,c_]:=Product[2c-2i-1,{i,0,(j-3)/2}]Product[2c-2i-1,{i,0,k-1}];
5  Ratio:=Function[{j,k,c},r=N[NumQ[j,k,c],4000];v=FindIntegerNullVector[{1,r,N[Catalan*r,4000]}];-v[[1]]/v[[2]]];
6  GenAB[j_,k_]:=({a,b}/.Solve[Table[{a,b}.{1,c-1}==g[j,k,c]/Ratio[j,k,c]/h[j,k,c],{c,1,2}],{a,b}])[[1]];
7  psi1[i_,j_]:={-1,14-j,-464+58j-3j^2,27936-4692j+432j^2-15j^3,-2659968+542256j-67836j^2+4260j^3-105j
       ^4,367568640-86278560j+13203480j^2-1139700j^3+51450j^4-945j^5}[[i]];
8  psi2[i_,j_]:={4j-15,306-95j+4j^2,-13360+4646j-357j^2+12j^3,999648-379692j+40368j^2-2457j^3+60j^4,-113885568+
9  46449360j-6124164j^2+513228j^3-22935j^4+420j^5,18333538560-7933530720j+1224286440j^2-126833100j^3+7864950j^4-266175j
       ^5+3780j^6}[[i]];
10 mu[i_,j_]:=-Product[j-2q-2,{q,1,i}]/(-2)^((3j-11-4i)/2);
```

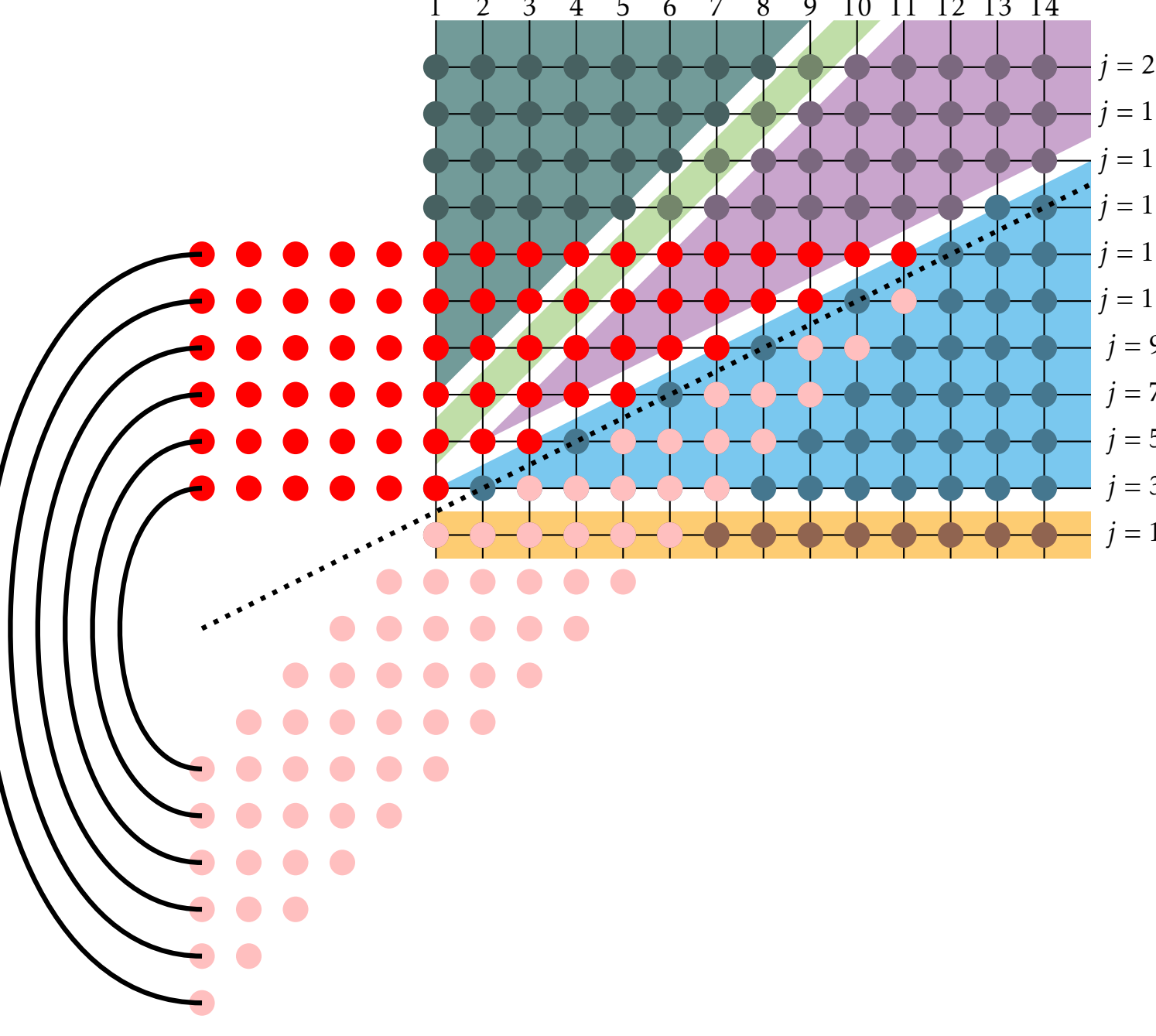

**Fig. 19.** The points on which the Kosovo-Serbia symmetry is tested at the $(\alpha_{j,\kappa}, \beta_{j,\kappa})$-level by the snippet "`6. symmetry`". Tested points are shown in red and their symmetrical correspondents in pink.

```
11 Print[Union[Table[Union[Table[GenAB[j,(j-3)/2-i]=={psi1[1+i,j],psi2[1+i,j]}/mu[i,j],{j,5+2i,37,2}]],{i,0,5}]]];ClearAll
    ["Global`*"];
```

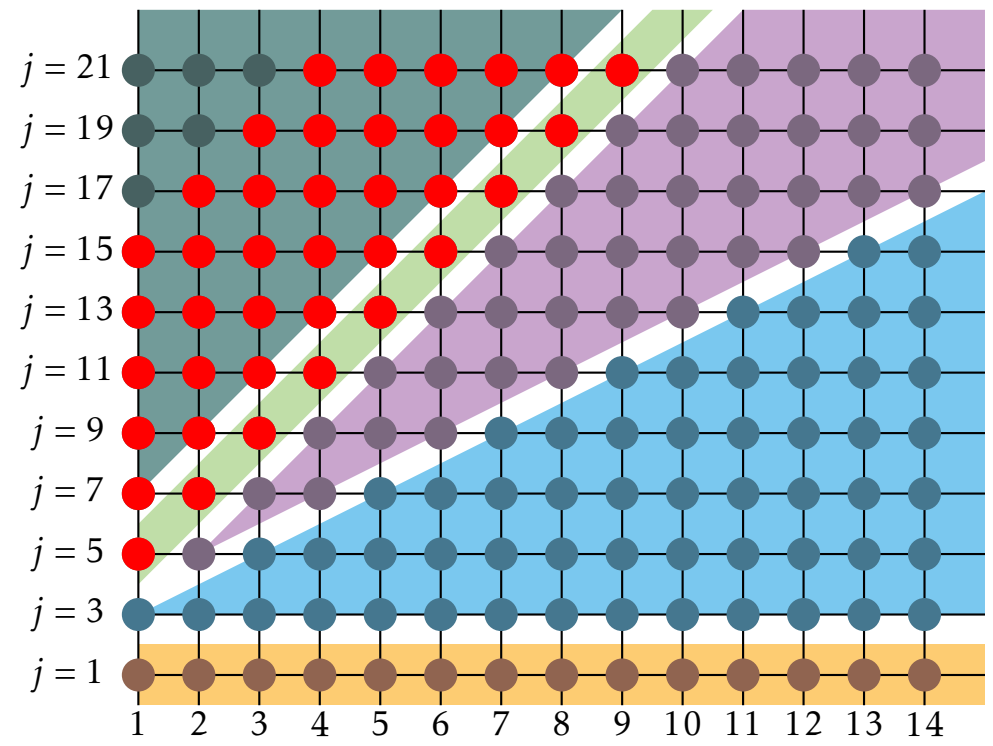

**Fig. 20.** The points on which Croatia is tested. Note that the test is performed directly at the $(\alpha_{j,\kappa}, \beta_{j,\kappa})$-level. Points tested by the snippet "7. Croatia" are shown in red. Points tested beyond the quadrant are not shown.

```
(* code snippet "8. coefficients" *)
coef=CoefficientList[{3198013886925+(145296572850+(5207427225+(4353102000+(877052475+(78210090+(3023055+41580(X-29))(X-27))(X-25))(X-23))(X-21))(X-19))(X-17),-14487726825-(104826150+(452605725+(121200300+(13697775+(640710+10395(X-27))(X-25))(X-23))(X-21))(X-19))(X-17)},X];Print[Table[GCD[coef[[1, i]], coef[[2, i]]], {i, 1, 7}]];ClearAll["Global‘*"];
```

```
(* code snippet "9. ratio" *)
NumQ[j_,k_,c_]:=j(2-j+2k)+ContinuedFractionK[-2n(c+n)(j+n-1)(1-j+2k+n),3n^2+(3+4k)n+j(2-j+2k),{n,1,50000}];
For[j=1,j<=7,For[k=Abs[j-2],k<=7,For[c=1,c<=7,r=N[NumQ[j,k,c],3000];v=FindIntegerNullVector[{1,r,N[Catalan*r,3000]}];
q=(-v[[1]]/v[[3]]==((2c)!2^(2k+(j-7)/2+Floor[1/j])/CatalanNumber[k-1]/(2k-1)/(2c-1)!!^2/(j-2)/Product[(2c-2k+2i-
1)(k-i+1)/(2c-2i+1)/(2k-2i+1),{i,1,(j-1)/2}]/CatalanNumber[(j-3)/2]));Print[{j,k,c,q}];c++];k++];j+=2];ClearAll["Global‘*"];
```

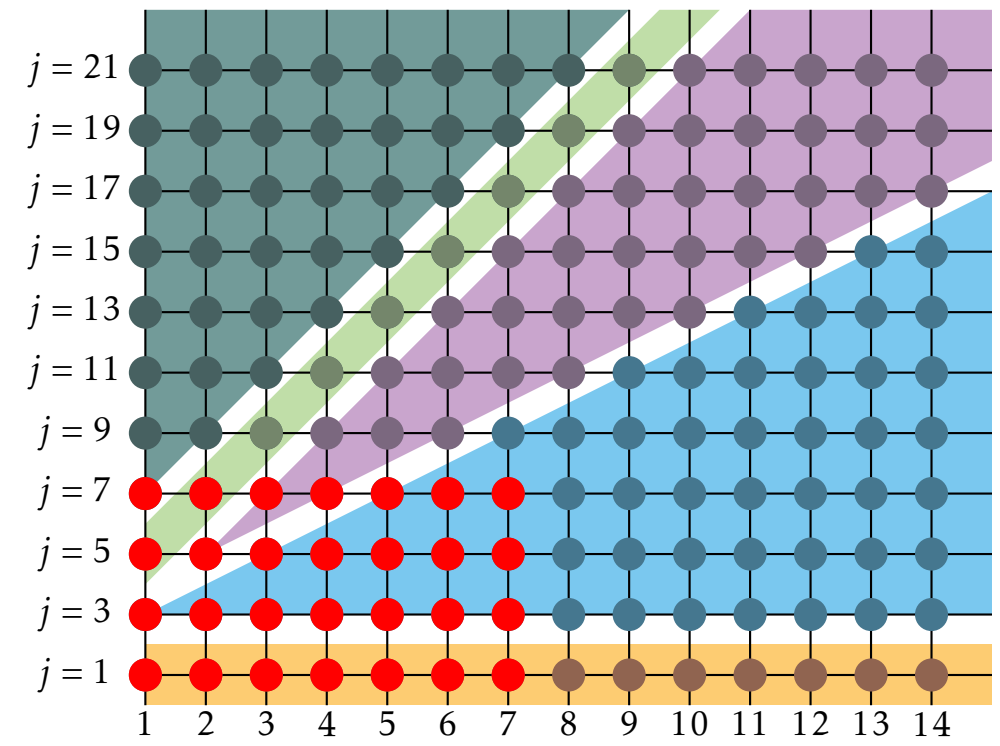

**Fig. 21.** The points on which the ratio conjecture is tested for $1 \leq c \leq 7$ by the snippet "`9. ratio`" are shown in red.

```
(* code snippet "10. log2-a" *)
Q[c_]:=c+ContinuedFractionK[-2n^2,3n+c,{n,1,2000}];
For[c=2,c<=25,r=N[Q[c],80];
v=FindIntegerNullVector[{1,r,N[Log[2]*r,80]}];
Print[{c,Simplify[2/LerchPhi[1/2,1,-1+c]==-v[[1]]/(v[[2]]+Log[2]v[[3]])==1/(Log[2]2^(c-2)-Sum[2^(c-j-2)/j,{j,c-2}])]}];
    c++];ClearAll["Global‘*"];
```

```
(* code snippet "11. log2-b"*)
R[c_]:=c+ContinuedFractionK[-2n^2-2n,3n+c,{n,1,2000}];
For[c=3,c<=20,r=N[R[c],200];v=FindIntegerNullVector[{1,r,N[Log[2]*r,200]}];
Print[{c,2^(c-4)*(c-3)*v[[1]]==v[[3]],-v[[1]]/(v[[2]]+Log[2]v[[3]])}];c++];ClearAll["Global‘*"];
```

```
(* code snippet "12. Inostranstvo1" *)
```

```
NumQ[j_]:=7+6j+ContinuedFractionK[-2n(n+1)^2(n+2j+1),3n^2+(9+4j)n+(7+6j),{n,1,10000}];
D1:=Function[i,If[i<2,2+15i,2(2i-1)^3(1-i)*D1[i-2]+(8i^2-2i+3)D1[i-1]]];For[j=0,j<=14,r=N[NumQ[j],80];
v=FindIntegerNullVector[{1,r,N[Catalan*r,80]}];Q=-v[[1]]/(v[[2]]+Catalan v[[3]]);Print[{j,(2*j+1)!/(D1[j]-2 Catalan(2*j
    +1)!!^2)==Q}];j++];ClearAll["Global`*"];
```

```
(* code snippet "13. Inostranstvo2" *)
NumQ[j_]:=23+10j+ContinuedFractionK[-2n(n+1)(n+3)(n+2j+3),3n^2+(17+4j)n+(23+10j),{n,1,30000}];
For[j=0,j<=30,r=N[NumQ[j],180];v=FindIntegerNullVector[{1,r,N[Catalan*r,180]}];Print[{j,(2j+5)!/(2j+4)/(2j+5)!!^2==v
    [[1]]/v[[3]]}];j++];ClearAll["Global`*"];
```

```
(* code snippet "14. Inostranstvo" *)
NumQ:=Function[{i,a,b,c},z=a+b+1;x=a*b+(c+1)z;x+2z*i+ContinuedFractionK[-2n(n+a)(n+b)(n+2i+c),
3n^2+(2z+2c+1+4i)n+(x+2z*i),{n,1,5000}]];lim=6;ToV:=Function[{i,a,b,c},r=N[NumQ[i,a,b,c],80];v=FindIntegerNullVector
    [{1,r,N[Catalan*r,80]}];
Q=-v[[1]]/(v[[2]]+Catalan v[[3]])];Table[ToV[i,a+p,b+p,c+p],{p,0,1},{a,0,lim,2},{i,0,lim},{c,0,lim,2},{b,a,lim,2}]
```

```
(* code snippet "15. Series" *)
all={{Pi,{0,0,2,2,2,2,0},{57153600,-33075,103904},0},{Pi,{0,1,1,2,2,1,1},{69854400,-24255,76192},0},{Pi
    ,{0,1,2,1,1,2,1},{62868960,-14553,45712},0},{Pi,{0,2,2,1,2,2,0},{-65318400,-4725,14848},0},{Pi
    ,{0,2,2,2,2,0,0},{-6350400,-3675,11552},0},{Pi,{1,1,2,2,1,1,0},{-907200,-315,992},0},{Pi
    ,{1,2,1,1,2,1,0},{-635040,-147,464},0},{Pi,{1,2,2,2,2,1,0},{16934400,-245,768},0},{Pi
    ,{2,1,1,1,1,2,0},{-59535000,-6615,21292},0},{Pi,{2,2,0,0,2,1,1},{-93139200,-2695,9344},0},{Pi
    ,{2,2,0,0,2,2,0},{-52920000,-2695,9056},0},{Pi,{2,2,1,2,2,0,0},{-50803200,3675,-11264},0},{Pi
    ,{2,2,2,2,0,0,0},{129600,-75,224},0},{Catalan,{0,1,2,2,2,1,0},{-50803200,66150,-60577},0},{Catalan
    ,{0,3,2,3,0,0,0},{-3456000,-6750,6197},0},{Catalan,{1,2,2,2,1,0,0},{-2419200,-3150,2909},0},{Catalan
    ,{1,3,0,3,1,0,0},{8064000,8750,-8109},0},{Catalan,{2,2,0,2,2,0,0},{-33868800,-22050,21131},0},{Catalan
    ,{3,2,3,0,0,0,0},{-27648,54,-25},0},{Log[2],{0,0,0,2,2,2,2},{-33177600,-38400,26617},1},{Log
    [2],{0,0,2,2,2,2,0},{1382400,-1600,1109},1},{Log[2],{0,1,1,2,2,1,1},{11059200,-7680,5323},1},{Log
    [2],{0,1,2,1,1,2,1},{-22118400,10240,-7097},1},{Log[2],{0,2,2,0,0,2,2},{17280000,-1760,1219},1},{Log
    [2],{0,2,2,2,2,0,0},{442368,512,-355},1},{Log[2],{1,1,2,0,0,2,2},{-88473600,5120,-3539},1},{Log
```

```
    [2],{1,1,2,2,1,1,0},{-230400,-160,111},1},{Log[2],{1,2,1,1,2,1,0},{-22118400,-10240,7109},1},{Log
    [2],{2,2,0,0,2,1,1},{-88473600,-5120,3627},1},{Log[2],{2,2,0,0,2,2,0},{69120000,7040,-4951},1},{Log
    [2],{2,2,1,2,2,0,0},{7077888,-1024,707},1},{Log[2],{2,2,2,2,0,0,0},{-13824,16,-11},1}};
G:=Function[{e,ep},(-1)^(n+1)/Product[(2n+2i-3+ep)^e[[i]],{i,1,Length[e]}]];F:=Function[{e,ep},Sum[G[e,ep]/.{n->j},{j
    ,1,100000}]];Union[Table[{F[all[[i,2]],all[[i,4]]],all[[i,1]],1}.all[[i,3]]<10^30,{i,1,Length[all]}]]
```

```
(* code snippet "16. files" *)
NumQ[j_,k_,c_]:=j(2-j+2k)+ContinuedFractionK[-2n(c+n)(j+n-1)(1-j+2k+n),3n^2+(3+4k)n+j(2-j+2k),{n,1,17000}];
f[j_,k_,c_]:=Product[(2c-2k+2i-1)(k-i+1),{i,1,(j-1)/2}]CatalanNumber[(j-3)/2](j-2)(2k-1)*(2c-1)!!^2*CatalanNumber[k-1];
g[j_,k_,c_]:=Product[(2c-2i+1)(2k-2i+1),{i,1,(j-1)/2}](2c)!2^(2k+(j-7)/2);
h[j_,k_,c_]:=Product[2c-2i-1,{i,0,(j-3)/2}]*Product[2c-2i-1,{i,0,k-1}];
l:=Function[{n,j,k},(-1)^(k+1)(2k)!^2/k!/2^(3k-2)/Product[(k-i)(2k-2i-1)^2,{i,0,(j-3)/2}]/(2k-j)/(2k-1)/(n((2k-j-2)(3-2
    k)-1)+1)];
AB:=Function[{j,k},
For[c=1,c<=2,r=N[NumQ[j,k,c],5000];
v=FindIntegerNullVector[{1,r,N[Catalan*r,5000]}];
d[c]=(-g[j,k,c](v[[3]]Catalan+v[[2]])/v[[1]]-f[j,k,c]Catalan)/h[j,k,c];c++];{d[1],d[2]-d[1]}];
z[1]=z[0]={};
For[j=3,j<=501,For[w=0,w<=1,e[w]=Table[AB[j,k].{1,w}*l[w,j,k],{k,j-2,j-1}];
ab[w]={a,b}/.Solve[Table[a+b(u-1)==e[w][[u]],{u,1,2}],{a,b}];
AppendTo[z[w],ab[w][[1]]];w++];j+=2];
z[1]>>file1.txt;
z[0]>>file2.txt;
ClearAll["Global‘*"];
```

```
(*code snippet "17. j-level"*)
l[1]=<<file1.txt;l[2]=<<file2.txt;
data[j_,x_,y_]:=l[x][[(j-1)/2,y]];
rho[j_]:=2^j CatalanNumber[1/2(j-3)]((1/2(j-1))!)^2;
th1[j_]:=(j-6)(j-4)(j-2)(j-1)j(2j-7)(2j-5)(-1-3j+2j^2);
th2[j_]:=4(-390+312j+653j^2-942j^3+442j^4-87j^5+6j^6);
```

```
7  th3[j_]:=(j-6)(j-4)^2(j-1)(1+j)(13-11j+2j^2);
8  test1=Union[Table[data[j+2,1,2]==4(data[j,1,2]*th1[j]-th2[j]*rho[j])/th3[j],{j,3,499,2}]];
9  test2=Union[Table[data[j+2,1,1]==(4(data[j,1,1]j(j-6)(j-4)(j-2)(j-1)(2j-7)(2j-5)-12(2-4j+j^2)rho[j]))/((j-6)(j-4)^2(j
     ^2-1)),{j,3,499,2}]];
10 test3=Union[Table[data[j+2,2,2]==(4(data[j,2,2]j(2j+1)(j-1)(2j-5)(j-4)-(18+11j-17j^2+3j^3)rho[j]))/((j-4)(j-3)(1+j)),{j
     ,5,499,2}]];
11 test4=Union[Table[data[j+2,2,1]==(4(data[j,2,1](j-4)(j-1)j(2j-5)(2j-3)-(3j-2)rho[j]))/((j-4)(j-1)(1+j)),{j,3,499,2}]];
12 Print[Union[Join[test1,test2,test3,test4]]];
13 ClearAll["Global‘*"];
```

```
1  (*code snippet "18. Altogether"*)
2  NumQ[j_,k_,c_]:=j(2-j+2k)+ContinuedFractionK[-2n(c+n)(j+n-1)(1-j+2k+n),3n^2+(3+4k)n+j(2-j+2k),{n,1,5000}];
3  rho[j_]:=2^(j-3)(j-5)!(j-3)/(j-6);
4  p[j_,r_,i_]:={(j-3)(j-2)(2j-9)(2j-r),(j-1)(j-r),(j-8)(j-r),(2j-7)(j-r)}[[i]];
5  aa[j_]:=If[j==3,-1,(4(aa[j-2]p[j,7,1]-(3j-8)rho[j]))/p[j,3,2]];
6  ab[j_]:=If[j<=5,115j-341,4(ab[j-2]p[j,3,1]-(3j^3-35j^2+115j-96)rho[j])/p[j,5,2]];
7  ba[j_]:=If[j==3,-1/3,(4(ba[j-2]p[j,11,1]p[j,4,3]-12(p[j,7,2]+7)rho[j]))/(p[j,6,3]p[j,3,2])];
8  bb[j_]:=If[j==3,-14/3,4(bb[j-2]p[j,11,1]p[j,4,3](p[j,2,4]-1)-4(6x^6-15x^5-68x^4+74x^3+89x^2-44x-18/.x->j-4)rho[j])/(p[j
     ,3,2]p[j,6,3](p[j,6,4]+1))];
9  vD[z_,j_,k_,c_]:=If[c<2,z.{1,c},2c(j-2c)(2c-2k+j-2)(2c-2k-1)vD[z,j,k,c-2]+(8c^2+(2-8k)c+(j-2)(2k-j))vD[z,j,k,c-1]];
10 f[j_,k_,c_]:=Product[(2c-2k+2i-1)(k-i+1),{i,1,(j-1)/2}]CatalanNumber[(j-3)/2](j-2)(2k-1)(2c-1)!!^2*CatalanNumber[k-1];
11 g[j_,k_,c_]:=Product[2(2c-2i+1)(2k-2i+1),{i,1,(j-1)/2}](2c)!2^(2k-3);
12 h[j_,k_,c_]:=Product[2c-2i-1,{i,0,(j-3)/2}]Product[2c-2i-1,{i,0,k-1}];
13 l[n_,j_,k_]:=(2k)!^2/(k!(-2)^(3k-2)Product[(k-i)(2k-2i-1)^2,{i,0,(j-3)/2}](j-2k)(2k-1)(n((2k-j-2)(3-2k)-1)+1));
14 kD[n_,z_,j_,k_]:=If[k<2,z.{1,k},(2k+2j-9-2n)(2k+j-8-2n)(5-2k-j)(2k+j-6)kD[n,z,j,k-2]+(8k^2+k(10j-48-8n)+(3j^2-(28+4n)j
     +68+18n))kD[n,z,j,k-1]];
15 M:=Function[{j,k},t=Table[kD[u,{{aa[j],ab[j]},{ba[j],bb[j]}}[[u+1]],j,k-j+2]/l[u,j,k],{u,0,1}];{t[[1]],t.{-1, 1}}];
16 Q[j_,k_,c_,z_]:=Together[g[j,k,c]/(vD[z,j,k,c-1]h[j,k,c]+f[j,k,c]Catalan)];
17 Print[Union[Flatten[Table[N[{NumQ[j,k,c]==Q[j,k,c,M[j,k]]},200],{j,3,11,2},{k,j-2,10},{c,1,7}]]]];
18 ClearAll["Global‘*"];
```

```
(*code snippet "19. j-level bis"*)
{aa[3],ab[3],ba[3],bb[3],aa[5],ab[5],ba[5],bb[5],rho[j_],b[j_],a[j_],p[j_],d[j_],e[j_
    ]}:={-1,4,-1/3,-14/3,19,234,-17,-8,2^{j+1}(j-1)!/((j-2)(j-4)),(j-6)(j-2)(j-1)j(2j-7)(2j-5),4(j-1)j(2j-5),(j-6)(j
    -4)(j-1)(j+1)/4,6j^6-15j^5-68j^4+74j^3+89j^2-44j-18,(3j+1)(j^2-7)+3};
For[j=5,j<=17,Print[{aa[j+2],ab[j+2],ba[j+2],bb[j+2]}={(aa[j]a[j](2j-3)-(3j-2)rho[j])/((j-1)(j+1)),(ab[j]a[j](2j+1)-e[j
    -2]rho[j])/((j-3)(j+1)),(ba[j]b[j]-3((j-3)(j-1)-1)rho[j])/p[j],(bb[j]b[j](j(2j-3)-1)-d[j-2]rho[j])/(p[j]((j-2)(2j
    -7)-1))}];j+=2]
ClearAll["Global'*"];
```

```
(*code snippet "20. identities"*)
l[1][c_,k_,j_]:=Product[(2c-2k+2i-1)*(k-i+1),{i,1,(j-1)/2}];
r[1][c_,k_,j_]:=(2c-2k+j-2)!!/(2c-2k-1)!!*k!/(k-(j-1)/2)!;
l[2][c_,k_,j_]:=Product[(2c-2i+1)*(2k-2i+1),{i,1,(j-1)/2}];
r[2][c_,k_,j_]:=(2c-1)!!/(2c-j)!!(2k-1)!!/(2k-j)!!;
l[3][c_,k_,j_]:=Product[(2c-2i-1),{i,0,(j-3)/2}]*Product[(2c-2i-1),{i,0,k-1}];
r[3][c_,k_,j_]:=(2c-1)!!/(2c-j)!!(2c-1)!!/(2c-2k-1)!!;
l[4][c_,k_,j_]:=Product[(k-i)(2k-2i-1)^2,{i,0,(j-3)/2}];
r[4][c_,k_,j_]:=k!/(k-(j-1)/2)!*(2k-1)!!^2/(2k-j)!!^2;
t[c_,k_,j_]:=And@@Table[r[u][c,k,j]==l[u][c,k,j],{u,1,4}];
And@@Flatten[Table[t[c,k,j],{c,1,10},{k,1,10},{j,3,19,2}]]
ClearAll["Global'*"];
```

```
(*code snippet "21. compact representation using snippet 20"*)
(*caution: functions changed w.r.t. paper and previous code.*)
NumQ[j_,k_,c_]:=j(2-j+2k)+ContinuedFractionK[-2n(c+n)(j+n-1)(1-j+2k+n),3n^2+(3+4k)n+j(2-j+2k),{n,1,5000}];
rho[j_]:=2^(j-1)(j-5)!(j-3)/(j-6);
p[r_,i_]:={4(y-3)(y-2)(2y-9)(2y-r),(y-1)(y-r),(y-8)(y-r),(2y-7)(y-r),(v-r)(v-4)(v-5)}[[i]];
aa[j_]:=If[j==3,-1,(aa[j-2]p[7,1]-(3j-8)rho[j])/p[3,2]/.{y->j}];
ab[j_]:=If[j<=5,115j-341,(ab[j-2]p[3,1]-(3j^3-35j^2+115j-96)rho[j])/p[5,2]/.{y->j}];
ba[j_]:=If[j==3,-1/3,(ba[j-2]p[11,1]p[4,3]-12(p[7,2]+7)rho[j])/(p[6,3]p[3,2])/.{y->j}];
```

```
9  bb[j_]:=If[j==3,-14/3,(bb[j-2]p[11,1]p[4,3](p[2,4]-1)-(6x^6-15x^5-68x^4+74x^3+89x^2-44x-18/.x->j-4)4rho[j])/(p[3,2]p
       [6,3](p[6,4]+1))/.{y->j}];
10 vD[d_,c_]:=({z,j,k}=d;If[c<2,z[[c+1]],2c(j-2c)(2c-2k+j-2)(2c-2k-1)vD[d,c-2]+(8c^2+(2-8k)c+(j-2)(2k-j))vD[d,c-1]]);
11 f[j_,k_,c_]:=2(2c-j)!!(2c-2k+j-2)!!(2k-1)!(j-2)!/((k-1)!((j-3)/2)!^2(j-1)(k-(j-1)/2)!);
12 g[j_,k_,c_]:=(2k-1)!2^(k+(j-5)/2+c)c!(2c-2k-1)!!/(k-1)!/(2k-j)!!;
13 l[n_,j_,k_]:=(-2)^(k-2)(j-2k)(2k-1)(2k-j-2)^n(3-2k)^n/(k-(j-1)/2)!/(2k-j)!!^2;
14 kD[d_,k_]:=({n,z,j}=d;If[k<2,z.{1,k},p[n,5](n-2k)kD[d,k-2]+((v-3)^2+v(2k-2n+1)+4j+7n-19)kD[d,k-1]/.{v->2k+j-1}]);
15 M[j_,k_]:=Table[kD[{2u+8-j,{{aa[j],ab[j]},{ba[j],bb[j]}}[[u+1]],j},k-j+2]l[u,j,k],{u,0,1}];
16 Q[j_,k_,c_]:=g[j,k,c]/(vD[{M[j,k],j,k},c-1]+f[j,k,c]Catalan);
17 Print[Union[Flatten[Table[N[{NumQ[j,k,c]==Q[j,k,c]},200],{j,3,19,2},{k,j-2,15},{c,1,12}]]]];
18 ClearAll["Global‘*"];
```

## C Altogether

This appendix recaps the entire algorithm. It takes as input $j, \kappa, c$ and returns the exact expression of $Q_{j,\kappa,c}$ where:

$$Q_{j,\kappa,c} = j(2-j+2\kappa) + \operatorname*{K}_{n=1}^{\infty}\left(\frac{-2n(c+n)(j+n-1)(1-j+2\kappa+n)}{j(2-j+2\kappa)+(3+4\kappa)n+3n^2}\right)$$

### C.1 If $j \geq 2\kappa + 3$

In this case $(j, \kappa, c) \in$ Bosnia & Herzegovina $\cup$ Croatia, hence $Q_{j,\kappa,c}$ is computed by straightforward finite summation.

$$Q_{j,\kappa,c} = j(2-j+2\kappa) + \operatorname*{K}_{n=1}^{j-2\kappa-1}\left(\frac{-2n(c+n)(j+n-1)(1-j+2\kappa+n)}{j(2-j+2\kappa)+(3+4\kappa)n+3n^2}\right)$$

### C.2 If $3 + \kappa \leq j \leq 2\kappa + 1$

In this case we are in Serbia. We hence use the symmetry relation:

$$Q_{j,\kappa,c} = Q_{2(\kappa+1)-j,\kappa,c}$$

Indeed,

$$3+\kappa \leq j \leq 2\kappa+1 \Rightarrow 1 \leq j' \leq \kappa - 1 \leq \kappa + 2 \Rightarrow (j', \kappa) \in \text{Montenegro} \cup \text{Kosovo}$$

Replace $j$ by $j' = 2(\kappa+1) - j$ and compute $Q_{j',\kappa,c}$ using the formula for Montenegro or for Kosovo given in the next subsections.

### C.3 If $j = 1$

In this case we are in Montenegro. We hence define:

$$\Delta_{\kappa,c}(\alpha,\beta) = \begin{cases} \alpha + \beta c & \text{if } c < 2 \\ -2c(2c-1)(2(c-\kappa)-1)^2\Delta_{\kappa,c-2}(\alpha,\beta) & \text{if } c \geq 2 \\ \quad + (8c^2 + (2-8\kappa)c - 2\kappa + 1)\Delta_{\kappa,c-1}(\alpha,\beta) & \end{cases}$$

$$\Gamma_{\kappa,c}(\alpha,\beta) = (2c-1)!!^2 G + \Delta_{\kappa,c-1}(\alpha,\beta) \cdot \prod_{i=0}^{\kappa-1}(2(c-i)-1)$$

$$\delta_\kappa = \frac{4^{\kappa-1}}{(2\kappa-1)C_{\kappa-1}} \quad \text{and} \quad \rho_\kappa = \frac{\delta_\kappa(-1)^\kappa(1-2\kappa)}{(2\kappa)!(2\kappa-3)!!}$$

$$\alpha_\kappa = \rho_\kappa \Delta_{1,\kappa-1}(1,-2) \quad \text{and} \quad \beta_\kappa = -\rho_\kappa(2\kappa-3)^2\Delta_{2,\kappa-1}(1,12) - \alpha_\kappa$$

$$\text{And output} \quad Q_{1,\kappa,c} = \frac{\delta_\kappa(2c)!}{\Gamma_{\kappa,c}(\alpha_\kappa,\beta_\kappa)}$$

### C.4 If $3 \leq j \leq \kappa + 2$

$3 \leq j \leq \kappa + 2$ is in Kosovo. In which case we proceed in three steps:

**Step 1:**
If $j = 3$ or $j = 5$ let:

$$\{\overset{\alpha}{\alpha}_3, \overset{\beta}{\alpha}_3, \overset{\alpha}{\beta}_3, \overset{\beta}{\beta}_3\} = \{-1, 4, -\frac{1}{3}, -\frac{14}{3}\} \text{ and } \{\overset{\alpha}{\alpha}_5, \overset{\beta}{\alpha}_5, \overset{\alpha}{\beta}_5, \overset{\beta}{\beta}_5\} = \{19, 234, -17, -8\}$$

If $j \geq 7$ define:

$$\varrho_j = \frac{2^{j+1}(j-1)!}{(j-2)(j-4)} \text{ and } b_j = (j-6)(j-2)(j-1)j(2j-7)(2j-5)$$

$$a_j = 4(j-1)j(2j-5) \text{ and } p_j = \frac{(j-6)(j-4)(j-1)(j+1)}{4}$$

$$d_j = 6j^6 - 15j^5 - 68j^4 + 74j^3 + 89j^2 - 44j - 18 \text{ and } e_j = (3j+1)(j^2-7) + 3$$

and iterate using the following formulae to compute $\{\overset{\alpha}{\alpha}_j, \overset{\beta}{\alpha}_j, \overset{\alpha}{\beta}_j, \overset{\beta}{\beta}_j\}$:

$$\overset{\alpha}{\alpha}_{j+2} = \frac{\overset{\alpha}{\alpha}_j a_j (2j-3) - (3j-2)\varrho_j}{(j-1)(j+1)} \text{ and } \overset{\beta}{\alpha}_{j+2} = \frac{\overset{\beta}{\alpha}_j a_j (2j+1) - e_{j-2}\varrho_j}{(j-3)(j+1)}$$

$$\overset{\alpha}{\beta}_{j+2} = \frac{\overset{\alpha}{\beta}_j b_j - 3\left((j-3)(j-1)-1\right)\varrho_j}{p_j} \text{ and } \overset{\beta}{\beta}_{j+2} = \frac{\overset{\beta}{\beta}_j b_j (j(2j-3)-1) - d_{j-2}\varrho_j}{p_j((j-2)(2j-7)-1)}$$

**Step 2:** Define:

$$s(j,\kappa) = \prod_{i=0}^{\frac{j-3}{2}} (\kappa - i)(2\kappa - 2i - 1)^2$$

$$\ell(n,j,\kappa) = \frac{(-1)^{\kappa+1}(2\kappa)!^2}{\kappa! 2^{3\kappa-2}(2\kappa-j)(2\kappa-1)(n((2\kappa-j-2)(3-2\kappa)-1)+1) \cdot s(j,\kappa)}$$

$$\eta(n,j,\kappa) = (2\kappa+2j-9-2n)(2\kappa+j-8-2n)(-2\kappa+5-j)(2\kappa+j-6)$$

$$\phi(n,j,\kappa) = 8\kappa^2 + \kappa(10j - 48 - 8n) + 3j^2 - (28+4n)j + 68 + 18n$$

$$\bar{\Delta}_{n,j,\kappa}(\alpha,\beta) = \begin{cases} \alpha + \beta\kappa & \text{if } \kappa < 2 \\ \eta(n,j,\kappa) \cdot \bar{\Delta}_{n,j,\kappa-2}(\alpha,\beta) + \phi(n,j,\kappa) \cdot \bar{\Delta}_{n,j,\kappa-1}(\alpha,\beta) & \text{if } \kappa \geq 2 \end{cases}$$

Using the values $\overset{\alpha}{\alpha}_j, \overset{\beta}{\alpha}_j, \overset{\alpha}{\beta}_j, \overset{\beta}{\beta}_j$ compute:

$$\alpha_{j,\kappa} = \frac{\bar{\Delta}_{0,j,\kappa-j+2}(\overset{\alpha}{\alpha}_j, \overset{\beta}{\alpha}_j)}{\ell(0,j,\kappa)} \quad \text{and} \quad \beta_{j,\kappa} = \frac{\bar{\Delta}_{1,j,\kappa-j+2}(\overset{\alpha}{\beta}_j, \overset{\beta}{\beta}_j)}{\ell(1,j,\kappa)} - \alpha_{j,\kappa}$$

**Step 3:** Define:

$$\Delta_{j,\kappa,c} = \begin{cases} \alpha_{j,\kappa} + \beta_{j,\kappa} c & \text{if } c < 2 \\ -2c(2c-j)(2c-2\kappa+j-2)(2c-2\kappa-1)\Delta_{j,\kappa,c-2} & \text{if } c \geq 2 \\ \quad + (8c^2 + (2-8\kappa)c + (j-2)(2\kappa-j))\Delta_{j,\kappa,c-1} & \end{cases}$$

$$f_{j,\kappa,c} = C_{\frac{j-3}{2}} C_{\kappa-1}(j-2)(2\kappa-1)(2c-1)!!^2 \prod_{i=1}^{\frac{j-1}{2}} (2c-2\kappa+2i-1)(\kappa-i+1)$$

$$g_{j,\kappa,c} = (2c)! 2^{\frac{j+4\kappa-7}{2}} \prod_{i=1}^{\frac{j-1}{2}} (2c-2i+1)(2\kappa-2i+1)$$

$$h_{j,\kappa,c} = \prod_{i=0}^{\frac{j-3}{2}} (2c-2i-1) \prod_{i=0}^{\kappa-1} (2c-2i-1)$$

Output:

$$Q_{j,\kappa,c} = \frac{g(j,\kappa,c)}{\Delta_{j,\kappa,c-1} \cdot h(j,\kappa,c) + f(j,\kappa,c) \cdot G}$$

Because the nontrivial part of the algorithm is Kosovo we only provide a code snippet (`"18. Altogether"`) testing the final formulae over Kosovo. It seems amazing that such as short yet complex code perfectly captures exactly the behavior of $Q_{j,\kappa,c}$.

## D Avoiding $x!!$

Note that:

$$\prod_{i=1}^{\frac{j-1}{2}}(2c-2\kappa+2i-1)\cdot(\kappa-i+1) = \frac{(2c-2\kappa+j-2)!!}{(2c-2\kappa-1)!!}\cdot\frac{\kappa!}{(\kappa-\frac{j-1}{2})!}$$

$$\prod_{i=1}^{\frac{j-1}{2}}(2c-2i+1)\cdot(2\kappa-2i+1) = \frac{(2c-1)!!}{(2c-j)!!}\cdot\frac{(2\kappa-1)!!}{(2\kappa-j)!!}$$

$$\prod_{i=0}^{\frac{j-3}{2}}(2c-2i-1)\prod_{i=0}^{\kappa-1}(2c-2i-1) = \frac{(2c-1)!!}{(2c-j)!!}\cdot\frac{(2c-1)!!}{(2c-2\kappa-1)!!}$$

$$\prod_{i=0}^{\frac{j-3}{2}}(\kappa-i)(2\kappa-2i-1)^2 = \frac{k!}{(\kappa-\frac{j-1}{2})!}\cdot\frac{(2\kappa-1)!!^2}{(2\kappa-j)!!^2}$$

And the well-known identities:

$$(2x-1)!! = \frac{(2x)!}{2^x x!} = \frac{(2x-1)!}{2^{x-1}(x-1)!} \quad\text{and}\quad C_x = \frac{(2x)!}{(x+1)!\,x!}$$

yield expressions that avoid double factorials. The first identity is always usable because $j$ is odd.

## E Lazy evaluation

A very useful trick allowing to speed-up the computation of $Q_{j,\kappa,c}$ is lazy evaluation. We illustrate it in the case of Kosovo.

$g, h, f$ are all smooth. Hence it is inefficient to compute them and then simplify the resulting fraction. Instead, it is best to keep for each of those functions a vector whose coordinates count the occurrence of small prime factors met during calculation – without actually multiplying integers.

When the three vectors $V_g, V_h, V_f$ are ready we can identify at each coordinate $i$ the value:

$$t_i = \min(V_g[i], V_h[i], V_f[i])$$

and update: $V_\bullet[i] = V_\bullet[i] - t_i$ before performing the actual multiplication.

Estimating the largest primes that potentially appearar in each of the functions is easy. $f$'s largest prime is $\max(2\kappa-1, 2c+j-2)$, $g$'s is $\max(2\kappa-1, 2c)$ and $h$'s is $2c-1$. Hence it suffices to consider primes smaller than $\max(2\kappa-1, 2c+j-2)$.

Given that the coefficients of $Q_{j,\kappa,c}$ would explode for $j, \kappa, c$ values exceeding 1000 it suffices to tabulate the factorization of integers between 2 and 2998 and use this table to distribute factors in the vectors.

## F Monetary prize

The authors offer four cash prizes and signed acknowledgment certificates to the first researchers who will prove or refute the following.

| Conjecture | Assume steps | Prove steps | Prize |
|---|---|---|---|
| **C.3** | | **C.3** | €100 |
| **C.4** | 1,2 | 3 | €100 |
| **C.4** | 1 | 2,3 | €100 |
| **C.4** | | 1,2,3 | €100 |

**Table 15.** Prizes

To claim the prizes, provers should publish their proofs in a peer-reviewed venue or journal.